\documentclass[11pt,reqno]{amsart}
\usepackage{cases}
\usepackage{amscd}
\usepackage{amsfonts}
\usepackage{amssymb}
\usepackage{amsmath}
\usepackage{graphicx}
\usepackage{epstopdf}
\usepackage{subfigure}
\usepackage{amsthm}
\usepackage{algorithm} 
\usepackage{algorithmic} 
\usepackage{multirow} 
\usepackage{stmaryrd}
\usepackage{caption} 
\theoremstyle{remark}
\newtheorem{example}{\textbf{Example}}[section]
\numberwithin{equation}{section}

\usepackage{color}
\usepackage{datetime}
\usepackage{tikz}
\usepackage{hyperref}
\usepackage{cleveref}

\makeatletter
\newcommand\figcaption{\def\@captype{figure}\caption}
\newcommand\tabcaption{\def\@captype{table}\caption}
\makeatother

\makeatletter
\newenvironment{breakablealgorithm}
  {
   \begin{center}
     \refstepcounter{algorithm}
     \hrule height.8pt depth0pt \kern2pt
     \renewcommand{\caption}[2][\relax]{
       {\raggedright\textbf{\ALG@name~\thealgorithm} ##2\par}%
       \ifx\relax##1\relax 
         \addcontentsline{loa}{algorithm}{\protect\numberline{\thealgorithm}##2}%
       \else 
         \addcontentsline{loa}{algorithm}{\protect\numberline{\thealgorithm}##1}%
       \fi
       \kern2pt\hrule\kern2pt
     }
  }{
     \kern2pt\hrule\relax
   \end{center}
  }
\makeatother
\usepackage{geometry}
\usepackage{algorithm}
\usepackage{multirow}

\def\bq{\begin{equation}}
\def\eq{\end{equation}}
\def\bqs{\begin{equation*}}
\def\eqs{\end{equation*}}
\def\bsqs{\begin{subequations}}
\def\esqs{\end{subequations}}
\def\ba{\begin{aligned}}
\def\ea{\end{aligned}}
\def\br{\begin{eqnarray}}
\def\er{\end{eqnarray}}
\def\brr{\bq\begin{array}{rlll}}
\def\err{\end{array}\eq}

\def\text#1{\hbox{#1}}
\newtheorem{thm}{Theorem}[section]
\newtheorem{lem}[thm]{Lemma}

\newtheorem{rem}[thm]{Remark}

\newtheorem{definition}{Definition}[section]

\newcommand{\bsub}{\begin{subequations}}
\newcommand{\esub}{\end{subequations}$\!$}

\makeatletter
\newcommand{\definetitlefootnote}[1]{%
  \newcommand\addtitlefootnote{%
    \makebox[0pt][l]{$^{*}$}%
    \footnote{\protect\@titlefootnotetext}
  }%
  \newcommand\@titlefootnotetext{\spaceskip=\z@skip $^{*}$#1}%
}
\makeatother

\title[Estimators for $C^0$ interior penalty]{A posteriori error estimators for fourth order elliptic problems with concentrated loads}

\author[H.~Cao, Y.~Huang, N.~Yi, P. Yin]{Huihui Cao$^{\dag}$, Yunqing Huang$^\dag$, Nianyu Yi$^\dag$, Peimeng Yin$^{\ddag,*}$}

\address{$\dag$ Hunan Key Laboratory for Computation and Simulation in Science and Engineering, School of Mathematics and Computational Science, Xiangtan University, Xiangtan 411105, Hunan, P.R.China} \email{caohh@sustech.edu.cn (H. Cao);\  huangyq@xtu.edu.cn (Y. Huang);\  yinianyu@xtu.edu.cn (N. Yi).}
\address{$^*$ Department of Mathematical Sciences, University of Texas at El Paso,  El Paso, Texas 79968, USA.}\email{pyin@utep.edu}


\thanks{$^*$ Corresponding author.}

\begin{document}

\begin{abstract} 
In this paper, we study two residual-based a posteriori error estimators for the $C^0$ interior penalty method in solving the biharmonic equation in a polygonal domain under a concentrated load. The first estimator is derived directly from the model equation without any post-processing technique. We rigorously prove the efficiency and reliability of the estimator by constructing bubble functions. Additionally, we extend this type of estimator to general fourth-order elliptic equations with various boundary conditions. The second estimator is based on projecting the Dirac delta function onto the discrete finite element space, allowing the application of a standard estimator. Notably, we additionally incorporate the projection error into the standard estimator. The efficiency and reliability of the estimator are also verified through rigorous analysis. We validate the performance of these a posteriori estimates within an adaptive algorithm and demonstrate their robustness and expected accuracy through extensive numerical examples.
\end{abstract}

\maketitle

\bigskip



\section{Introduction}

In this paper, we are interested in an adaptive $C^0$ interior penalty method for the biharmonic problem in  polygonal domain $\Omega \subset \mathbb{R}^2$ with a concentrated load \cite{TSW1959, Camp87}
\begin{equation}\label{eq:bh1}
\Delta^2 u  = \delta_{\mathbf{x}_0}  \quad \text{in }  \Omega, \quad
u =0 \quad \text{and} \quad \partial_\mathbf{n} u  =0 \quad \text{on }  \partial \Omega.
\end{equation}
The boundary conditions are known as homogeneous Dirichlet boundary conditions or clamped boundary conditions \cite{KN2014}, where $\partial_\mathbf{n} u$ denotes the outward normal derivative of $u$ on $\partial \Omega$. 
$\delta_{\mathbf{x}_0}$ is a Dirac delta function concentrated at a point $\mathbf{x}_0 \in \Omega_0 \subset \Omega$ satisfying
\[
\langle \delta_{\mathbf{x}_0}, v \rangle =  v(\mathbf{x}_0), \qquad \forall\ v \in C(\Omega_0).
\]

Elliptic problems with Dirac delta source terms are encountered in various applications, such as the electric field generated by a point charge, transport equations for effluent discharge in aquatic media, modeling of acoustic monopoles {\cite{J1975,ABR2007,LV2013,GHZ2014}. 
The biharmonic problem can be used to study the small deflections of a thin plate, especially the biharmonic problem \eqref{eq:bh1} with the Dirac delta source term describes the deflections for thin plates with a concentrated load \cite{TSW1959, Camp87}.


The Dirac measure in \eqref{eq:bh1} does not belong to $H^{-1}(\Omega)$, resulting in the solution exhibiting low regularity.
Analytical solutions for biharmonic problem \eqref{eq:bh1} are typically challenging, though they exist for some special cases of geometry and loads. 
For example, analytical methods for the biharmonic problem \eqref{eq:bh1} have primarily focused on circular and annular domains (see, e.g., \cite{TSW1959, Camp87, CL2011}). Consequently, numerical methods have garnered widespread attention for solving \eqref{eq:bh1}, such as the boundary element method \cite{Camp87}, the regular hybrid boundary node method \cite{TZ2013}. Discussion on various numerical methods for biharmonic problem \eqref{eq:bh1} can be found in \cite{TZ2013} and the references therein. Among various numerical methods, the finite element method is the most popular.

Finite element methods developed for the general biharmonic equation with Dirichlet boundary conditions can typically be applied to solve the specific biharmonic problem \eqref{eq:bh1}.
The conforming finite element method is one approach, where the presence of high-order derivatives necessitates finite element spaces that belong to $H^2$, such as the $C^1$ Argyris finite element method \cite{A1968}. Additionally, the general biharmonic problem can be decomposed into Poisson and Stokes equations, which are then solved using the $C^0$ finite element method \cite{LWY22}. However, whether this decomposition strategy is effective for problem \eqref{eq:bh1} remains to be explored since such decomposition requires the source term in $H^{-1}(\Omega)$ or its subset.
Another option that utilizes the $C^0$ finite element space, yet can still accommodate singular source terms not in $H^{-1}(\Omega)$, is the $C^0$ interior penalty method \cite{EGHLMT2002}. Its stability is ensured by penalty terms enforced across the mesh cell interfaces.

For the biharmonic problem \eqref{eq:bh1}, some finite element methods and error analyses are available in the literature. 
A $C^1$ finite element approximation was proposed in \cite{Scott73}, and optimal error estimates were studied on quasi-uniform meshes, in which the $H^2$ error estimate is of order $h$ when using polynomials of degree greater than $2$.
More recently, a $C^0$ interior penalty method was studied in \cite{L2020}, and a local $H^2$ error estimate of order $|\ln h|^{\frac{3}{2}}$ was given on quasi-uniform meshes. 
Due to the low regularity of the solution, the convergence rates on quasi-uniform meshes are inherently limited. To improve the convergence rate in the finite element approximation, adopting an adaptive finite element method becomes necessary. Thus, the primary objective of this paper is to develop an adaptive $C^0$ interior penalty method for the biharmonic problem \eqref{eq:bh1}.

Many adaptive finite element methods are available for second order elliptic equations with Dirac delta source term. Tough $\delta_{\mathbf{x}_0}$ is not in $H^{-1}(\Omega)$ and $u \notin H^1(\Omega)$, the $C^0$ finite element method can still approximate the equations. However, direct application of the residual-based a posteriori error estimator using standard energy norms is not viable. In the literature, two typical strategies have been studied to address this issue. One approach is to utilize norms weaker than $H^1(\Omega)$ for error estimation. For instance, Araya et al. \cite{ABR2006} derived a posteriori error estimators in $L^p(1<p<\infty)$ norm and $W^{1,p}(p_0<p<2,p_0\in [1,2))$ seminorms for a Poisson problem with a Dirac delta source term on two-dimensional domains.  Gaspoz et al. \cite{GMV2016} provided a posteriori error estimates in $H^{1 - s},s \in (0,\frac{1}{2})$ norm. Additionally, a global upper bound and a local lower bound of residual type a posteriori error estimators in a weighted Sobolev norm $\|\cdot\|_{H_{\alpha}^1}$ with $\alpha \in (\frac{d}{2} -1,\frac{d}{2})$ ($d$ is the spatial dimension) for elliptic problems were obtained by Agmon et al. in \cite{AGM2014}. Another is to regularize the source term to an $L^2(\Omega)$ function by projecting it onto a polynomial space, potentially introducing a projection error. This regularization allows the application of standard residual-based a posteriori error estimators for general Poisson problems. For further insights, readers are referred to early review articles, such as \cite{HW2012,MMRZ2022}.

Results on a posteriori error estimates of the $C^0$ interior penalty method for biharmonic problems with $L^2(\Omega)$ source terms can be found in \cite{GPJ2009}.
However, no result is available for the fourth order elliptic equation with the Dirac delta source term, which does not belong to $H^{-1}(\Omega)$, not to mention in $L^2(\Omega)$. Therefore, to improve the accuracy of the numerical solution while optimizing the distribution of computational resources, we propose two types of residual-based a posteriori error estimators for the biharmonic problem \eqref{eq:bh1} to guide mesh adaptive refinement around singular points.

The first type of a posterior error estimator is derived based on the primal equation \eqref{eq:bh1}. Depending on the location of $\mathbf{x}_0$ in the computational element, this error estimator can take different forms. 
Specifically, if $\mathbf{x}_0$ is not a vertex, an additional term that depends on the size of the element will be required.
We rigorously prove the upper and lower bounds of the proposed estimator to ensure its reliability and efficiency. Moreover, we extend this residual-type a posteriori error estimator to fourth-order elliptic equations with various boundary conditions.

The second type of a posterior error estimator is proposed based on the projection techniques. We first project $\delta_{\mathbf{x}_0}$ onto $\delta_h$ in the finite element space, and then use this projection $\delta_h$ to construct a residual-type a posteriori error estimator. This method introduces an additional error between $\delta_h$ and $\delta_{\mathbf{x}_0}$, which is shown to be of the same order as the finite element approximation. Therefore, it does not compromise the accuracy of the numerical solution. This is further supported by our error analysis and numerical experimental results.

The rest of the paper is organized as follows. In Section \ref{sec2}, we establish the well-posedness and discrete problem of \eqref{eq:bh1} by the $C^0$ interior penalty method. The main results are presented in Section \ref{sec3}, where we propose two types of residual-based a posteriori error estimators, upper and lower bounds are proved in order to guarantee the reliability and the efficiency of the proposed estimators. In Section \ref{sec4}, we extend our results to a broader class of fourth-order elliptic equations with various boundary conditions. Section \ref{sec5} provides numerous numerical examples to illustrate the robustness of our estimators and the corresponding adaptive $C^0$ interior penalty method. Finally, we draw some conclusions in Section \ref{sec6}.

Throughout the paper, the generic constant $C>0$ in our estimates may differ at different occurrences. It will depend on the computational domain, but not on the functions involved or on the mesh level in the finite element algorithms.

\section{Preliminaries and $C^0$ interior penalty method}\label{sec2}

Denote by $H^m(\Omega)$, $m$ is a non-negative integer,  the Sobolev space that consists of functions whose $i$th ($0\leq i\leq m$) derivatives are square integrable. Let $L^2(\Omega):=H^0(\Omega)$.
Denote by $H^1_0(\Omega)\subset H^1(\Omega)$  the subspace consisting of functions with  zero trace on the boundary  $\partial\Omega$.
For $0<t<1$ and $s=m+t$, the fractional order Sobolev space $H^s(D)$ consists of distributions $v\in D \subset \mathbb{R}^d$ ($d=1,2$) satisfying
$$
\|v\|^2_{H^s(D)}:=\|v\|^2_{H^m(D)} + \sum_{|\nu|= m}\int_{D}\int_{D} \frac{|\partial^\nu v(x) - \partial^\nu v(y)|^2 }{|x-y|^{d+2t}} dxdy <\infty,
$$
where $\nu = (\nu_1, \cdots, \nu_d) \in \mathbb{Z}^d_{\geq 0}$ is a multi-index such that $\partial^\nu=\partial_{x_1}^{\nu_1}\cdots\partial^{\nu_d}_{x_d}$ and $|\nu|=\sum_{i=1}^d\nu_i$.

\subsection{Well-posedness and regularity}

We first show the well-posedness of the problem \eqref{eq:bh1}. 
\begin{lem}\label{lemma2-2}
For any $\epsilon>0$, it follows that the point Dirac delta function $\delta_{\mathbf{x}_0} \in H^{-1-\epsilon}(\Omega)$ and satisfies
$$
\|\delta_{\mathbf{x}_0}\|_{H^{-1-\epsilon}(\Omega)} \leq C.
$$
\end{lem}
\begin{proof}
For any $v\in H^{1+\epsilon}(\Omega)$, the embedding theorem \cite{Ciarlet74} implies that $v \in C^{0,\epsilon}(\Omega) \subset C^0(\Omega)$. Then, 
$$
|\langle \delta_{\mathbf{x}_0}, v\rangle| = |v(\mathbf{x}_0)|\leq \|v\|_{L^\infty(\Omega)} \leq C\|v\|_{H^{1+\epsilon}(\Omega)},
$$
and
$$
\|\delta_{\mathbf{x}_0}\|_{H^{-1-\epsilon}(\Omega)} := \sup_{v\in H^{1+\epsilon}(\Omega)\backslash{\{0\}}}\frac{|\langle \delta_{\mathbf{x}_0}, v\rangle|}{\|v\|_{H^{1+\epsilon}(\Omega)}} \leq C.
$$
\end{proof}

The variational formulation for problem (\ref{eq:bh1}) is to find $u\in H_0^2(\Omega)$, such that
\begin{eqnarray}\label{eqn.weak}
a(u, v):=\int_\Omega\Delta u \Delta v \,\mathrm{d}\mathbf{x}=\left\langle \delta_{\mathbf{x}_0}, v \right\rangle, \quad \forall\ v\in H_0^2(\Omega).
\end{eqnarray}
The Sobolev imbedding theorem \cite{M20} implies $v \in C(\Omega)$ for $v \in H_0^2(\Omega)$, thus the variational formulation (\ref{eqn.weak}) is well-posed.

We sketch a drawing of the domain $\Omega$ with a singular point $\mathbf{x}_0$ in Figure \ref{fig:Omega2}. We assume the largest interior angle $\omega \in [\frac{\pi}{3}, 2\pi)$ of the domain associated with the vertex $Q$. 
For simplicity of the analysis, we assume that 
\bq\label{wcond}
\sin\sqrt{\frac{\omega^2}{\sin^2 \omega}-1} \not = \sqrt{1-\frac{\sin^2 \omega}{\omega^2}}.
\eq

\begin{figure}[H]
\centering
\begin{tikzpicture}[scale=0.25]
\draw[thick]
(-6,-11) -- (2,-11) -- (0,-2) -- (10,-2) -- (8,7) -- (-8,6) -- (-11,-2) -- (-6,-11);
\draw (7,5.5) node {$\Omega$};
\draw (1.4,-1.3) node {$\omega$};
\draw (1.1,-2.8) node {$Q$};

\draw[ultra thick]  (-5,1) node {$\bullet$};


\draw (-6,1.8) node {$\mathbf{x}_0$};
\draw[thick][densely dotted] (0.5,-2) arc (0:283:0.5);
\end{tikzpicture}
\vspace*{-10pt}
    \caption{Domain $\Omega$ with interior angle $\omega$ contains  a singular point $\mathbf{x}_0$.}
    \label{fig:Omega2}
\end{figure}
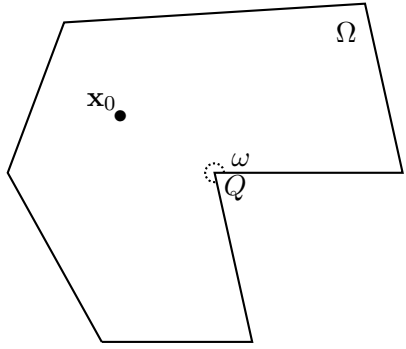

Let $z_\ell$, $\ell=1,2,\ldots$ satisfying $\rm{Re}(z_\ell)>0$ be the solutions of the following characteristic equation
\bq\label{fourthchar}
\sin^2(z\omega) = z^2\sin^2(\omega).
\eq
Then there exists a threshold
\bq\label{alpha0}
\alpha_0 := \min\{\rm{Re}(z_\ell), \  \ell=1,2,\ldots, \} > \frac{1}{2},
\eq
such that the following regularity result holds.
\begin{lem}\label{udecompthm}
For any $\epsilon>0$, let $ u$ be the solution of the biharmonic problem \eqref{eq:bh1}. Then it follows $u \in H^{\min\{3- \epsilon, 2+\alpha\}}(\Omega) \cap H^2_0(\Omega)$ with $\frac{1}{2}< \alpha <\alpha_0$. Moreover, if $\omega<\pi$, it holds $u \in H^{3- \epsilon}(\Omega) \cap H^2_0(\Omega)$; and if $\omega>\pi$, it holds $u \in H^{2+\alpha}(\Omega) \cap H^2_0(\Omega)$.
\end{lem}
The graph of $\alpha_0$ in terms of the largest interior angle $\omega$ is shown in \Cref{Regularity}, and some numerical values of $\alpha_0$ are shown in Table \ref{alpha0tab} \cite{LWY22}.
In \Cref{udecompthm}, when $\omega<\pi$, the regularity is dominated by the singularity of the Dirac delta source; and when $\omega>\pi$, the regularity is dominated by the singularity of the domain \cite{kozlov2001, bacuta2002, bourlard1992, Grisvard92}. To design high-order accurate numerical methods, one has to handle the singularities introduced by these two singular sources: domain corner, and Dirac delta source.

 \begin{table}[!htbp]\tabcolsep0.04in
 \caption{Values of $\alpha_0$ for different interior angles $\omega$.}
 \begin{tabular}[c]{|c|c|c|c|c|c|c||c|c|c|c|c|c|}
 \hline
 $\omega$ & $\frac{\pi}{3}$ & $\frac{\pi}{2}$ & $\frac{2\pi}{3}$ & $\frac{3\pi}{4}$ & $\frac{5\pi}{6}$ & $\frac{11\pi}{12}$ & $\frac{7\pi}{6}$ & $\frac{6\pi}{5}$ & $\frac{5\pi}{4}$ & $\frac{4\pi}{3}$ & $\frac{3\pi}{2}$ & $\frac{7\pi}{4}$   \\
 \hline
 $\alpha_0 \approx$  &  4.0593 & 2.7396 & 2.0941 & 1.8854 & 1.5339 & 1.2006 & 0.7520 & 0.7178 & 0.6736 & 0.6157 & 0.5445 & 0.5050 \\
 \hline
 \end{tabular}\label{alpha0tab}
 \end{table}

 \begin{figure}
 \centering
 \subfigure{\includegraphics[width=0.49\textwidth]{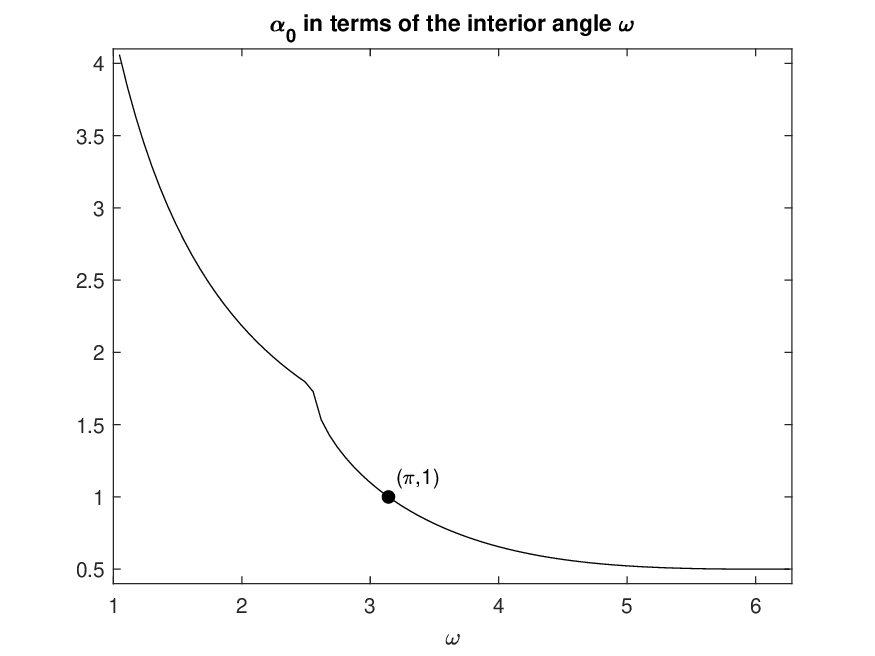}}
\caption{$\alpha_0$ in terms of the largest interior angle $\omega$. }\label{Regularity}
 \end{figure}

\subsection{$C^0$ interior penalty method} 
Let $\mathcal{T}_h$ be a triangulation of domain $\Omega$ satisfying $\overline{\Omega} = \sum_{K \in \mathcal{T}_h}\overline{K}$. We denote the sets of  interior and boundary edges of $\mathcal{T}_h$ by $\mathcal{E}_I$ and $\mathcal{E}_B$, respectively. We also set $\mathcal{E}_h = \mathcal{E}_I\cup\mathcal{E}_B$.  We further denote the mesh size of $K \in \mathcal{T}_h$ by $h_K = \text{diam}(K)$, and denote by $h = \max_{K \in \mathcal{T}_h}h_K$. The length of an edge $e \in \mathcal{E}_h$ is denoted by $h_e$. Here special attention has to be paid that the elements of $\mathcal{T}_h$ are shape-regularity, it implies that mesh $\mathcal{T}_h$ is locally quasi-uniform, i.e., if two elements $K_i$ and $K_j$ satisfy $\overline{K_i} \cap \overline{K_j} \neq \emptyset$, there exists a constant $C> 1$ such that, 
\begin{align}\label{meshcond}
C^{-1}h_{K_i} \le h_{K_j} \le C h_{K_i}.
\end{align}
Throughout the paper, we denote $K_0 \in \mathcal{T}_h$ by one element such that the singular point $\mathbf{x}_0 \in \overline{K}_0$, where $\overline{K}_0$ is the closure of $K_0$. If $\mathbf{x}_0$ lies on an inner edge, either of the two triangles sharing that edge can be chosen as $K_{0}$. Similarly, if $\mathbf{x}_0$ is a vertex of several triangles, any one of these triangles can be chosen as $K_0$. The diameter of $K_0$ is denoted by $h_{K_0}$.

We define the broken Sobolev space $H^r(\Omega,\mathcal{T}_h)$ associated with the triangulation $\mathcal{T}_h$ by
\begin{equation*}
\begin{aligned}
H^r(\Omega,\mathcal{T}_h) = \{v\in L^2(\Omega):\,v|_K \in H^r(K),\quad \forall K \in \mathcal{T}_h\}.
\end{aligned}
\end{equation*}
The space $H^r(\Omega,\mathcal{T}_h)$ equipped with the broken Sobolev norm and seminorm 
\begin{equation*}
\begin{aligned}
\|u\|_{r,\mathcal{T}_h} = \left(\sum_{K \in\mathcal{T}_h} \|u\|^2_{H^r(K)}\right)^{\frac{1}{2}},\quad |u|_{r,\mathcal{T}_h} = \left(\sum_{K \in\mathcal{T}_h} |u|^2_{H^r(K)}\right)^{\frac{1}{2}}
\end{aligned}
\end{equation*}

The $C^0$ finite element space is 
\begin{equation}\label{space}
\begin{aligned}
V_{h,0}^m = \{v_h \in H_0^1(\Omega) \cap C^0(\Omega):\,v_h|_K \in P_m(K),\,m\ge 2,\quad \forall K \in \mathcal{T}_h\},
\end{aligned}
\end{equation}
where $P_m(K)$ is the space of polynomials of degree less than or equal to $m$ on element $K$. 

For each $e \in \mathcal{E}_{I}$, we denote $K^+$ and $K^-$ by two adjacent triangles that share one common edge $e$. The unit outward normal vector $\mathbf{n}$ is oriented from $K^+$ to $K^-$. We may designate as $K^+$ that with the higher of the indices. When $e \in \mathcal{E}_B$, let $K^+$ be the element with the edge $e$ and denote by $\mathbf{n}$ a unit outward normal vector to $\partial K^+$. For any $e \in \mathcal{E}_h$, denote by $v^+$ and $v^-$ the two traces of $v$ along the edge $e$. For a scalar function $v$ and a vector function $\mathbf{q}$ that may be discontinuous across $e$, we define the following jumps:
\begin{align*}
[\![\mathbf{q}]\!]=
\begin{cases}
(\mathbf{q}^+- \mathbf{q}^-)\cdot\mathbf{n}, \quad & e \in \mathcal{E}_I,\\
\mathbf{q}^+\cdot\mathbf{n},\quad  & e \in \mathcal{E}_B,
\end{cases}
\qquad 
[\![v]\!]=
\begin{cases}
(v^+- v^-)\mathbf{n}, \quad & e \in \mathcal{E}_I,\\
v^+\mathbf{n},\quad  & e \in \mathcal{E}_B,
\end{cases}
\end{align*}
and averages 
\begin{align*}
\{\!\{\mathbf{q}\}\!\} = \begin{cases}
\frac{1}{2}\,(\mathbf{q}^+ + \mathbf{q}^-),\quad &e \in \mathcal{E}_I,\\
\mathbf{q}^+,\quad &e \in \mathcal{E}_B.
\end{cases}
\qquad 
\{\!\{v\}\!\} =
\begin{cases}
\frac{1}{2}\,(v^+ + v^-),\quad &e \in \mathcal{E}_I,\\
v^+,\quad &e \in \mathcal{E}_B.
\end{cases}
\end{align*}
According to above definition, for $\forall v \in H^r(\Omega,\mathcal{T}_h)$ and $\forall \mathbf{q} \in [H^r(\Omega,\mathcal{T}_h)]^2$, it is clearly that
\begin{align}\label{magic formula}
[\![ \mathbf{q}v]\!] = \{\!\{v\}\!\}[\![ \mathbf{q}]\!] + \{\!\{\mathbf{q}\}\!\}\cdot[\![ v]\!].
\end{align}
The following identity can be verified by simple algebraic manipulation
\begin{align}\label{magic formula2}
\sum_{K \in \mathcal{T}_h}\int_{\partial K}v\,\mathbf{q}\cdot\mathbf{n} \,\mathrm{d}s = \sum_{e \in \mathcal{E}_I}\int_{e}\{\!\{\mathbf{q}\}\!\}\cdot[\![v]\!]\,\mathrm{d}s + \sum_{e \in \mathcal{E}_h}\int_{e}\{\!\{v\}\!\}[\![\mathbf{q}]\!]\,\mathrm{d}s.
\end{align}

The $C^0$ interior penalty method for \eqref{eq:bh1} is to find $u_h \in V_{h,0}^m$ such that \cite{RH2012} 
\begin{align}\label{discrete form}
A_h(u_h,v_h) = v_h(\mathbf{x}_0)\quad \forall v_h \in V_{h,0}^m,
\end{align}
where the bilinear form 
\begin{align}\label{bilinear form}
A_h(u_h,v_h):= \sum_{K \in \mathcal{T}_h}&\int_K \Delta u_h \Delta v_h \,\mathrm{d}\mathbf{x} - \sum_{e \in \mathcal{E}_h}\left(\int_{e} \{\!\{\Delta u_h\}\!\} \,[\![\nabla v_h]\!]\,\mathrm{d}s +  \int_{e}\{\!\{\Delta v_h\}\!\} \,[\![\nabla u_h]\!]\,\mathrm{d}s \right)\nonumber\\
&+ \sum_{e \in \mathcal{E}_h}\frac{\beta}{h_e}\int_{e}[\![\nabla u_h]\!]\,[\![\nabla v_h]\!]\,\mathrm{d}s.
\end{align}
Here, the penalty parameter $\beta$ needs to be large enough to ensure the stability of the $C^0$ interior penalty method.
Define the energy norm by
\begin{equation}
\begin{aligned}
|||v|||^2 := a_h(v,v) + \sum_{e \in \mathcal{E}_h}\frac{\beta}{h_e}\|[\![\nabla v]\!]\|_{L^2(e)}^2, \quad \forall v \in H^2(\Omega,\mathcal{T}_h),
\end{aligned}\label{energy norm}
\end{equation}
where 
\begin{equation*}
\begin{aligned}
a_h(u,v) := \sum_{K \in \mathcal{T}_h}\int_K \Delta u \Delta v\,\mathrm{dx}.
\end{aligned}
\end{equation*}
It can be observed that $|||\cdot|||$ defines a norm on the space $H^2(\Omega,\mathcal{T}_h)$. 

Recall that $a(\cdot,\cdot)$ is defined in \eqref{eqn.weak}, we can observe that 
\begin{equation}
\begin{aligned}
a_h(v,v) = a(v,v)  \approx |v|^2_{H^2(\Omega)},\quad \forall v \in H^2_0(\Omega).\label{H2 property}
\end{aligned}
\end{equation}
The statement $ \approx $ represents equivalence.  Then the following inequalities hold.

\begin{lem}[Continuity and coercivity \cite{RH2012}]\label{property of Ah}
For sufficiently large $\beta$, there exists positive constants $C_s$ and $C_b$ for the bilinear form \eqref{bilinear form}, such that
\begin{align}
&|A_h(u_h,v_h)|\le C_b|||u_h|||*|||v_h|||,\quad \forall u_h,\,v_h \in V_{h,0}^m,\\
&A_h(v_h,v_h) \ge C_s|||v_h|||^2,\quad \forall v_h \in V_{h,0}^m.
\end{align}
\end{lem}
By \Cref{property of Ah} and the Lax–Milgram Theorem, the discretized problem \eqref{discrete form} admits a unique solution.

\subsection{A priori error estimate}
Given the regularity outlined in \Cref{udecompthm}, we review the following results, which are extensively used in the a priori and a posteriori error estimates.
\begin{lem}[Trace inequality \cite{BS2008}]\label{Trace inequality}
For any element $K \in \mathcal{T}_h$ and $e  \subset \partial K$, it follows
\begin{align*}
&\|v\|_{L^2(e)} \le C h_K^{-1/2}\left(\|v\|_{L^2(K)} + h_K\|\nabla v\|_{L^2(K)}\right),\qquad\qquad\, \forall v \in H^{1}(K),\\
&\|\partial_\mathbf{n} v\|_{L^2(e)} \le C h_K^{-1/2}\left(\|\nabla v\|_{L^2(K)} + h_K\|\nabla^2 v\|_{L^2(K)}\right),\qquad \forall v \in H^{2}(K).
\end{align*}
\end{lem}

\begin{lem}[Inverse inequality \cite{BS2008}]\label{Inverse inequality}
For any element $K\in\mathcal{T}_h$, $v \in P_{m}(K)$, and $e \subset \partial K$, it follows
\begin{align*}
&\|v\|_{L^2(e)} \le C h_K^{-1/2}\|v\|_{L^2(K)},\\
&\|\partial_\mathbf{n} v\|_{L^2(e)} \le C  h_K^{-1/2}\|\nabla v\|_{L^2(K)},\\
&\|\nabla^j v\|_{L^2(K)}\le C h_K^{-j}\|v\|_{L^2(K)},\qquad \forall \, 0\le j\le m.
\end{align*}
\end{lem}
\begin{lem}[Interpolation error estimate \cite{BGS2010}]\label{interpolation estimate}
Let $\Pi_h: H^2_0(\Omega) \to V_{h,0}^m$ be the standard Lagrange nodal interpolation operator, then it follows
\begin{align*}
&|\phi - \Pi_h \phi|_{H^l(K)} \le C h_K^{2-l}|\phi|_{H^2(K)},\qquad \,\,0 \le l \le 2,\\
&|\phi- \Pi_h \phi|_{H^r(\partial K)}\le C h_K^{3/2-r}|\phi|_{H^2(K)},\quad r=0,1.
\end{align*}
\end{lem}

Based on the preparations above and the analysis in \cite{BS2005}, the following a priori error estimate can be derived for the solution of the $C^0$ interior penalty method.
\begin{lem}\label{prioriest}
Let $u \in H^2_0(\Omega)$ be the solution of equation \eqref{eq:bh1}, and $u_h$ be the approximation solution of \eqref{discrete form}. Then it follows
\begin{align}
|||u - u_h||| \le C h^{\min\{ 1-\epsilon, \alpha \}}\|u\|_{H^{\min\{3- \epsilon, 2+\alpha\}}(\Omega)},
\end{align}
where $\alpha< \alpha_0$ with $\alpha_0$ given in \eqref{alpha0}.
\end{lem}

\section{Residual-based a posteriori error estimators}\label{sec3}
To improve the convergence rate of the $C^0$ interior penalty method in \Cref{prioriest}, we propose an adaptive $C^0$ interior penalty method in this section.
Specifically, we introduce two residual-type a posteriori error estimators for problem \eqref{eq:bh1}. Based on the derived error estimators and a bisection mesh refinement method, we then develop an adaptive $C^0$ interior penalty algorithm.

\subsection{A posteriori error estimation based on primal problem} 
 
The first type of error estimator is obtained in a straightforward manner based on the problem \eqref{eq:bh1}. Theoretically, we establish upper and lower bounds to ensure the reliability and efficiency of the proposed estimator.

Let $u_h \in V_{h,0}^m$ be the approximation solution obtained by the $C^0$ interior penalty method \eqref{discrete form} for problem \eqref{eq:bh1}. For each $K \in \mathcal{T}_h$, $\mathcal{E}_K$ represents the set of three edges of element $K$. Denote the set of all mesh nodes of the triangulation $\mathcal{T}_h$ by $\mathcal{N}$. The number of nodes is equal to the degrees of freedom. For example, $\mathcal{N}$ includes vertices and edge center points for the quadratic polynomial approximation. 
We propose the following residual-based a posteriori error estimator on $K \in \mathcal{T}_h$ involving the location of the Dirac point in the mesh
\begin{align}
\eta_{K}(u_h) = 
\begin{cases}
\left(h_{K_0}^2 +\overline{\eta}_{K_0}^2\right)^{1/2},&\qquad \text{if }K = K_0\text{ and } \mathbf{x}_0 \notin \mathcal{N}, \\
\overline{\eta}_K,&\qquad \text{otherwise},
\label{lobal indicator}
\end{cases}
\end{align}
where
\begin{align}\label{lobal indicator part}
\overline{\eta}_K(u_h) &= \left(\eta_{1,K}^2 +\sum_{e \in \mathcal{E}_K \cap \mathcal{E}_h}\alpha_e\eta_{2,e}^2 +  \sum_{e \in \mathcal{E}_K \cap \mathcal{E}_I}\alpha_e\eta_{3,e}^2 + \sum_{e \in \mathcal{E}_K \cap \mathcal{E}_I}\alpha_e\eta_{4,e}^2\right)^{1/2},
\end{align}
with  $\alpha_e = 1$ for $e \in \mathcal{E}_B$, $\alpha_e = 1/2$ for $e \in \mathcal{E}_I$, and 
\begin{align}
&\eta_{1,K} = h_K^2\|\Delta^2 u_h\|_{L^2(K)},\qquad\qquad\qquad
\eta_{2,e} = \beta h_e^{-1/2}\|[\![\nabla u_h]\!]\|_{L^2(e)},\label{part 1}\\
&\eta_{3,e} = h_e^{1/2}\|[\![\Delta u_h]\!]\|_{L^2(e)},\qquad\qquad\qquad
\eta_{4,e} = h_e^{3/2}\|[\![\nabla \Delta u_h]\!]\|_{L^2(e)}.\label{part 2}
\end{align}
Then the corresponding global error estimator is given by
\begin{align}\label{global indicator}
\eta(u_h) =
\begin{cases}
\left(h_{K_0}^2 + \sum\limits_{K \in \mathcal{T}_h}\eta_K^2(u_h)\right)^{1/2},&\qquad \text{if } \mathbf{x}_0 \notin \mathcal{N}, \\
\left(\sum\limits_{K \in \mathcal{T}_h}\eta_K^2(u_h)\right)^{1/2},&\qquad \text{if } \mathbf{x}_0 \in \mathcal{N}. \\
\end{cases}
\end{align}
If $\mathbf{x}_0$ is not a vertex of the triangulation, an additional term $h_{K_0}$ appears in the indicators corresponding to the triangle $\mathbf{x}_0 \in \overline{K}_0$.

To derive the reliability bound of a posteriori error estimator, we introduce the linear operator $E_h$ mapping elements in $V_{h,0}^m$ onto a $C^1$ conforming macro-elements space $S_h^{m+2}$ of degree $m+2$. 
For the detailed definition of this $C^1$ conforming macro-elements, refer to \cite{GPJ2009, DDPS1979}.
For the convenience of readers, we provide a brief review of the high-order versions of the classical Hsieh-Clough-Tocher macro-element.

\begin{definition}[\cite{GPJ2009}]
Let element $K \in \mathcal{T}_h$. For $m\ge 2$, a macro-element of degree $m+2$ is a nodal finite element $(K,\,\widetilde{P}_{m+2},\,\widetilde{N}_{m+2})$. Here, the element $K$ consists of subtriangles $K_i$, $i = 1,2,3$ satisfying $\overline K = \cup_{i=1}^3 \overline{K_i}$ as shown in \Cref{P5cme}. The local element space $\widetilde{P}_{m+2}$  on $K$ is defined by
\begin{align}
\widetilde{P}_{m+2}: = \{v \in C^1(K):\,v|_K \in P_{m+2}(K_i),\,i =1,2,3\}.
\end{align}
The degrees of freedom $\widetilde{N}_{m+2}$ on $K$ consist of all the following values:
\begin{itemize}
\item {The value and the first (partial) derivatives at the vertices of $K$;}
\item {The value at $m-1$ distinct points in the interior of each exterior edge of $K$;}
\item {The normal derivative at $m$ distinct points in the interior of each exterior edge of $K$;}
\item {The value and the first (partial) derivatives at the common vertex of all $K_i$, where $i = 1,2,3$;}
\item {The value and the normal derivative at $m-2$ distinct points in the interior of each edge of the $K_i$, where $i = 1,2,3$, that is not an edge of $K$;}
\item {The value at $(m-2)(m-3)/2$ distinct points in the interior of each $K_i$ chosen so that, if a polynomial of degree $m-4$ vanishes at those points, then it vanishes identically.}
\end{itemize}
\end{definition}
For example, the $\widetilde{P}_5$ macro-element is a $C^1$ extension of the $C^0$ Lagrange element that consists of $P_3$ polynomials.
These elements are illustrated in Figure \ref{C1 conforming element}, where we use the solid dot ($ \bullet$) to denote the value of the shape functions, the circle ($ \bigcirc$) to denote the value of all the first (partial) derivatives of the shape functions, and the arrow ($\uparrow$) to denote the value of the normal derivatives. 
\begin{figure}
\centering
\subfigure[]{
\begin{tikzpicture}[scale=0.2]
\draw[thick]
(-12,-6) -- (12,-6) -- (0,12) -- (-12,-6);

\draw[ultra thick]  (-12,-6) node {$\bullet$};
\draw[ultra thick]  (12,-6) node {$\bullet$};
\draw[ultra thick]  (0,12) node {$\bullet$};
\draw[ultra thick]  (-4,6) node {$\bullet$};
\draw[ultra thick]  (4,6) node {$\bullet$};
\draw[ultra thick]  (-4,-6) node {$\bullet$};
\draw[ultra thick]  (4,-6) node {$\bullet$};
\draw[ultra thick]  (-8,0) node {$\bullet$};
\draw[ultra thick]  (8,0) node {$\bullet$};
\draw[ultra thick]  (0,0) node {$\bullet$};
\end{tikzpicture}
}
\subfigure[ ]{
\begin{tikzpicture}[scale=0.2]
\draw[thick]
(-12,-6) -- (12,-6) -- (0,12) -- (-12,-6);
\draw[thick] (-12,-6) -- (0,0);
\draw[thick] (12,-6) -- (0,0);
\draw[thick] (0,12) -- (0,0);

\draw[ultra thick]  (-12,-6) node {$\bullet$};
\draw[ultra thick]  (12,-6) node {$\bullet$};
\draw[ultra thick]  (0,12) node {$\bullet$};
\draw[ultra thick]  (-4,6) node {$\bullet$};
\draw[ultra thick]  (4,6) node {$\bullet$};
\draw[ultra thick]  (-4,-6) node {$\bullet$};
\draw[ultra thick]  (4,-6) node {$\bullet$};
\draw[ultra thick]  (-8,0) node {$\bullet$};
\draw[ultra thick]  (8,0) node {$\bullet$};
\draw[ultra thick]  (0,0) node {$\bullet$};
\draw[ultra thick]  (0,6) node {$\bullet$};
\draw[ultra thick]  (-6,-3) node {$\bullet$};
\draw[ultra thick]  (6,-3) node {$\bullet$};

\draw[->,thick](-8,-6)--(-8,-8.3);
\draw[->,thick](0,-6)--(0,-8.3);
\draw[->,thick](8,-6)--(8,-8.3);

\draw[->,thick](-2,9)--(-4,10.5);
\draw[->,thick](-6,3)--(-8,4.5);
\draw[->,thick](-10,-3)--(-12,-1.5);

\draw[->,thick](2,9)--(4,10.5);
\draw[->,thick](6,3)--(8,4.5);
\draw[->,thick](10,-3)--(12,-1.5);

\draw[->,thick](0,6)--(2.5,6);
\draw[->,thick](-6,-3)--(-7,-1.1);
\draw[->,thick](6,-3)--(5.2,-5);

\draw (0,0) circle [radius = 1];
\draw (-12,-6) circle [radius = 1];
\draw (12,-6) circle [radius = 1];
\draw (0,12) circle [radius = 1];
\end{tikzpicture}\label{P5cme}
}
\vspace*{-15pt}
    \caption{(a) A $P_3$ Lagrange element. (b) A $\widetilde{P}_5$ $C^1$ conforming macro element.} \label{C1 conforming element}
\end{figure}
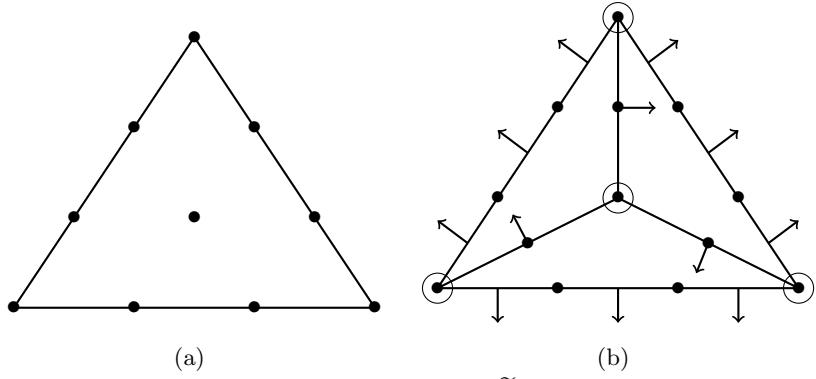

Denote by $\omega_v$ the set of elements containing a node $v \in \mathcal{N}$, and let $\sharp \omega_v$ denote the number of elements in $\omega_v$. 
We construct $E_h$ by averaging the nodal function values as follows:
\begin{equation}\label{linear operator}
\begin{aligned}
N_v(E_h(u_h)) =
\begin{cases}
\frac{1}{\sharp \omega_v}\sum\limits_{K \in\omega_v}N_v(u_h|_K),\quad &\text{if} \,\, v  \not \in \partial{\Omega};\\
0,\quad &\text{if} \,\, v \in \partial{\Omega}.
\end{cases}
\end{aligned}
\end{equation}
Here, $N_v$ represents either the nodal value of a shape function, its first partial derivatives, or its normal derivative at $v$, where $v$ is any node in the macro-elements space $S_h^{m+2}$. 

\begin{lem}\label{recovery operator estimate} \cite[Lemma 2.9]{GPJ2009}
Let $u_h \in V_{h,0}^m$ be the solution of \eqref{discrete form}. Then there exists an operator $E_h:\,V_{h,0}^m\to S_h^{m+2}\cap H^2_0(\Omega)$ that satisfies the following bound
\begin{equation}
\begin{aligned}
\sum_{K \in \mathcal{T}_h}|u_h - E_h(u_h)|_{H^{l}(K)}^2  \le C \sum_{e \in \mathcal{E}_h} h_e^{3-2l}\|[\![\nabla u_h]\!]\|_{L^2(e)}^2, \quad l=0,1,2. \label{recovery estimate}
\end{aligned}
\end{equation}
\end{lem}
\begin{rem}
Different from the estimate in \cite{GPJ2009}, the right-hand side of \eqref{recovery estimate} here does not include the term $\sum_{e \in \mathcal{E}_h} h_e^{1-2l}\|[\![u_h]\!]\|_{L^2(e)}^2$. This is because $V_h^m \subset C^0(\Omega)$, then the jump of $u_h$ across any edge $e \in \mathcal{E}_I$ is zero, and $u_h|_e = 0$ for $e \in \mathcal{E}_B$. 
\end{rem}

We further introduce the following result.
\begin{lem}[Weak continuity]\label{Weak continuity}
Let $u$ be the solution of the problem (\ref{eq:bh1}).
Then for any $e \in \mathcal{E}_I$, 
\begin{align}
&\int_e[\![u]\!]\cdot\mathbf{q}\,\mathrm{d}s =0,\qquad \,\,\,\forall \mathbf{q} \in [L^2(e)]^2,\label{intuq1}\\
&\int_e[\![\nabla u]\!]\ v\,\mathrm{d}s =  0,\qquad\,\, \forall v \in L^2(e),\label{intuq2}\\
&\int_e[\![\Delta  u]\!]\cdot\mathbf{q}\,\mathrm{d}s =  0,\qquad \forall \mathbf{q} \in [L^2(e)]^2,\label{intuq3}\\
& \int_e [\![\nabla \Delta u]\!]\ v \,\mathrm{d}s=0, \qquad \forall v \in L^2(e).\label{intuq4}
\end{align}
\end{lem}
\begin{proof}
Given that $u\in H^{\min\{3- \epsilon, 2+\alpha\}}(\Omega) \cap H_0^2(\Omega)$, \eqref{intuq1} and \eqref{intuq2} follow immediately. We denote the set of the neighborhoods of the domain corners by $\Omega_{C}$ and of the singular point $\mathbf{x}_0$ by $\Omega_{\mathbf{x}_0}$. Then $u \in H^4(\Omega \setminus (\Omega_C \cup \Omega_{\mathbf{x}_0} ))$. Therefore, \eqref{intuq3} and \eqref{intuq4} hold for any $e \in \Omega \setminus (\Omega_C \cup \Omega_{\mathbf{x}_0} )$. To this end, we show that they also hold for $e \in \Omega_C \cup \Omega_{\mathbf{x}_0}$. In the neighborhood of $\mathbf{x}_0$, the solution $u$ can be decomposed as $u = u_R + u_{\mathbf{x}_0}$, where $u_R \in H^4(\Omega_{\mathbf{x}_0})$ and $u_{\mathbf{x}_0} = -\frac{|\mathbf{x} - \mathbf{x}_0|}{8 \pi} \ln |\mathbf{x} - \mathbf{x}_0|$ is the fundamental solution of the biharmonic equation with the Dirac delta function $\delta_{\mathbf{x}_0}$. It is straightforward to verify that \eqref{intuq3} and \eqref{intuq4} hold for any $e \in \Omega_{\mathbf{x}_0}$, and they can be similarly proved for any $e \in \Omega_C$.
\end{proof}
Let $\omega_e$ be the collection of two adjacent elements that share the common edge $e$. Specially, we define
\begin{align*}
\omega_e^0 = \{K \in \mathcal{T}_h: \partial K \cap \partial K_0 =e\}.
\end{align*}
For any $K \in \omega_e$, by \eqref{meshcond} and the shape regular assumption, 
there exist positive constants $C_1$ and $C_2$ such that
$$
C_1h_K \le h_e \le C_2h_K.
$$

Next, we are ready to present one of the main results.
\begin{thm}[Reliability]\label{th1}
Let $u$ be the solution of \eqref{eq:bh1} and $u_h \in V_{h,0}^m$ be the solution of \eqref{discrete form}. Then the residual-based a posteriori error estimator $\eta$ satisfies the global bound
\begin{equation}
|||u-u_h||| \le C \eta.
\end{equation}
\end{thm}
\begin{proof}
Recall the energy norm $|||\cdot|||$ in \eqref{energy norm}. It follows
\begin{equation*}
\begin{aligned}
|||u-u_h|||^2= \sum_{K \in \mathcal{T}_h}|u-u_h|_{H^2(K)}^2 + \sum_{e \in \mathcal{E}_h}\frac{\beta}{h_e}\|[\![\nabla (u-u_h)]\!]\|_{L^2(e)}^2. 
\end{aligned}
\end{equation*}
According to the Lemma \ref{Weak continuity}, it is clearly that
\begin{equation*}
\begin{aligned}
\|[\![\nabla (u-u_h)]\!]\|_{L^2(e)}^2 = \int_e [\![\nabla u]\!]^2 - 2[\![\nabla u]\!][\![\nabla u_h]\!] +  [\![\nabla u_h]\!]^2\,\mathrm{d}s  = \int_e [\![\nabla u_h]\!]^2\,\mathrm{d}s = \|[\![\nabla u_h]\!]\|_{L^2(e)}^2,
\end{aligned}
\end{equation*}
then
\begin{equation*}
\begin{aligned}
\sum_{e \in \mathcal{E}_h}\frac{\beta}{h_e}\|[\![\nabla (u-u_h)]\!]\|_{L^2(e)}^2 = \sum_{e \in \mathcal{E}_h}\frac{\beta}{h_e}\|[\![\nabla u_h]\!]\|_{L^2(e)}^2 \le C \sum_{e \in \mathcal{E}_h}\eta_{2,e}^2.
\end{aligned}
\end{equation*}

Let $\chi = E_hu_h \in H_0^2(\Omega)$, the triangle inequality gives
\begin{align}\label{th1 pf1}
\sum_{K \in \mathcal{T}_h}|u-u_h|^2_{H^2(K)} \le  \sum_{K \in \mathcal{T}_h}|\chi-u_h|^2_{H^2(K)} + \sum_{K \in \mathcal{T}_h}|u-\chi|^2_{H^2(K)}.
\end{align}
By Lemma \ref{recovery operator estimate}, it holds
\begin{align*}
\sum_{K \in \mathcal{T}_h}|\chi-u_h|_{H^2(K)}^2 \le C\sum_{e \in \mathcal{E}_h}\eta_{2,e}^2.
\end{align*}
To this end, it suffices to show that the second term on the right-hand side of \eqref{th1 pf1} satisfies
\begin{align}\label{uchitoeta}
\left(\sum_{K \in \mathcal{T}_h}|u-\chi|^2_{H^2(K)}\right)^{1/2}\le C \eta.
\end{align}
By \eqref{H2 property} and duality argument,
\begin{align}\label{th1 pf2}
\left(\sum_{K \in \mathcal{T}_h}|u-\chi|^2_{H^2(K)}\right)^{1/2} = |u-\chi|_{H^2(\Omega)} \le \sup_{\phi \in H^2_0(\Omega)\backslash{\{0\}}}\frac{a(u-\chi,\phi)}{|\phi|_{H^2(\Omega)}}.
\end{align}
Denote the continuous interpolation polynomial of $\phi$ by $\phi_I = \Pi_h \phi \in V_{h,0}^m$. By \eqref{H2 property}, \eqref{eqn.weak} and \eqref{discrete form},
\begin{align}
a(u-\chi,\phi) & = a(u,\phi) + a_h(u_h,\phi)- a(\chi,\phi) - a_h(u_h,\phi)\nonumber\\
&= a(u,\phi) +a_h(u_h-\chi,\phi)  - a_h(u_h,\phi)\nonumber\\
&=\langle \delta_{\mathbf{x}_0}, \phi \rangle+a_h(u_h-\chi,\phi) -a_h(u_h,\phi_I)-a_h(u_h,\phi-\phi_I) \nonumber\\
&=\langle \delta_{\mathbf{x}_0}, \phi \rangle -A_h(u_h,\phi_I) +a_h(u_h-\chi,\phi)+ A_h(u_h,\phi_I)- a_h(u_h,\phi_I) - a_h(u_h,\phi-\phi_I) \nonumber\\
& = \langle \delta_{\mathbf{x}_0}, \phi-\phi_I \rangle +a_h(u_h-\chi,\phi)+ A_h(u_h,\phi_I)- a_h(u_h,\phi_I) - a_h(u_h,\phi-\phi_I).\label{th1 pf0}
\end{align}

To estimate the first term on the right-hand side of \eqref{th1 pf0}, we consider the following three possibilities based on the different locations of $\mathbf{x}_0$. Recall that $\mathcal{N}$ is the set of all mesh nodes of the triangulation.\\
(1) If $\mathbf{x}_0 \in \mathcal{N}$, the values of $\phi_I$ and $\phi$ are equal at the node $\mathbf{x}_0$, i.e., $\phi(\mathbf{x}_0) = \phi_I(\mathbf{x}_0)$, then
\begin{align}\label{deltaest1}
\langle \delta_{\mathbf{x}_0}, \phi-\phi_I \rangle = \phi(\mathbf{x}_0)-\phi_I(\mathbf{x}_0) = 0.
\end{align}
(2) If $\mathbf{x}_0\notin \mathcal{N}$, but it is located inside one element $K_0\in \mathcal{T}_h$, it follows
\begin{align}
\langle \delta_{\mathbf{x}_0}, \phi-\phi_h \rangle \le \|\phi-\phi_I\|_{L^{\infty}(K_0)} \le C h_{K_0}|\phi|_{H^2(K_0)} \le Ch_{K_0}|\phi|_{H^2(\Omega)}.
\end{align}
(3) If $\mathbf{x}_0\notin \mathcal{N}$, but it belongs to an internal edge $e$, then 
\begin{align}
\langle \delta_{\mathbf{x}_0}, \phi-\phi_I \rangle \le  \|\phi-\phi_I\|_{L^{\infty}(\omega_e)} \le C h_{K_0}|\phi|_{H^2(\omega_e)} \le C h_{K_0}|\phi|_{H^2(\Omega)}.
\end{align}

By the Cauchy-Schwarz inequality and Lemma \ref{recovery operator estimate}, the second term in \eqref{th1 pf0} follows
\begin{align}\label{th1 pf3}
a_h(u_h - \chi,\phi) \le \left(\sum_{K \in \mathcal{T}_h}|u_h-\chi|_{H^2(K)}^2\right)^{1/2}|\phi|_{H^2(\Omega)} \le C \left(\sum_{e \in \mathcal{E}_h}\eta_{2,e}^2\right)^{1/2}|\phi|_{H^2(\Omega)}.
\end{align}

According to the definition of $A_h(\cdot,\cdot)$ and $a_h(\cdot,\cdot)$, it is clear that
\begin{align*}
A_h(u_h,\phi_I) - a_h(u_h,\phi_I) = &-  \sum_{e \in \mathcal{E}_h}\int_{e} \{\!\{\Delta u_h\}\!\} \,[\![\nabla \phi_I]\!]\,\mathrm{d}s -  \sum_{e \in \mathcal{E}_h}\int_{e} [\![\nabla u_h]\!]\,\{\!\{\Delta \phi_I\}\!\} \,\mathrm{d}s+ \sum_{e \in \mathcal{E}_h}\frac{\beta}{h_e}\int_{e}[\![\nabla u_h]\!]\,[\![\nabla \phi_I]\!]\,\mathrm{d}s.
\end{align*}
Recall that $\phi \in H_0^2(\Omega),\,\phi_I \in V_{h,0}^m$. Then $\phi-\phi_I$ is continuous, and $(\phi-\phi_I)|_{e} = 0$ for any $e \in \mathcal{E}_B$. By Lemma \ref{Weak continuity}, it holds $\int_e[\![\nabla \phi]\!]v \, \mathrm{d}s = 0,\,\forall v \in L^2(e)$ for any $e \in \mathcal{E}_I$. 
Therefore,
\begin{align*}
\sum_{K \in \mathcal{T}_h} \int_{\partial K}\nabla(\Delta u_h)\cdot \mathbf{n}(\phi-\phi_I)\,\mathrm{d}s = \sum_{e \in \mathcal{E}_I}\int_e [\![\nabla(\Delta u_h)]\!](\phi-\phi_I)\,\mathrm{d}s,
\end{align*}
which together with the integration by parts, and \eqref{magic formula2} yields
\begin{align*}
-a_h(u_h,\phi-\phi_I) &= -\sum_{K \in \mathcal{T}_h}\int_K\Delta u_h \Delta (\phi-\phi_I)\,\mathrm{d}\mathbf{x}\\
& = \sum_{K \in \mathcal{T}_h}\left(\int_K\nabla (\Delta u_h)\cdot\nabla(\phi-\phi_I)\,\mathrm{d}\mathbf{x} - \int_{\partial K}\Delta u_h \nabla(\phi-\phi_I)\cdot \mathbf{n}\,\mathrm{d}s\right)\\
& = -\sum_{K \in \mathcal{T}_h}\left(\int_K\Delta^2u_h(\phi-\phi_I)\,\mathrm{d}\mathbf{x} + \int_{\partial K}\nabla(\Delta u_h)\cdot \mathbf{n}(\phi-\phi_I)\,\mathrm{d}s-  \int_{\partial K}\Delta u_h \nabla(\phi-\phi_I)\cdot \mathbf{n}\,\mathrm{d}s\right) \\
& = -\sum_{K \in \mathcal{T}_h}\int_K\Delta^2u_h(\phi-\phi_I)\,\mathrm{d}\mathbf{x} + \sum_{e \in \mathcal{E}_I}\int_e [\![\nabla(\Delta u_h)]\!](\phi-\phi_I)\,\mathrm{d}s\\
&\quad - \sum_{e \in \mathcal{E}_I}\int_e[\![\Delta u_h]\!] \cdot \{\!\{\nabla(\phi-\phi_I)\}\!\}\,\mathrm{d}s+\sum_{e \in \mathcal{E}_h}\int_e\{\!\{\Delta u_h\}\!\}[\![\nabla\phi_I]\!]\,\mathrm{d}s.
\end{align*}
The sum of the two qualities above gives
\begin{align*}
&\left(A_h(u_h,\phi_I) - a_h(u_h,\phi_I)\right) - a_h(u_h,\phi-\phi_I) \nonumber\\
&\qquad = -\sum_{K \in \mathcal{T}_h}\int_K\Delta^2u_h(\phi-\phi_I)\,\mathrm{d}\mathbf{x}+ \sum_{e \in \mathcal{E}_I}\int_e[\![\nabla(\Delta u_h)]\!](\phi-\phi_I)\,\mathrm{d}s\nonumber\\
 &\qquad \quad- \sum_{e \in \mathcal{E}_I}\int_e[\![\Delta u_h]\!]\cdot\{\!\{\nabla(\phi-\phi_I)\}\!\}\,\mathrm{d}s-  \sum_{e \in \mathcal{E}_h}\int_{e} [\![\nabla u_h]\!]\,\{\!\{\Delta \phi_I\}\!\} \,\mathrm{d}s + \sum_{e \in \mathcal{E}_h}\frac{\beta}{h_e}\int_{e}[\![\nabla u_h]\!]\,[\![\nabla \phi_I]\!]\,\mathrm{d}s.
\end{align*}
Then, we estimate the five terms on the right hand side of the above equation one by one.
Using the Cauchy-Schwarz inequality and Lemma \ref{interpolation estimate} gives
\begin{align}\label{th1 pf4}
\sum_{K \in \mathcal{T}_h}\int_K\Delta^2u_h(\phi-\phi_I)\,\mathrm{d}\mathbf{x} & \le \sum_{K \in \mathcal{T}_h}\|\Delta^2u_h\|_{L^2(K)}\|\phi-\phi_I\|_{L^2(K)} \le C  \sum_{K \in \mathcal{T}_h}h_K^2\|\Delta^2u_h\|_{L^2(K)}|\phi|_{H^2(K)} \nonumber \\
& \le C\left(\sum_{K \in \mathcal{T}_h}\eta_{1,K}^2\right)^{1/2}|\phi|_{H^2(\Omega)},
\end{align}
\begin{align}\label{th1 pf6}
\sum_{e \in \mathcal{E}_I}\int_e[\![\nabla(\Delta u_h)]\!](\phi-\phi_I)\,\mathrm{d}s & \le \left(\sum_{e \in \mathcal{E}_I}h_e^3\int_e[\![\nabla(\Delta u_h)]\!]^2\,\mathrm{d}s\right)^{1/2}\left(\sum_{e \in \mathcal{E}_I}h_e^{-3}\int_e(\phi-\phi_I)^2\,\mathrm{d}s\right)^{1/2}\nonumber\\
& \le C\left(\sum_{e \in \mathcal{E}_I}\eta_{4,e}^2\right)^{1/2}|\phi|_{H^2(\Omega)},
\end{align}
and
\begin{align}\label{th1 pf7}
\sum_{e \in \mathcal{E}_h}\frac{\beta}{h_e}\int_{e}[\![\nabla u_h]\!]\,[\![\nabla \phi_I]\!]\,\mathrm{d}s & \le \left(\sum_{e \in \mathcal{E}_h}\int_e\beta^2h_e^{-1}[\![\nabla u_h]\!]^2\,\mathrm{d}s\right)^{1/2}\left(\sum_{e \in \mathcal{E}_h}h_e^{-1}\int_e [\![\nabla \phi_I]\!]^2\,\mathrm{d}s\right)^{1/2}\nonumber\\
& = \left(\sum_{e \in \mathcal{E}_h}\int_e\beta^2h_e^{-1}[\![\nabla u_h]\!]^2\,\mathrm{d}s\right)^{1/2} \left(\sum_{e \in \mathcal{E}_h}h_e^{-1}\int_e [\![\nabla (\phi-\phi_I)]\!]^2\,\mathrm{d}s\right)^{1/2}\nonumber\\
& \le C\left(\sum_{e \in \mathcal{E}_h}\eta_{2,e}^2\right)^{1/2}|\phi|_{H^2(\Omega)}.
\end{align}
For any $e = \partial K^+ \cap \partial K^-\in \mathcal{E}_I$, applying Lemma \ref{interpolation estimate} and Lemma \ref{Inverse inequality} gives
\begin{align*}
\|\{\!\{\nabla(\phi-\phi_I)\}\!\}\|_{L^2(e)} & \le \frac{1}{2}\left(\|\nabla(\phi-\phi_I)|_{K^+}\|_{L^2(e)} + \|\nabla(\phi-\phi_I)|_{K^-}\|_{L^2(e)} \right)\\
& \le C\left(\|\nabla(\phi-\phi_I)\|_{L^2(\partial K_+)} + \|\nabla(\phi-\phi_I)\|_{L^2(\partial K^-)} \right)\\
& \le C \left(h_{K^+}^{1/2}|\phi|_{H^2(K^+)}+ h_{K^-}^{1/2}|\phi|_{H^2(K^-)}\right),
\end{align*}
\begin{align*}
\|\{\!\{\Delta \phi_I\}\!\}\|_{L^2(e)} & \le \frac{1}{2} \left(\|\Delta \phi_I|_{K^+}\|_{L^2(e)} + \|\Delta \phi_I|_{K^-}\|_{L^2(e)} \right)\\
& \le C \left(h_{K^+}^{-1/2}\|\Delta \phi_I\|_{L^2(K^+)}+ h_{K^-}^{-1/2}\|\Delta \phi_I\|_{L^2(K^-)}\right)\\
& \le C\left[h_{K^+}^{-1/2}\left(\|\Delta (\phi_I-\phi)\|_{L^2(K^+)} + \|\Delta \phi\|_{L^2(K^+)} \right)\right.\\
&\qquad \left.+h_{K^-}^{-1/2} \left(\|\Delta (\phi_I-\phi)\|_{L^2(K^-)}  + \|\Delta \phi\|_{L^2(K^-)} \right)\right]\\
& \le C \left(h_{K^+}^{-1/2}|\phi|_{H^2(K^+)}+ h_{K^-}^{-1/2}|\phi|_{H^2(K^-)}\right).
\end{align*}
Recall that $h_e$ is equivalent to $h_{K^+}$ and $h_{K^-}$. Summing up over all interior edges $e \in \mathcal{E}_I$ leads to
\begin{align*}
\sum_{e \in \mathcal{E}_I}h_e^{-1} \|\{\!\{\nabla(\phi-\phi_I)\}\!\}\|_{L^2(e)}^2
& \le C \sum_{e \in \mathcal{E}_I} h_e^{-1}\left(h_{K^+}^{1/2}|\phi|_{H^2(K^+)}+ h_{K^-}^{1/2}|\phi|_{H^2(K^-)}\right)^2\\
& \le C \sum_{e \in \mathcal{E}_I}h_e^{-1} \left(h_{K^+}|\phi|_{H^2(K^+)}^2+ h_{K^-}|\phi|_{H^2(K^-)}^2\right)\\
& \le C |\phi|_{H^2(\Omega)}^2, 
\end{align*}
\begin{align*}
\sum_{e \in \mathcal{E}_I}h_e\|\{\!\{\Delta \phi_I\}\!\}\|_{L^2(e)}^2 & \le C \sum_{e \in \mathcal{E}_I} h_e\left(h_{K^+}^{-1/2}|\phi|_{H^2(K^+)}+ h_{K^-}^{-1/2}|\phi|_{H^2(K^-)}\right)^2\\
& \le C \sum_{e \in \mathcal{E}_I}h_e \left(h_{K^+}^{-1}|\phi|_{H^2(K^+)}^2+ h_{K^-}^{-1}|\phi|_{H^2(K^-)}^2\right)\\
& \le C |\phi|_{H^2(\Omega)}^2.
\end{align*}
Similarly, for any $e \in \partial K \cap \mathcal{E}_B$, applying Lemma \ref{interpolation estimate} and Lemma \ref{Inverse inequality} gives
\begin{align*}
\sum_{e \in \mathcal{E}_B}h_e\|\{\!\{\Delta \phi_I\}\!\}\|_{L^2(e)}^2 \le C |\phi|_{H^2(\Omega)}^2
\end{align*}
Then, by using the Cauchy-Schwarz inequality, it holds
\begin{align}\label{th1 pf5}
\sum_{e \in \mathcal{E}_I}\int_e[\![\Delta u_h]\!]\cdot\{\!\{\nabla(\phi-\phi_I)\}\!\}\,\mathrm{d}s & \le \left(\sum_{e \in \mathcal{E}_I}h_e\int_e[\![\Delta u_h]\!]^2\,\mathrm{d}s\right)^{1/2}\left(\sum_{e \in \mathcal{E}_I}h_e^{-1}\int_e\{\!\{\nabla(\phi-\phi_I)\}\!\}^2\,\mathrm{d}s\right)^{1/2}\nonumber\\
& \le C\left(\sum_{e \in \mathcal{E}_I}\eta_{3,e}^2\right)^{1/2}|\phi|_{H^2(\Omega)},
\end{align}
\begin{align}\label{th1 pf8}
\sum_{e \in \mathcal{E}_h}\int_{e} [\![\nabla u_h]\!]\,\{\!\{\Delta \phi_I\}\!\}\,\mathrm{d}s & \le \left(\sum_{e \in \mathcal{E}_h}h_e^{-1}\int_e[\![\nabla u_h]\!]^2\,\mathrm{d}s\right)^{1/2}\left(\sum_{e \in \mathcal{E}_h}h_e\int_e \{\!\{\Delta \phi_I\}\!\}^2\,\mathrm{d}s\right)^{1/2}\nonumber\\
& \le C\left(\sum_{e \in \mathcal{E}_h}\eta_{2,e}^2\right)^{1/2}|\phi|_{H^2(\Omega)}.
\end{align}
The estimates \eqref{deltaest1} - \eqref{th1 pf8} imply that  \eqref{uchitoeta} holds.

\end{proof}
To prove the efficiency of the a posteriori error estimator, we will define four types of ``bubble" functions and introduce some properties of these functions that are frequently used in error estimation and analysis. Let $\hat{K}$ be a reference triangular or rectangular element. 
If $\hat{K}$ is a triangle with the barycentric coordinates $\lambda_1, \lambda_2$ and $\lambda_3$, we denote the standard ``bubble" function in $\hat{K}$ by 
\begin{align}
b_{\hat{K}} = 27\lambda_{1}\lambda_2\lambda_3.
\end{align}
If $\hat{K}$ is a rectangle with the corresponding coordinates $x$ and $y$, we denote the ``bubble" function in $\hat{K}$ by 
\begin{align}
b_{\hat{K}} = (1-x^2)(1-y^2).
\end{align}

For a triangle $K \in \mathcal{T}_h$, we define $F_K: \hat{K} \to K$ as an affine element mapping, where $\hat{K}$ is the reference triangle. Then the ``bubble'' function $b_K$ on $K$ is defined as
\begin{align}\label{bubbbK}
b_K = 
\begin{cases}
b_{\hat{K}}\circ F_K,& \quad \text{if } K \in \mathcal{T}_h \backslash K_0,\\
\frac{|\mathbf{x} - \mathbf{x}_0|^2}{h_{K_0}^2}\,b_{\hat{K}}\circ F_{K_0},& \quad \text{if } K = K_0.
\end{cases}
\end{align}

\begin{lem}\label{bK properties}
For the element ``bubble" function $b_K$ on a triangle $K \in \mathcal{T}_h$, it holds
\begin{align}
b_K(\mathbf{x}_0) =0,\qquad\qquad
b_K(\mathbf{x}) = 0,\quad \forall \mathbf{x} \in \Omega \backslash K.\label{bK property 1}
\end{align}
Moreover, for $\forall v \in P_m(K)$, it follows
\begin{align}
\|v\|_{L^2(K)}&\le C\|b_K v\|_{L^2(K)},\qquad \|b_K^2v\|_{L^2(K)}\le C \|v\|_{L^2(K)}.\label{bK property 2}
\end{align}
\end{lem} 


For each internal edge $e \in \mathcal{E}_I$, let $T \subset \omega_e$ be the largest rhombus contained in $\omega_e$ that has $e$ as one of its diagonals. We define $F_T: \hat{K} \to T$ as an affine element mapping, where $\hat{K}$ is the reference rectangle. The ``bubble'' function $b_T$ on element $T$ is defined by
\begin{align}\label{bubbbwe}
b_{T} = 
\begin{cases}
b_{\hat{K}}\circ F_{T},& \quad \text{if } \mathbf{x}_0 \notin \overline{T},\\
\frac{|\mathbf{x} - \mathbf{x}_0|^2}{h_e^2}\,b_{\hat{K}}\circ F_{T},& \quad \text{if } \mathbf{x}_0 \in \overline{T},
\end{cases}
\end{align}
where $\overline{T}$ is the closure of $T$.

\begin{lem}\label{bwe properties}
For the ``bubble" function $b_{T}$ in rhombus $T$, it holds
\begin{align}\label{bwe property 1}
&b_{T}^3 \in C^2(\Omega)\cap H_0^2(\Omega),\quad\quad \,\, b_T(\mathbf{x}_0) = 0,\qquad
\qquad \{\!\{b_{T}^3\}\!\} = 0\quad \text{on} \quad \mathcal{E}_h \backslash e ,\\
&[\![b_{T}^3]\!] = \mathbf{0} \qquad \text{and} \qquad [\![\nabla b_{T}^3]\!] = \{\!\{\nabla b_{T}^3\}\!\}\cdot \mathbf{n}  =0\quad \text{on}\quad \mathcal{E}_h.
\end{align}
Moreover, for $\forall v \in P_m(K)$, it follows
\begin{align}
\|v\|_{L^2(e)} \le C\|b_{T}^{3/2}v\|_{L^2(e)},\qquad \|vb_{T}^3\|_{L^2(\omega_e)} \le h_e^{1/2}\|v\|_{L^2(e)}.\label{bwe property 2}
\end{align}
\end{lem} 

\begin{rem}
Inspired by the ``bubble'' functions in \cite{GHV2011}, where it is used to estimate the a posterior estimator of the discontinuous Galerkin method for the fourth-order elliptic problem with a source term in $L^2(\Omega)$. Here we modify the definition of ``bubble'' functions $b_K$ and $b_T$ to accommodate two specific cases. Compared to \cite{GHV2011}, we further consider the effects of singularity point $\mathbf{x}_0$. As a result, the new ``bubble" functions retain the favorable properties of the original ones and also obtain values of $0$ at $\mathbf{x}_0$, which is crucial for the subsequent proof of lower bounds.
\end{rem}

Let $b_l$ be an affine function that satisfies $b_l|_e =0$ and $(\nabla b_l \cdot \mathbf{n})|_e = h^{-1}_e$, where $\mathbf{n}$ is the unit normal to the edge $e$. Using the element ``bubble" function definition given above, we define an edge ``bubble" function as
\begin{align}\label{bubbbe}
b_e= b_lb_{T}^3.
\end{align}

\begin{lem}\label{be properties}
For the edge ``bubble" function $b_e$ in (\ref{bubbbe}), it holds
\begin{align}
&b_e \in C^2(\Omega)\cap H_0^2(\Omega),\qquad\qquad b_e(\mathbf{x}_0) = 0,\qquad\qquad b_e =0\quad \text{in} \quad \Omega\backslash T,\\
&[\![b_e]\!] = \mathbf{0}\qquad \text{and}\qquad [\![\nabla b_e]\!]  = \{\!\{b_e\}\!\} =0\quad \text{on}\quad \mathcal{E}_h, \\
&(\{\!\{\nabla b_e\}\!\}\cdot \mathbf{n})|_e = h_e^{-1} b_{T}^3|_e,\qquad \{\!\{\nabla b_e\}\!\} = \mathbf{0}\quad \text{on} \quad \mathcal{E}_h \backslash e ,\label{be property 1}
\end{align}
Moreover, for $\forall v \in P_m(K)$, it follows
\begin{align}
\|b_e v\|_{L^2(\omega_e)} &\le C h_e^{1/2}\|v\|_{L^2(e)}.\label{be property 2}
\end{align}
\end{lem} 

Denote $\omega_{K}$ by the collection of elements in $\mathcal{T}_h$ that share a common edge or vertex with $K$. Specially, we define the set
\begin{align}\label{omgx0}
\omega_{K_0} &= \{K \in \mathcal{T}_h: \overline{K} \cap \overline{K}_0 \neq \varnothing \},
\end{align}
and the distance of $\mathbf{x}_0$ to the boundary of $\omega_{K_0}$ is defined by $d: = \text{dist}(\mathbf{x}_0,\partial \omega_{K_0})$. We define the smooth ``bubble'' function associated with the point $\mathbf{x}_0$ by convolution of the characteristic function of the set $\{\mathbf{x} \in \Omega:|\mathbf{x}-\mathbf{x}_0|< d/4\}$ satisfying
\begin{align*}
&0 \le b_{\mathbf{x}_0}(\mathbf{x}) \le 1,\quad \forall \mathbf{x} \in \Omega,\qquad |b_{\mathbf{x}_0}(\mathbf{x})|_{m,\infty,\omega_{K_0}} \le Cd^{-m},\quad m=1,2,\\
&b_{\mathbf{x}_0}(\mathbf{x}) = 1,\quad  \forall \mathbf{x} \in \Omega: |\mathbf{x}-\mathbf{x}_0|\le \frac{d}{4},\\
&b_{\mathbf{x}_0}(\mathbf{x}) = 0,\quad \forall \mathbf{x}\in \Omega: |\mathbf{x}-\mathbf{x}_0| \ge \frac{3d}{4}.
\end{align*}

\begin{lem}\label{bx0 properties}
Assume that $\mathbf{x}_0 \in \overline{K}_0$. For the ``bubble'' function $b_{\mathbf{x}_0}$, it holds,
\begin{align}
&|b_{\mathbf{x}_0}(\mathbf{x})|_{H^m(\omega_{K_0})} \le C h_{K_0}^{1-m},\quad m=0,1,2.\label{bx0 property 1}\\
&\|b_{\mathbf{x}_0}(\mathbf{x})\|_{L^2(e)} \le C h_{e}^{1/2},\qquad\,\, e \in \partial \omega_{K_0}.\label{bx0 property 2}
\end{align}
\end{lem} 
\begin{proof}
The proof of \eqref{bx0 property 1} follow form Lemma 3.2 \cite{ABR2006}. For \eqref{bx0 property 2}, the definition of $\|\cdot\|_{L^2(e)}$ yields
\begin{align*}
\|b_{\mathbf{x}_0}(\mathbf{x})\|_{L^2(e)} = \left(\int_e b_{\mathbf{x}_0}^2(s) \mathrm{d}s\right)^{1/2} \le C \left(\int_e \mathrm{d}s\right)^{1/2} = C h_{e}^{1/2}.
\end{align*}
\end{proof}

Then we are ready to present our next main result.
\begin{thm}[Efficiency]\label{Efficiency}
For the local indicator $\eta_{K}$ defined in \eqref{lobal indicator}, there exists a positive constant $C$ independent of the mesh size satisfying
\begin{align}\label{th2}
\eta_{K} \le C|||u-u_h|||_{\omega_{K}}.
\end{align}
\end{thm}
\begin{proof}
We first prove the estimate \eqref{th2} on the element $\overline{K}_0$ whose closure contains $\mathbf{x}_0$, but $\mathbf{x}_0 \notin \mathcal{N}$.
By the definition of the energy norm, it holds
\begin{align*}
|||u-u_h|||_{\omega_{K_0}}^2= \sum_{{K} \in \omega_{K_0}}\left(|u-u_h|_{H^2({K})}^2 + \sum_{e \in \mathcal{E}_{K}}\frac{\beta}{h_e}\|[\![\nabla (u-u_h)]\!]\|_{L^2(e)}^2\right),
\end{align*}
the estimation \eqref{th2} is equivalent to
\begin{align}\label{th2 pf1}
\left(h_{K_0}^2 + \eta_{1,{K_0}}^2 +  \sum_{e \in \mathcal{E}_{K_0} \cap \mathcal{E}_I}\alpha_e\eta_{3,e}^2 + \sum_{e \in \mathcal{E}_{K_0} \cap \mathcal{E}_I}\alpha_e\eta_{4,e}^2\right)^{1/2}\le C|||u-u_h|||_{\omega_{K_0}}.
\end{align}
Then, we prove \eqref{th2 pf1} in four steps.

(i) To prove $\eta_{1,{K_0}} = h^{2}_{K_0}\|\Delta^2 u_h\|_{L^2(K_0)} \le C|||u-u_h|||_{\omega_{K_0}}$. For $\forall v \in V_{h,0}^m$, using integration by part gives
\begin{align}\label{th2 pf2}
v(\mathbf{x}_0)&-\int_{K_0}(\Delta^2 u_h)v\,\mathrm{d}\mathbf{x} = \int_{K_0}  \Delta^2(u-u_h) v\,\mathrm{d}\mathbf{x} \nonumber\\
&= -\int_{K_0} \nabla \Delta (u-u_h) \cdot\nabla v\,\mathrm{d}\mathbf{x} + \int_{\partial {K_0}}\nabla \Delta (u-u_h) \cdot \mathbf{n}v\,\mathrm{d}s\nonumber\\
& = \int_{K_0}\Delta (u-u_h)\Delta v\,\mathrm{d}\mathbf{x} - \int_{\partial {K_0}}\Delta (u-u_h) \mathbf{n}\cdot \nabla v\,\mathrm{d}s + \int_{\partial {K_0}}\nabla \Delta (u-u_h) \cdot \mathbf{n}v\,\mathrm{d}s.
\end{align}
We set $v|_{K_0} = (\Delta^2 u_h)b_{K_0}^2$ in \eqref{th2 pf2}, where $v \in H^2_0(\Omega)\cap H^2_0(K_0)$ and satisfies $v=0$ in $\Omega \backslash K_0$. Additionally, $v$ is a polynomial on $K$ with $b_{K_0}(\mathbf{x}_0) = 0$. Then, for $\forall e \in \partial K_0$, it holds that $v(\mathbf{x}_0) = v|_e = \nabla v|_e =0$. Consequently, \eqref{th2 pf2} yields
\begin{align*}
\int_{K_0}\Delta (u-u_h)\Delta v\,\mathrm{d}\mathbf{x} =-\int_{K_0}\Delta^2 u_hv\,\mathrm{d}\mathbf{x}.
\end{align*}
According to Lemma \ref{bK properties}, it holds
\begin{align*}
\|v\|_{L^2({K_0})} \le
C \|\Delta^2 u_h\|_{L^2({K_0})}.
\end{align*}
Then, by Cauchy-Schwarz inequality and Lemma \ref{Inverse inequality}, 
\begin{align}\label{Lapuhv}
\int_{K_0}(\Delta^2 u_h)v\,\mathrm{d}\mathbf{x} 
&\le \|\Delta (u-u_h)\|_{L^2(K_0)}|v|_{H^2({K_0})}  \le Ch^{-2}_{K_0}|u-u_h|_{H^2(\omega_{K_0})}\|v\|_{L^2({K_0})}\nonumber\\
& \le Ch^{-2}_{K_0}|u-u_h|_{H^2(\omega_{\mathbf{x}_0})}\|\Delta^2 u_h\|_{L^2({K_0})}.
\end{align}
Since $\Delta^2 u_h$ is a piecewise polynomial over $K_0$, according to \eqref{bK property 2} and \eqref{Lapuhv}, it holds
\begin{align*}
\|\Delta^2 u_h\|_{L^2({K_0})}^2 & \le C\int_{K_0}(\Delta^2u_h)^2b_{K_0}^2\,\mathrm{d}\mathbf{x} = C\int_{K_0}(\Delta^2u_h)v\,\mathrm{d}\mathbf{x}\nonumber\\
& \le Ch^{-2}_{K_0}|u-u_h|_{H^2(\omega_{K_0})}\|\Delta^2 u_h\|_{L^2({K_0})},
\end{align*}
which implies
\begin{align}\label{th2 pf3}
h^{2}_{K_0}\|\Delta^2 u_h\|_{L^2({K_0})}\le C|u-u_h|_{H^2({\omega_{K_0}})} \le C|||u-u_h|||_{\omega_{K_0}}.
\end{align}

(ii) To prove $\sum_{e \in \mathcal{E}_{K_0} \cap \mathcal{E}_I}\alpha_e\eta_{3,e}^2 = \sum_{e \in \mathcal{E}_{K_0}}h_e^{1/2}\|[\![\Delta u_h]\!]\|_{L^2(e)}  \le C|||u-u_h|||_{{\omega_{K_0}}}$.
For $\forall e \in \mathcal{E}_{K_0}\cap \mathcal{E}_I$, denote by $\mathcal{F}_e = \{e \in \partial K: K \in \omega_e^0\}$. By \eqref{magic formula}, the summation of \eqref{th2 pf2} over all elements $K \in \omega_e^0$ gives
\begin{align}\label{th2 pf4}
v(\mathbf{x}_0) - \sum_{K \in \omega_e^0}\int_{K}\Delta^2 u_hv\,\mathrm{d}\mathbf{x} 
& = \int_{\omega_e^0}\Delta (u-u_h)\Delta v\,\mathrm{d}\mathbf{x} +\sum_{e\in \mathcal{F}_e \cap \mathcal{E}_I}\int_e[\![v]\!]\cdot\{\!\{\nabla \Delta (u-u_h)\}\!\}\,\mathrm{d}s \nonumber\\
&\quad+ \sum_{e\in \mathcal{F}_e}\int_e\{\!\{v\}\!\}[\![\nabla \Delta (u-u_h)]\!]\,\mathrm{d}s-\sum_{e\in \mathcal{F}_e}\int_e[\![\nabla v]\!]\{\!\{\Delta (u-u_h)\}\!\}\,\mathrm{d}s\nonumber\\
&\quad- \sum_{e\in \mathcal{F}_e \cap \mathcal{E}_I}\int_e \{\!\{\nabla v\}\!\}\cdot[\![\Delta (u-u_h)]\!]\,\mathrm{d}s.
\end{align}
We take $v =\phi b_e$ in \eqref{th2 pf4}, where $\phi$ is continuous on $e \in \mathcal{E}_h$, and $\phi$ is a constant function in the normal direction to $e$ (i.e., $(\nabla \phi \cdot \mathbf{n})|_e = 0$). By $v(\mathbf{x}_0) = 0$, Lemma \ref{Weak continuity}, and Lemma \ref{be properties}, the equality \eqref{th2 pf4} reduces to
\begin{align*}
-\int_{\omega_e^0}\Delta^2 u_hv\,\mathrm{d}\mathbf{x}  = \int_e \{\!\{\nabla v\}\!\}\cdot [\![\Delta u_h]\!]\,\mathrm{d}s+\int_{\omega_e^0}\Delta (u-u_h)\Delta v\,\mathrm{d}\mathbf{x}.
\end{align*}
By using Cauchy-Schwarz inequality and Lemma \ref{Inverse inequality}, the equation above gives
\begin{align}\label{th2 pf5}
\int_e [\![\Delta u_h]\!]\cdot\{\!\{\nabla v\}\!\}\,\mathrm{d}s
&\le \|\Delta^2 u_h\|_{L^2(\omega_e^0)}\|v\|_{L^2(\omega_e^0)}
+|u-u_h|_{H^2(\omega_{e}^0)}\|\Delta v\|_{L^2(\omega_e^0)}\nonumber\\
&\le C\left(h_e^2\| \Delta^2 u_h\|_{L^2(\omega_e^0)}
+|u-u_h|_{H^2(\omega_{e}^0)}\right)h_e^{-2}\|v\|_{L^2(\omega_e^0)}
\end{align}
We extend $([\![\Delta u_h]\!]\cdot \mathbf{n})|_e$ from edge $e$ to $\omega_e^0$ by taking constants along the normal on $e$. The resulting extension $ E([\![\Delta u_h]\!]\cdot \mathbf{n})$ is a piecewise polynomial in $\omega_e$. Setting $\phi = E([\![\Delta u_h]\!]\cdot \mathbf{n})$, and using \eqref{bwe property 2} and \eqref{be property 2} yield
\begin{align}\label{th2 pf6}
&\|[\![\Delta u_h]\!]\|_{L^2(e)}^2 \le C\int_e [\![\Delta u_h]\!]^2b_{T}^3\,\mathrm{d}s= Ch_e\int_e [\![\Delta u_h]\!]\cdot\{\!\{\nabla v\}\!\}\,\mathrm{d}s\nonumber,\\
&\|v\|_{L^2(\omega_e^0)}\le  Ch^{1/2}_e\|\phi\|_{L^2(e)} = Ch^{1/2}_e\|[\![\Delta u_h]\!]\|_{L^2(e)}.
\end{align}
By \eqref{th2 pf3}, \eqref{th2 pf5}-\eqref{th2 pf6}, it follows
\begin{align*}
h_e^{1/2}\|[\![\Delta u_h]\!]\|_{L^2(e)}^2 &\le Ch_e^{3/2}\int_e [\![\Delta u_h]\!]\cdot\{\!\{\nabla v\}\!\}\,\mathrm{d}s\\
&\le C h_e^{-1/2}\left(h_e^2\|\Delta^2u_h\|_{L^2(\omega_e^0)} + |u-u_h|_{H^2(\omega_{e}^0)}\right)\|v\|_{L^2(\omega_e^0)}\\
&\le  C |||u-u_h|||_{\omega_{e}^0}\|[\![\Delta u_h]\!]\|_{L^2(e)},
\end{align*}
which gives
\begin{align}\label{hjumpdeltauh}
h_e^{1/2}\|[\![\Delta u_h]\!]\|_{L^2(e)} &\le  C |||u-u_h|||_{\omega_{e}^0}.
\end{align}
Summing up \eqref{hjumpdeltauh} over all edges $e \in  \mathcal{E}_{K_0} \cap \mathcal{E}_{I}$ yields
\begin{align}\label{th2 pf7}
\sum_{e \in \mathcal{E}_{K_0}\cap \mathcal{E}_{I}}h_e^{1/2}\|[\![\Delta u_h]\!]\|_{L^2(e)} &\le C |||u-u_h|||_{\omega_{K_0}}.
\end{align}

(iii) To prove $\sum_{e \in \mathcal{E}_{K_0} \cap \mathcal{E}_I}\alpha_e\eta_{4,e}^2 = \sum_{e \in \mathcal{E}_{K_0}}h_e^{3/2}\|[\![\nabla \Delta u_h]\!]\|_{L^2(e)}  \le C|||u-u_h|||_{\omega_{K_0}}$.
We take $v = \Phi b_{T}^3$ in \eqref{th2 pf4}, where $\Phi$ is continuous on $e \in \mathcal{E}_h$ and $(\nabla \Phi \cdot \mathbf{n})|_e = 0$. By $b_{T}^3(\mathbf{x}_0) = 0$, \Cref{Weak continuity}, and \Cref{bwe properties}, \eqref{th2 pf4} can be written as
\begin{align}\label{egradlap}
\int_{\omega_e^0}\Delta^2 u_hv\,\mathrm{d}\mathbf{x} + \int_{\omega_e^0}\Delta (u-u_h)\Delta v\,\mathrm{d}x = \int_e [\![\nabla \Delta u_h]\!]\{\!\{v\}\!\}\,\mathrm{d}s.
\end{align}
By \eqref{egradlap}, Cauchy-Schwarz inequality, and Lemma \ref{Inverse inequality},
\begin{align}\label{th2 pf8}
\int_e[\![\nabla \Delta u_h]\!]\{\!\{v\}\!\}\,\mathrm{d}s \le & C \left(|u-u_h|_{H^2(\omega_{e}^0)} + h_e^2\|\Delta^2 u_h\|_{L^2(\omega_e^0)}\right)h_e^{-2}\|v\|_{L^2(\omega_e^0)}.
\end{align}
We extend $([\![\nabla \Delta u_h]\!])|_e$ to a function $E([\![\nabla \Delta u_h]\!])$, defined over $\omega_e^0$, by taking it to be constants along lines normal to $e$. Setting $\Phi = E([\![\nabla \Delta u_h]\!])$ and using \eqref{bwe property 2} yield
\begin{align}\label{th2 pf9}
&\|v\|_{L^2(\omega_e^0)}  \le  Ch^{1/2}_e\|\Phi\|_{L^2(e)} = Ch^{1/2}_e\|[\![\nabla \Delta u_h]\!]\|_{L^2(e)},\nonumber \\
&\|[\![\nabla \Delta u_h]\!]\|^2_{L^2(e)} \le  C\int_e [\![\nabla \Delta u_h]\!]^2b_{T}^3\,\mathrm{d}s = C\int_e [\![\nabla \Delta u_h]\!]\{\!\{v\}\!\}\,\mathrm{d}s,\nonumber\\
& \qquad\qquad\qquad\,\,\,\le  C h_e^{-3/2}\left(|u-u_h|_{H^2(\omega_{e}^0)} + h_e^2\|\Delta^2 u_h\|_{L^2(\omega_e^0)}\right)\|[\![\nabla \Delta u_h]\!]\|_{L^2(e)}.
\end{align}
Form \eqref{th2 pf3} and \eqref{th2 pf8}-\eqref{th2 pf9}, 
\begin{align}\label{hegradlap}
h_e^{3/2}\|[\![\nabla \Delta u_h]\!]\|_{L^2(e)} &\le C \left(h_e
^2\|\Delta^2u_h\|_{L^2(\omega_e^0)}+|u-u_h|_{H^2(\omega_{e}^0)}\right)\le C |||u-u_h|||_{\omega_{e}^0},
\end{align}
Summing up \eqref{hegradlap} over all edges $e \in  \mathcal{E}_{K_0}\cap \mathcal{E}_{I}$ gives
\begin{align}\label{th2 pf10}
\sum_{e \in \mathcal{E}_{K_0}\cap \mathcal{E}_{I}}h_e^{3/2}\|[\![\nabla \Delta u_h]\!]\|_{L^2(e)} &\le C|||u-u_h|||_{\omega_{K_0}}.
\end{align}

(iv) To prove $h_{K_0}  \le C|||u-u_h|||_{\omega_{K_0}}$.
By the weak formulation of \eqref{eq:bh1} and Lemma \ref{bx0 properties}, it can be observed 
\begin{align}\label{th2 pf11}
1 = \langle \delta_{\mathbf{x}_0},b_{\mathbf{x}_0} \rangle &\le \left|\int_{\Omega}\Delta(u-u_h)\Delta b_{\mathbf{x}_0}\,\mathrm{d}\mathbf{x} + \int_{\Omega}\Delta u_h\Delta b_{\mathbf{x}_0}\,\mathrm{d}\mathbf{x}\right|\nonumber\\
& \le |u-u_h|_{H^2(\omega_{K_0})}|b_{\mathbf{x}_0}|_{H^2(\omega_{K_0})} + \left|\int_{\Omega}\Delta u_h\Delta b_{\mathbf{x}_0}\,\mathrm{d}\mathbf{x}\right|\nonumber\\
& \le h_{K_0}^{-1}|u-u_h|_{H^2(\omega_{K_0})} + \left|\int_{\omega_{K_0}}\Delta u_h\Delta b_{\mathbf{x}_0}\,\mathrm{d}\mathbf{x}\right|.
\end{align}
Let $\mathcal{F}_I = \{e \in \partial K\cap \mathcal{E}_I : K \in \omega_{K_0}\}$. By integration by parts, Cauchy-Schwarz inequality, Lemma \ref{bx0 properties} and Lemma \ref{Inverse inequality},
\begin{align}\label{th2 pf12}
&\int_{\omega_{K_0}}\Delta u_h\Delta b_{\mathbf{x}_0}\,\mathrm{d}\mathbf{x} =\sum_{{K} \in \omega_{K_0}} \int_{{K}}\Delta^2 u_h b_{\mathbf{x}_0}\,\mathrm{d}\mathbf{x} - \int_{\partial{K}}\nabla(\Delta u_h)\cdot\mathbf{n} b_{\mathbf{x}_0}\,\mathrm{d}s + \int_{\partial {K}}\Delta u_h \nabla b_{\mathbf{x}_0}\cdot\mathbf{n}\,\mathrm{d}s\nonumber\\
&\qquad= \sum_{{K} \in \omega_{K_0}} \int_{{K}}\Delta^2 u_h b_{\mathbf{x}_0}\,\mathrm{d}\mathbf{x} - \sum_{e \in \mathcal{F}_I}\left(\int_{e}[\![\nabla(\Delta u_h)]\!] \,b_{\mathbf{x}_0}\,\mathrm{d}s - \int_{e}[\![\Delta u_h ]\!]\cdot\nabla b_{\mathbf{x}_0}\,\mathrm{d}s\right)\nonumber\\
& \qquad\le \|\Delta^2 u_h\|_{L^2(\omega_{K_0})}\|b_{\mathbf{x}_0}\|_{L^2(\omega_{K_0})} +\sum_{e \in \mathcal{F}_I}\|[\![\nabla(\Delta u_h)]\!]\|_{L^2(e)}\|b_{\mathbf{x}_0}\|_{L^2(e)}+\sum_{e \in \mathcal{F}_I}\|[\![\Delta u_h]\!]\|_{L^2(e)}\|\nabla b_{\mathbf{x}_0}\|_{L^2(e)}\nonumber\\
&\qquad \le Ch_{K_0}^{-1}\left(h_{K_0}^2\|\Delta^2 u_h\|_{L^2(\omega_{K_0})} + \sum_{e \in \mathcal{F}_I}h_e^{3/2}\|[\![\nabla(\Delta u_h)]\!]\|_{L^2(e)} + \sum_{e \in \mathcal{F}_I}h_e^{1/2}\|[\![\Delta u_h]\!]\|_{L^2(e)}\right).
\end{align}
Inserting \eqref{th2 pf12} into \eqref{th2 pf11} yields
\begin{align}\label{th2 pf13}
1 &\le Ch_{K_0}^{-1}\bigg(|u-u_h|_{H^2(\omega_{K_0})} + h_{K}^2\|\Delta^2 u_h\|_{L^2(\omega_{K_0})} \bigg.\nonumber\\
&\quad \bigg.+ \sum_{e \in \mathcal{F}_I}h_{e}^{3/2}\|[\![\nabla(\Delta u_h)]\!]\|_{L^2(e)} + \sum_{e \in \mathcal{F}_I}h_{e}^{1/2}\|[\![\Delta u_h]\!]\|_{L^2(e)}\bigg),
\end{align}
which, together with \eqref{th2 pf3}, \eqref{th2 pf7}, and \eqref{th2 pf10}, implies
\begin{align}\label{th2 pf14}
h_{K_0} \le C |||u-u_h|||_{\omega_{K_0}}.
\end{align}
The estimate \eqref{th2 pf1} follows from \eqref{th2 pf3}, \eqref{th2 pf7}, \eqref{th2 pf10}, and \eqref{th2 pf14}.

If $\mathbf{x}_0 \in \mathcal{N}$, for $ \forall K \in \mathcal{T}_h$, the term $h_K$ vanishes in $\eta_K$ and the estimation \eqref{th2 pf1} reduces to
\begin{align}\label{th2 pf15}
\left(\eta_{1,{K}}^2 +  \sum_{e \in \mathcal{E}_{K} \cap \mathcal{E}_I}\alpha_e\eta_{3,e}^2 + \sum_{e \in \mathcal{E}_{K} \cap \mathcal{E}_I}\alpha_e\eta_{4,e}^2\right)^{1/2}\le C|||u-u_h|||_{\omega_{K}}.
\end{align}
Then, the estimate \eqref{th2 pf15} is given directly by \eqref{th2 pf3}, \eqref{th2 pf7} and \eqref{th2 pf10}.
\end{proof}
\subsection{A posteriori error estimation based on regularized problem} 
Inspired by the technique of second-order elliptic equations with a Dirac delta source term, as discussed in \cite{HW2012}, where the $L^2$  projection of the Dirac delta function is used to produce a regular solution, allowing adaptive procedures based on standard a posteriori error estimators to work efficiently.
This regularization approach using projection techniques can also be applied for the problem \eqref{eq:bh1}. 

More specifically, the Dirac delta function can be approximated by $\delta_h \in V_{h,0}^m \subset H_0^1(\Omega)$, defined as
\begin{align*}
\delta_h = 
\begin{cases}
0, \quad &\text{in } \Omega\backslash \overline{K_0},\\
\delta_{K_0},\quad &\text{in }K_0,
\end{cases}
\end{align*}
where $\delta_{K_0} \in P_m(K_0)$ satisfying
\begin{align*}
\int_{K_0} \delta_{K_0} v\,\mathrm{d}\mathbf{x} = v(\mathbf{x}_0),\quad \forall v \in P_m(K_0).
\end{align*}
Therefore, it holds
\begin{align*}
\int_{\Omega} \delta_h v\,\mathrm{d}\mathbf{x} = v(\mathbf{x}_0).
\end{align*}
Let us consider the ensuing auxiliary problem:
\begin{equation}\label{eq:bh}
\Delta^2 \overline{u}  = \delta_h  \quad \text{in }  \Omega, \quad
\overline{u}=0 \quad \text{and} \quad \partial_\mathbf{n} \overline{u}  =0 \quad \text{on }  \partial \Omega.
\end{equation}
The $C^0$ interior penalty method for problem \eqref{eq:bh} is to find $\overline{u}_h \in V_{h,0}^m$ such that 
\begin{align}
A_h(\overline{u}_h,v_h) = v(\mathbf{x}_0) = \int_{\Omega} \delta_h v_h\,\mathrm{d}\mathbf{x} \quad \forall v_h \in V_{h,0}^m.\label{discrete form2}
\end{align}
The well-posedness of the scheme \eqref{discrete form2} follows from the Lax-Milgram theorem.
Let $u$ and $\overline{u}$ be exact solutions for the original problem \eqref{eq:bh1} and auxiliary problem \eqref{eq:bh}, respectively. Let $\overline{u}_h$ be corresponding numerical solutions of \eqref{discrete form2}. The error estimate can be decomposed into two parts using the triangular inequality
\begin{align}\label{error}
|u - \overline{u}_h|_{H^2(\Omega)} \le |u - \overline{u}|_{H^2(\Omega)} + |\overline{u} - \overline{u}_h|_{H^2(\Omega)}.
\end{align}
The first term represents the regularization error, while the second term represents the discretization error. We estimate the total error by summing the independent contributions from each part.

Based on the regularity of the solution to \eqref{eq:bh}, the discretization error can be estimated as \cite{Scott73, BS2005}
\begin{align}\label{solprojerr1}
|\overline{u} - \overline{u}_h|_{H^2(\Omega)} \le Ch^{\min\{1, \alpha\}}.
\end{align}
where $\alpha< \alpha_0$ with $\alpha_0$ given in \eqref{alpha0}.

The error estimate \eqref{error} will be dominated by the regularization error. Referring to \cite{Scott73}, the following projection error bound holds
\begin{align}\label{ineq2}
|\delta_{\mathbf{x}_0} - \delta_h|_{H^{-r}(\Omega)} \le C h^{r-1}, \qquad 1 < r \le m,
\end{align}
where $m$ is the degree of the polynomial. Using the elliptic regularity theory and \eqref{ineq2} yield
\begin{align}\label{solprojerr2}
|u - \overline{u}|_{H^2(\Omega)} \le C | \delta_{\mathbf{x}_0} - \delta_h|_{H^{-2}(\Omega)} \le C h_{K_0}.
\end{align}

Based on \eqref{error}, \eqref{solprojerr1}, and \eqref{solprojerr2}, we have the following result.
\begin{lem}\label{projerr}
Let $\mathcal{T}_h$ be the quasi-uniform triangulation with mesh size $h$, and let $u$ and $\overline{u}_h$ be solutions of \eqref{eq:bh1} and \eqref{discrete form2}, respectively. Then the following error estimate holds
\begin{align*}
|u - \overline{u}_h|_{H^2(\Omega)} \le C h^{\min\{1, \alpha \}}.
\end{align*}
\end{lem}
\begin{rem}
Two techniques are empolyed in the $C^0$ interior penalty method to solve problem \eqref{eq:bh1}: a direct method \eqref{discrete form} and a method using the projection technique \eqref{discrete form2}. By comparing the error estimates of the solutions obtained with the direct method (see \Cref{prioriest}) and with the projection technique (see \Cref{projerr}), we observe that solutions from both techniques exhibit the same convergence rate on quasi-uniform meshes.
\end{rem}

Based on the projection technique, we propose the second residual-based a posteriori error estimator for the $C^0$ interior penalty method solving problem \eqref{eq:bh1} as
\begin{align}
\xi(\overline{u}_h) = \left(h_{K_0}^2 + \sum_{K \in \mathcal{T}_h}\xi_K^2(\overline{u}_h)\right)^{1/2},\label{global indicator2}
\end{align}
the local indicator $\xi_K$ is given by
\begin{align}\label{lobal indicator2}
\xi_{K}(\overline{u}_h) = 
\begin{cases}
\left(h_{K_0}^2 +\overline{\xi}_{K_0}^2\right)^{1/2},&\qquad \text{if }K = K_0, \\
\overline{\xi}_K,&\qquad \text{if }K \in \Omega \backslash K_0,
\end{cases}
\end{align}
where
\begin{align}
\xi_K(\overline{u}_h) = \left(\xi_{1,K}^2 +\sum_{e \in \mathcal{E}_K \cap \mathcal{E}_h}\alpha_e\xi_{2,e}^2 +  \sum_{e \in \mathcal{E}_K \cap \mathcal{E}_I}\alpha_e\xi_{3,e}^2 + \sum_{e \in \mathcal{E}_K \cap \mathcal{E}_I}\alpha_e\xi_{4,e}^2\right)^{1/2},
\end{align}
with $\alpha_e = 1$ for $e \in \mathcal{E}_B$, $\alpha _e= 1/2$ for $e \in \mathcal{E}_I$, and
\begin{equation*}
\begin{aligned}
\xi_{1,K} &= h_K^2\|\delta_h-\Delta^2 \overline{u}_h\|_{L^2(K)},\qquad\qquad\qquad
\xi_{2,e} = \beta h_e^{-1/2}\|[\![\nabla \overline{u}_h]\!]\|_{L^2(e)},\\
\xi_{3,e} &= h_e^{1/2}\|[\![\Delta \overline{u}_h]\!]\|_{L^2(e)},\qquad\qquad\qquad\qquad
\xi_{4,e} = h_e^{3/2}\|[\![\nabla \Delta \overline{u}_h]\!]\|_{L^2(e)}.
\end{aligned}
\end{equation*}
The corresponding global upper and local lower bounds are given as follows.
\begin{thm}[Reliability]
Let $u$ and $\overline{u}$ be exact solutions for the original problem \eqref{eq:bh1} and auxiliary problem \eqref{eq:bh}, respectively. And $\overline{u}_h \in V_{h,0}^m$ is the solution of \eqref{discrete form2}. Then the residual-based a posteriori error estimator $\xi$ satisfies the global bound
\begin{align}\label{theq3}
|||u-\overline{u}_h||| \le C \xi.
\end{align}
\end{thm}
\begin{proof}
By triangular inequality, the left hand side of \eqref{theq3} is bounded from above by
\begin{align*}
|||u-\overline{u}_h||| \le |||u-\overline{u}||| + |||\overline{u}-\overline{u}_h|||.
\end{align*}
Using the elliptic regularity bound, Lemma \ref{Weak continuity} and \eqref{solprojerr2}, the first term can be estimated as follows
\begin{align}\label{th3 ph0}
|||u-\overline{u}|||^2 = \sum_{K \in \mathcal{T}_h}|u-\overline{u}|^2_{H^2(K)} + \sum_{e \in \mathcal{E}_h}\frac{\beta}{h_e}\|[\![\nabla (u-\overline{u})]\!]\|_{L^2(e)}^2= |u-\overline{u}|_{H^2(\Omega)}^2 \le Ch_{K_0}^2.
\end{align}
Similar to the proof of Theorem \ref{th1}, to prove $|||\overline{u}-\overline{u}_h||| \le C \xi$, it is sufficient to verify that 
\begin{align}\label{th3 ph1}
\left(\sum_{K \in \mathcal{T}_h}|\overline{u}-\chi|^2_{H^2(K)}\right)^{1/2}\le \sup_{\phi \in H^2_0(\Omega)\backslash{\{0\}}}\frac{a(\overline{u}-\chi,\phi)}{|\phi|_{H^2(\Omega)}}\le C \xi,
\end{align}
where $\chi = E_h\overline{u}_h \in H_0^2(\Omega)$. Let $\phi_I \in V_{h,0}^m$ be the continuous interpolation polynomial of $\phi$. Then,
\begin{align}\label{th3 ph2}
a(\overline{u}-\chi,\phi)=\sum_{K \in \mathcal{T}_h}\int_K(\delta_h-\Delta^2\overline{u}_h)(\phi-\phi_I)\,\mathrm{d}\mathbf{x} + L_h,
\end{align}
where
\begin{align*}
L_h &= a_h(\overline{u}_h-\chi,\phi)+ \sum_{e \in \mathcal{E}_I}\int_e[\![\nabla(\Delta \overline{u}_h)]\!](\phi-\phi_I)\,\mathrm{d}s- \sum_{e \in \mathcal{E}_I}\int_e[\![\Delta \overline{u}_h]\!]\cdot\{\!\{\nabla(\phi-\phi_I)\}\!\}\,\mathrm{d}s\nonumber\\
 &\quad - \sum_{e \in \mathcal{E}_h}\int_{e} [\![\nabla \overline{u}_h]\!]\,\{\!\{\Delta \phi_I\}\!\} \,\mathrm{d}s + \sum_{e \in \mathcal{E}_h}\frac{\beta}{h_e}\int_{e}[\![\nabla \overline{u}_h]\!]\,[\![\nabla \phi_I]\!]\,\mathrm{d}s
\end{align*}
satisfies
\begin{align}\label{th3 ph3}
L_h \le C\left(\sum_{e \in \mathcal{E}_h}\xi_{2,e}^2 + \sum_{e \in \mathcal{E}_I}\xi_{3,e}^2 + \sum_{e \in \mathcal{E}_I}\xi_{4,e}^2\right)^{1/2}|\phi|_{H^2(\Omega)}.
\end{align}
By the Cauchy-Schwarz inequality and Lemma \ref{interpolation estimate}, the first term in \eqref{th3 ph2} follows
\begin{align}
& \sum_{K \in \mathcal{T}_h}\int_K\left(\delta_h-\Delta^2\overline{u}_h\right)(\phi-\phi_I)\,\mathrm{d}\mathbf{x}  \le \sum_{K \in \mathcal{T}_h}\|\delta_h-\Delta^2\overline{u}_h\|_{L^2(K)}\|\phi-\phi_I\|_{L^2(K)} \nonumber\\
&\quad \le C  \sum_{K \in \mathcal{T}_h}h_K^2\|\delta_h-\Delta^2\overline{u}_h\|_{L^2(K)}|\phi|_{H^2(K)} \le C\left(\sum_{K \in \mathcal{T}_h}\xi_{1,K}^2\right)^{1/2}|\phi|_{H^2(\Omega)},
\end{align}
which together with \eqref{th3 ph0} and \eqref{th3 ph3} yields the conclusion.
\end{proof}
\begin{thm}[Efficiency]\label{Efficiency2}
For the local indicator $\xi_{K}$ defined in \eqref{lobal indicator2}, there exists a positive constant $C$ independent of the mesh size such that
\begin{align}\label{theq4}
\xi_{K} \le 
\begin{cases}
C\left(h_{K_0} + |||u-\overline{u}_h|||_{\omega_{K_0}} \right), \qquad &K=K_0,\\
C|||u-\overline{u}_h|||_{\omega_{K}} , \qquad &K \in \mathcal{T}_h \backslash K_0.
\end{cases}
\end{align}
\end{thm}
\begin{proof}
Let's first show the element residual term $\xi_{1,{K}} = h^{2}_K\|\delta_h -\Delta^2 \overline{u}_h\|_{L^2(K)}$ satisfies the estimate \eqref{theq4}.
For $\forall v \in V_{h,0}^m$, using integration by part gives
\begin{align}\label{th4 pf2}
\int_{K}(\delta_h-\Delta^2 \overline{u}_h)v\,\mathrm{d}\mathbf{x} &= \int_{K}  (\delta_h - \delta_{\mathbf{x}_0}) v\,\mathrm{d}\mathbf{x} + \int_{K}  \Delta^2(u-\overline{u}_h) v\,\mathrm{d}\mathbf{x} \nonumber\\
& =  \int_{K}  (\delta_h - \delta_{\mathbf{x}_0}) v\,\mathrm{d}\mathbf{x} +\int_{K}\Delta (u-\overline{u}_h)\Delta v\,\mathrm{d}\mathbf{x} \nonumber\\
&\quad- \int_{\partial {K}}\Delta (u-\overline{u}_h) \mathbf{n}\cdot \nabla v\,\mathrm{d}s + \int_{\partial {K}}\nabla \Delta (u-\overline{u}_h) \cdot \mathbf{n}v\,\mathrm{d}s.
\end{align}
We set $v|_{K} = (\delta_h -\Delta^2 \overline{u}_h)b_{K}^2$ in \eqref{th4 pf2}, where $v \in H^2_0(\Omega)\cap H^2_0(K)$ and satisfies $v=0$ in $\Omega \backslash K$. Noticing that $v|_e = \nabla v|_e =0$. Consequently, \eqref{th4 pf2} yields
\begin{align*}
\int_{K}(\delta_h-\Delta^2 \overline{u}_h)v\,\mathrm{d}\mathbf{x} =\int_{K}  (\delta_h - \delta_{\mathbf{x}_0}) v\,\mathrm{d}\mathbf{x} +\int_{K}\Delta (u-\overline{u}_h)\Delta v\,\mathrm{d}\mathbf{x} .
\end{align*}
By Lemma \ref{bK properties}, it holds
\begin{align*}
\|v\|_{L^2({K})}\le
C \|\delta_h-\Delta^2 \overline{u}_h\|_{L^2({K})}.
\end{align*}
We prove the first case in \eqref{theq4}, i.e., $K = K_0$. By Cauchy-Schwarz inequality, Lemma \ref{Inverse inequality} and \eqref{solprojerr2},
\begin{align}\label{Lapuhv2}
\int_{K_0}(\delta_h -\Delta^2 \overline{u}_h)v\,\mathrm{d}\mathbf{x} 
&\le \left(|\delta_h - \delta_{\mathbf{x}_0}|_{H^{-2}(K_0)}+ \|\Delta (u-\overline{u}_h)\|_{L^2(K_0)}\right)|v|_{H^2({K_0})} \notag\\
&\le C h^{-2}_{K_0}\left(h_{K_0} + |u-\overline{u}_h|_{H^2(\omega_{K_0})}\right)\|v\|_{L^2({K_0})}.
\end{align}
By \eqref{bK property 2} and \eqref{Lapuhv2}, it holds
\begin{align*}
\|\delta_h -\Delta^2 \overline{u}_h\|_{L^2({K_0})}^2 & \le C\int_{K_0}(\delta_h -\Delta^2\overline{u}_h)^2b_{K_0}^2\,\mathrm{d}\mathbf{x} = C\int_{K_0}(\delta_h-\Delta^2\overline{u}_h)v\,\mathrm{d}\mathbf{x}\nonumber\\
& \le Ch^{-2}_{K_0}\left(h_{K_0} + |u-\overline{u}_h|_{H^2(\omega_{K_0})}\right)\|\delta_h -\Delta^2 \overline{u}_h\|_{L^2({K_0})},
\end{align*}
which implies
\begin{align}\label{th4 pf3}
h^{2}_{K_0}\|\delta_h -\Delta^2 \overline{u}_h\|_{L^2({K_0})}\le C\left(h_{K_0} + |u-\overline{u}_h|_{H^2(\omega_{K_0})}\right).
\end{align}
Note that $(\delta_h - \delta_{\mathbf{x}_0})|_K = 0$ for any $K \in \mathcal{T}_h \backslash K_0$. Similar to the proof of the first case in \eqref{theq4}, it follows
\begin{align}\label{th4 pf4}
h^{2}_{K}\|\delta_h -\Delta^2 \overline{u}_h\|_{L^2({K})}\le C|u-\overline{u}_h|_{H^2(\omega_{K})}.
\end{align}
Verifying other terms in $\xi_K$ can refer to the proof of Theorem \ref{Efficiency}.
\end{proof}

\subsection{Adaptive finite element algorithm}


The adaptive finite element algorithm based on the residual-based a posteriori error estimator \eqref{global indicator} or \eqref{global indicator2} is summarized as follows.
\begin{breakablealgorithm}
\caption{The adaptive finite element algorithm. } 
\label{alg2} 
\begin{algorithmic}[1] 
\STATE Input: an initial mesh $\mathcal{T}^0_h$; a constant $0<\theta \le 1$; the maximum number of mesh refinements $n$.
\STATE Output: the numerical solution $u_h^n$ (\text{resp.} $\overline{u}_h^n$); a new refined mesh $\mathcal{T}^n_h$.
\STATE for $i=0$ to $n$ do\\
       \qquad Solve the discrete equation for the finite element solution $u_h^i$ (resp. $\overline{u}_h^i$) on $\mathcal{T}^i_h$;\\
       \qquad Computing the local and total error estimation $\eta_{K}^i(u_h^i)$ (resp. $\xi_{K}^i(\overline{u}_h^i)$) and $\eta^i(u_h^i)$ (\text{resp.} $\xi^i(\overline{u}_h^i)$);\\
       \qquad if $i< n$ then\\
       \qquad \qquad Select a subset $\widetilde{\mathcal{T}_h}^i \subset \mathcal{T}^i_h$ of marked elements to refined such that,
       \begin{align*}
       \left(\sum_{K\in \widetilde{\mathcal{T}_h}^i}{\eta_{K}^i(u_h^i)}^2\right)^{1/2}\ge \theta \eta^i(u_h^i),\quad \left(\text{resp.} \quad \left(\sum_{K\in \widetilde{\mathcal{T}_h}^i}{\xi_{K}^i(\overline{u}_h^i)}^2\right)^{1/2}\ge \theta \xi^i(\overline{u}_h^i) \right);
       \end{align*}
       \qquad \qquad Refine the each element $K \in \widetilde{\mathcal{T}_h}^i$ by longest edge bisection to obtain a new mesh $\mathcal{T}^{i+1}_h$.\\
       \qquad end if\\
       end for
\end{algorithmic}
\end{breakablealgorithm}

\section{Extensions}\label{sec4}
In this section, we extend our results to cover a broader class of fourth-order elliptic equations with various boundary conditions. 
Specifically, we generalize the biharmonic operator $\Delta^2$ in (\ref{eq:bh1}) to a more general fourth order operator $\mathcal{L}$ that includes low order terms
\begin{equation}\label{operator}
\mathcal{L}u = \Delta^2 u - \mu_1 \Delta u + \mu_2 u,
\end{equation}
where $\mu_\ell \geq 0$ (for $\ell=1,2$) are constants. 
Both residual-based a posteriori error estimators, \eqref{lobal indicator} and \eqref{lobal indicator2}, for problem \eqref{eq:bh1} can be extended to the new operator. We only present the first type residual-based a posteriori error estimator for simplicity. The proofs of the efficiency and reliability are expected to be similar to those provided in \Cref{sec3}. Therefore, we will focus solely on presenting the corresponding a posteriori error estimators and leave their performance to be verified numerically.

\subsection{Non-homogeneous Dirichlet boundary conditions}

Consider the following problem:
\begin{equation}\label{eq:bh4}
\Delta^2 u - \mu_1 \Delta u + \mu_2 u  = \delta_{\mathbf{x}_0} + f \quad \text{in }  \Omega, 
\quad
u =g \quad \text{and} \quad \partial_\mathbf{n} u  =g_N \quad \text{on }  \partial \Omega,
\end{equation}
where $g_D$ and $g_N$ are given functions on the boundary $\partial \Omega$, and $f$ is a given function in the domain $\Omega$. We assume that $g$, $g_N$, and $f$ are sufficiently smooth. This problem is a scalar analog of the variational problem for the strain gradient theory in elasticity and plasticity \cite{SKF1999} when the coefficients satisfy $\mu_1>0$ and $\mu_2 =0$. 

In \cite{BS2005}, a $C^0$ interior penalty method was proposed to solve the general fourth-order elliptic equation without the Dirac delta function term $\delta_{\mathbf{x}_0}$.
In the subsequent section, we extend and adapt this method to handle the presence of the Dirac delta function term $\delta_{\mathbf{x}_0}$ on the right-hand side of equation \eqref{eq:bh4}.

We define the space $V_h^m$ associated with the triangulation $\mathcal{T}_h$ by
\begin{align}
V_h^m =\{v_h \in H^1(\Omega) \cap C^0(\Omega):\,v_h|_K \in P_m(K),\,m\ge 2,\quad \forall K \in \mathcal{T}_h\}.
\end{align}
Denote the subspace of $V_h^m$ by
\begin{align}
V_{h,g}^m =\{v_h \in V_h^m\,:\,v_h|_{\partial \Omega}= \mathcal{I}_h g\},
\end{align}
where $\mathcal{I}_hg \in V_h^m$ is a polynomial approximation function of $g$ by the interpolation on the boundary.

The $C^0$ interior penalty method for \eqref{eq:bh4} is to find $u_h \in V_{h,g}^m$ such that
\begin{align}\label{discrete form3}
\sum_{K \in \mathcal{T}_h}&\left(\int_K \Delta u_h \Delta v_h \,\mathrm{d}\mathbf{x}+\mu_1\int_K \nabla u_h\cdot \nabla v_h \,\mathrm{d}\mathbf{x} + \mu_2\int_K  u_h v_h \,\mathrm{d}\mathbf{x}\right)+ \sum_{e \in \mathcal{E}_h}\frac{\beta}{h_e}\int_{e}[\![\nabla u_h]\!]\,[\![\nabla v_h]\!]\,\mathrm{d}s  \nonumber\\
&- \sum_{e \in \mathcal{E}_h}\Big(\int_{e} \{\!\{\Delta u_h\}\!\} \,[\![\nabla v_h]\!]\,\mathrm{d}s +  \int_{e}\{\!\{\Delta v_h\}\!\} \,[\![\nabla u_h]\!]\,\mathrm{d}s \Big)\nonumber\\
= &v(\mathbf{x}_0)+\sum_{K \in \mathcal{T}_h}\int_Kfv\,\mathrm{d}\mathbf{x} -  \sum_{e \in \mathcal{E}_B}\int_{e} g_N \,\Delta v\,\mathrm{d}s + \sum_{e \in \mathcal{E}_B}\frac{\beta}{h_e}\int_{e} g_N \,\nabla v \cdot \mathbf{n}\,\mathrm{d}s, \quad \forall v_h \in V_{h,0}^m,
\end{align}
where $V_{h,0}^m$ is defined as \eqref{space}.

Unlike the $C^0$ interior penalty method for problems with homogeneous Dirichlet boundary conditions described in \eqref{discrete form}, the modified method \eqref{discrete form3} incorporates additional boundary term. 
These additional terms account for the non-homogeneous nature of the boundary conditions.
The local error estimator on $K \in \mathcal{T}_h$ for $C^0$ interior penalty method \eqref{eq:bh4} is given by 
\begin{align}
\eta_{K}(u_h) = 
\begin{cases}
\left(h_{K_0}^2 +\widetilde{\eta}_{K_0}^2\right)^{1/2},&\qquad \text{if } K = K_0 \text{ and }\mathbf{x}_0 \notin \mathcal{N},\\
\widetilde{\eta}_K,&\qquad \text{otherwise},
\end{cases}\label{local indicator2}
\end{align}
where
\begin{align}
\widetilde{\eta}_K(u_h) &= \left(\widetilde{\eta}_{1,K}^2 +\sum_{e \in \mathcal{E}_K \cap \mathcal{E}_I}\gamma_e\eta_{2,e}^2 +  \sum_{e \in \mathcal{E}_K \cap \mathcal{E}_I}\alpha_e\eta_{3,e}^2 + \sum_{e \in \mathcal{E}_K \cap \mathcal{E}_I}\alpha_e\eta_{4,e}^2 \right.\notag \\
&\qquad \left.+ \sum_{e \in \mathcal{E}_K \cap \mathcal{E}_B}\gamma_e\eta_{5,e}^2 + \sum_{e \in \mathcal{E}_K \cap \mathcal{E}_B}\gamma_e\eta_{6,e}^2 \right)^{1/2},\label{lobal indicator part2}
\end{align}
with $\gamma_e = (1+\mu_1h_e^2 + \mu_2h_e^4)\alpha_e$, $\alpha_e = 1$ for $e \in \mathcal{E}_B$, $\alpha_e = 1/2$ for $e \in \mathcal{E}_I$, and
\begin{align}
&\widetilde{\eta}_{1,K}^2 = h_K^4\|f - \Delta^2u_h + \mu_1\Delta u_h - \mu_2u_h\|_{L^2(K)}^2,\\
&\eta_{5,e}^2 = \beta h^{-1/2}_e\|g_N-\nabla u_h\cdot \mathbf{n}\|_{L^2(e)}^2,\\
&\eta_{6,e}^2 = h^{-3/2}_e\|g-u_h\|_{L^2(e)}^2. \label{lobal 2}
\end{align}
and $\eta_{2,e},\eta_{3,e}$ and $\eta_{4,e}$ are given in \eqref{part 1}-\eqref{part 2}.


\subsection{Navier boundary conditions}
We also extend the results to the fourth order equation with Navier boundary conditions, 
\begin{equation}\label{eq:bh5}
\Delta^2 u - \mu_1 \Delta u + \mu_2 u  = \delta_{\mathbf{x}_0} + f \quad \text{in }  \Omega, 
\quad
u =g \quad \text{and} \quad \Delta u  =g_B \quad \text{on }  \partial \Omega,
\end{equation}
where $g$ and $g_B$ are given functions on the boundary $\partial \Omega$, and $f$ is a given function within the domain $\Omega$. We also assume that $g$, $g_B$, and $f$ are sufficiently smooth. 
The problem \eqref{eq:bh5} models a simply supported plate problem \cite{S1992}. When the Dirac delta function term $\delta_{\mathbf{x}_0}$ vanishes, \cite{BN2011} studied a symmetric $C^0$ interior penalty method for a fourth-order singular perturbation elliptic problem with these types of boundary conditions in two dimensions on polygonal domains.
The $C^0$ interior penalty method for \eqref{eq:bh5} is to find $u_h \in V_{h,g}^m$ such that
\begin{align}
\sum_{K \in \mathcal{T}_h}&\left(\int_K \Delta u_h \Delta v_h \,\mathrm{d}\mathbf{x}+\mu_1\int_K \nabla u_h\cdot \nabla v_h \,\mathrm{d}\mathbf{x} + \mu_2\int_K  u_h v_h \,\mathrm{d}\mathbf{x}\right) + \sum_{e \in \mathcal{E}_I}\frac{\beta}{h_e}\int_{e}[\![\nabla u_h]\!]\,[\![\nabla v_h]\!]\,\mathrm{d}s \nonumber\\
&-  \sum_{e \in \mathcal{E}_I}\Big(\int_{e} \{\!\{\Delta u_h\}\!\} \,[\![\nabla v_h]\!]\,\mathrm{d}s +  \int_{e}\{\!\{\Delta v_h\}\!\} \,[\![\nabla u_h]\!]\,\mathrm{d}s \Big)\nonumber\\
=& v(\mathbf{x}_0)+ \sum_{K \in \mathcal{T}_h}\int_Kfv\,\mathrm{d}\mathbf{x}+ \sum_{e \in \mathcal{E}_B}\int_{e} g_B \,\nabla v \cdot \mathbf{n}\,\mathrm{d}s.\qquad \forall v_h \in V_{h,0}^m,\label{discrete form4}
\end{align}
where $V_h^m$ is defined as \eqref{space}.
The local error estimator for the $C^0$ interior penalty method \eqref{discrete form4} is given by
\begin{align}\label{BiNavier}
\eta_{K}(u_h) = 
\begin{cases}
\left(h_{K_0}^2 +\zeta_{K_0}^2\right)^{1/2},&\qquad \text{if } K =K_0 \text{ and } \mathbf{x}_0 \notin \mathcal{N},\\ 
\zeta_K,&\qquad \text{otherwise},
\end{cases}
\end{align}
where
\begin{align}
\zeta_K(u_h) &= \left(\widetilde{\eta}_{1,K}^2 +\sum_{e \in \mathcal{E}_K \cap \mathcal{E}_h}\gamma_e\eta_{2,e}^2 +  \sum_{e \in \mathcal{E}_K \cap \mathcal{E}_I}\alpha_e\eta_{3,e}^2 + \sum_{e \in \mathcal{E}_K \cap \mathcal{E}_I}\alpha_e\eta_{4,e}^2 \right.\notag \\
&\qquad \left.+ \sum_{e \in \mathcal{E}_K \cap \mathcal{E}_B}\zeta_{5,e}^2+ \sum_{e \in \mathcal{E}_K \cap \mathcal{E}_B}\gamma_e\eta_{6,e}^2\right)^{1/2},\label{lobal indicator part3}
\end{align}
the element residual $\widetilde{\eta}_{1,K}$ and the edge residuals $\eta_{2,e},\eta_{3,e}$, $\eta_{4,e}$, $\eta_{6,e}$ are defined as \eqref{lobal indicator part2}. The boundary residual term $\zeta_{5,e}$ is defined by 
\begin{align}\label{zeta5e}
\zeta_{5,e}^2 = h^{1/2}_e||g_B-\Delta u_h||_{L^2(e)}^2.
\end{align}

\subsection{Homogeneous Neumann boundary conditions}
For fourth-order elliptic equations with homogeneous Neumann boundary conditions, we consider the following model
\begin{equation}\label{eq:bh6}
\Delta^2 u - \mu_1 \Delta u + \mu_2 u  = \delta_{\mathbf{x}_1}-\delta_{\mathbf{x}_0} \quad \text{in }  \Omega, 
\quad
 \partial_{\mathbf{n}} u  = 0 \quad \text{and} \quad  \partial_{\mathbf{n}} (\Delta u)  = 0 \quad \text{on }  \partial \Omega.
\end{equation}
where $\mathbf{x}_0,\,\mathbf{x}_1$ 
represent two distinct points strictly contained in the domain $\Omega$. These boundary-value problems can arise in the Cahn-Hilliard model, which describes phase-separation phenomena \cite{CH1958}.
\cite{BGGS2012} proposed a quadratic $C^0$ interior penalty method for fourth-order boundary value problems with similar types of boundary conditions, assuming the right-hand side function $f\in  L^2(\Omega)$. Under these conditions, a unique solution $u$ can be found, satisfying $u \in H^{2+\alpha}(\Omega)$ for $\alpha \in (0,2]$.



It is straightforward to validate the solvability condition:
\begin{align*}
\int_{\Omega} \delta_{\mathbf{x}_1}-\delta_{\mathbf{x}_0}\,\mathrm{d}x = 0.
\end{align*}
To obtain a unique solution, a common approach is to impose an additional constraint
\begin{align*}
\int_{\Omega} u\,\mathrm{d}x = 0.
\end{align*}
We define a subspace of $H^1(\Omega)$ by
\begin{align*}
V =\{v_h \in H^1(\Omega):\partial_{\mathbf{n}} v_h = 0,\,\int_{\Omega} v_h\,\mathrm{d}x = 0\}.
\end{align*}
The finite element space $X_h^m$ associated with the triangulation $\mathcal{T}_h$ is defined as
\begin{align}
X_h^m =\{v_h \in V \cap C^0(\Omega):\,v_h|_K \in P_m(K),\,m\ge 2,\quad \forall K \in \mathcal{T}_h\}.
\end{align}
The $C^0$ interior penalty method for \eqref{eq:bh6} is then to find $u_h \in X_h^m$ such that
\begin{align}\label{C0IPHNM}
\sum_{K \in \mathcal{T}_h}&\left(\int_K \Delta u_h \Delta v_h \,\mathrm{d}\mathbf{x}+\mu_1\int_K \nabla u_h\cdot \nabla v_h \,\mathrm{d}\mathbf{x} + \mu_2\int_K  u_h v_h \,\mathrm{d}\mathbf{x}\right)+ \sum_{e \in \mathcal{E}_h}\frac{\beta}{h_e}\int_{e}[\![\nabla u_h]\!]\,[\![\nabla v_h]\!]\,\mathrm{d}s  \nonumber\\
&- \sum_{e \in \mathcal{E}_h}\Big(\int_{e} \{\!\{\Delta u_h\}\!\} \,[\![\nabla v_h]\!]\,\mathrm{d}s +  \int_{e}\{\!\{\Delta v_h\}\!\} \,[\![\nabla u_h]\!]\,\mathrm{d}s \Big) = v(\mathbf{x}_1) -v(\mathbf{x}_0).\quad \forall v_h \in X_h^m,
\end{align}
Note that the bilinear in \eqref{C0IPHNM} h is identical to that in \eqref{discrete form3}, with the exception that the finite element space now includes only functions with a zero mean value.

Denote $K_1 \in \mathcal{T}_h$ by one element such that the singular point $\mathbf{x}_1 \in \overline{K}_1$, where $\overline{K}_1$ is the closure of $K_1$.The diameter of $K_1$ is denoted by $h_{K_1}$. The local error estimator for problem \eqref{eq:bh5} is given by 
\begin{align}
\eta_{K}(u_h) = 
\begin{cases}
\left(h_{K_0}^2 +{\chi}_{K_0}^2\right)^{1/2},&\qquad \text{if }K = K_0,\,\text{ and } \mathbf{x}_0 \notin \mathcal{N}, \\
\left(h_{K_1}^2 +{\chi}_{K_1}^2\right)^{1/2},&\qquad \text{if }K = K_1,\,\text{ and } \mathbf{x}_1 \notin \mathcal{N}, \\
{\chi}_K,&\qquad \text{otherwise}.
\end{cases}
\end{align}
where
\begin{align}
{\chi}_K(u_h) &= \left(\widetilde{\eta}_{1,K}^2 + \sum_{e \in \mathcal{E}_K \cap \mathcal{E}_h}\gamma_e\eta_{2,e}^2 +  \sum_{e \in \mathcal{E}_K \cap \mathcal{E}_I}\alpha_e\eta_{3,e}^2 + \sum_{e \in \mathcal{E}_K \cap \mathcal{E}_I}\alpha_e\eta_{4,e}^2\right).\label{lobal indicator part4}
\end{align}
The definition of ${\chi}_K(u_h)$ is same as $\zeta_K(u_h)$ in \eqref{lobal indicator part3}, except for discarding the the boundary residual term $\zeta_{5,e}$.

\section{Numerical examples}\label{sec5}

In this section, we present numerical test results to verify the accuracy of the $C^0$ interior penalty method and demonstrate the robustness of the proposed residual-type a posteriori estimators. 
If the exact solution $u$ is given, the convergence rate is calculated by
\begin{eqnarray}\label{rate1}
\mathcal R=\log_2 \frac{|u-u_h^{j-1}|_{H^2(\Omega)}}{|u-u_h^j|_{H^2(\Omega)}},
\end{eqnarray}
where $u_h^j$ is the $C^0$ finite element solution on the mesh $\mathcal{T}^j_h$ obtained after $j$th refinements of the initial triangulation $\mathcal{T}_h^0$.  
When the exact solution is unavailable or difficult to obtain, we instead use the following numerical convergence rate
\begin{eqnarray}\label{rate}
\mathcal R=\log_2\frac{|u_h^j-u_h^{j-1}|_{H^2(\Omega)}}{|u_h^{j+1}-u_h^j|_{H^2(\Omega)}}.
\end{eqnarray}

Due to the lack of regularity, the $C^0$ interior penalty method with high-degree polynomial approximations on quasi-uniform triangular meshes may not achieve optimal convergence rates. To address this, we apply an adaptive $C^0$ interior penalty method to improve the convergence order.

The convergence rate of the a posteriori error estimator $\eta$ (resp. $\xi$)  for $P_m$ polynomials with $m \geq 2$ is quasi-optimal if
$$
\eta \approx N^{-0.5(m-1)} \quad \text{(resp. } \xi \approx N^{-0.5(m-1)} \text{)}.
$$
Here and in what follows, we abuse the notation $N$ to represent the total number of degrees of freedom. 

\begin{example}[L-shape domain]\label{exam2}
We consider problem \eqref{eq:bh1} with homogeneous clamped boundary
conditions in an L-shaped domain $\Omega = (-2\pi,2\pi)^2 \backslash [0,2\pi)\times (-2\pi,0]$ with a largest interior interior angle $\omega = 3\pi/2$. The Dirac point $\mathbf{x}_0 = (-\pi,\pi)$, then the solution show singularities at two points $(0,0)$ and $\mathbf{x}_0$. We start with an initial mesh $\mathcal{T}_h^0$ as Figure \ref{fig:exam2 initial mesh}(a).

\begin{table}[H]
\centering
\vspace{-3mm}
\caption{Example \ref{exam2}: Convergence rate of numerical solution in uniform meshes.}
\begin{tabular}{|c|c|c|c|c||c|c|c|c||c|c|c|c|c|c|c|c|}
\hline
&  \multicolumn{3}{c}{$\qquad \quad P_2$} & &\multicolumn{3}{c}{$\qquad  \quad P_3$} & &\multicolumn{3}{c}{$\qquad  \quad P_4$}&\\
\hline
j &    4 &   5 &   6 &7 &    3 &    4 &   5 &   6  &     3 &    4 &   5 &   6\\
\hline
$\mathcal{R}$&  0.88 &  0.85 &  0.80  &  0.73 &  1.00   &  0.73 &  0.60  &  0.56 &  0.74 &  0.58 &  0.56 &  0.55\\
\hline
\end{tabular}\label{table:exam2}
\end{table}

Table \ref{table:exam2} shows the $H^2$ convergence history using the $P_m$-based $C^0$ interior penalty method on quasi-uniform meshes, with $m = 2,\,3,\,4$. From the results, we observe that the convergence rates are $\mathcal{R} < 1$ on coarse meshes, and $\mathcal{R} \approx 0.5445$ on sufficiently refined meshes. This indicates that the singularity in the solution is primarily influenced by the reentrant corner of the polygonal domain when $\omega>\pi$. The results in Table \ref{table:exam2} align with the theoretical expectations outlined in Lemma \ref{udecompthm}.

We then assess the efficiency of the a posteriori error estimators. Meshes generated by $\eta$ in \eqref{global indicator} and $\xi$ in \eqref{global indicator2} are shown Figure \ref{fig:exam2 test1 meshes} and Figure \ref{fig:exam2 test2 meshes}, respectively. It is evident that the error estimators effectively guide mesh refinements around the points $(0,0)$ and $(-\pi,\pi)$, where the solution shows singularities. As shown in Figure \ref{fig:exam2 errors}, the slopes for the estimators are close to $-0.5(m-1)$ when there are sufficient grid points, indicating optimal decay of the error with respect to the number of unknowns. The contours of the corresponding numerical solutions based on  $\eta$ and $\xi$ are shown in Figure \ref{fig:exam2 initial mesh}(b)-(c), they are very similar.
\begin{figure}
\centering
\subfigure[Initial mesh]
{\includegraphics[width=0.29\textwidth]{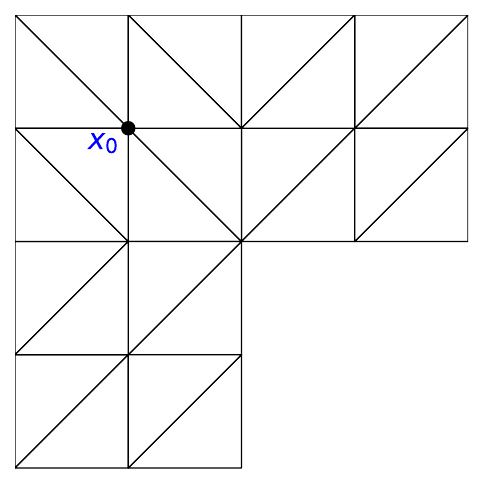}}
\subfigure[Solution obtained with $\eta$]
{\includegraphics[width=0.34\textwidth]{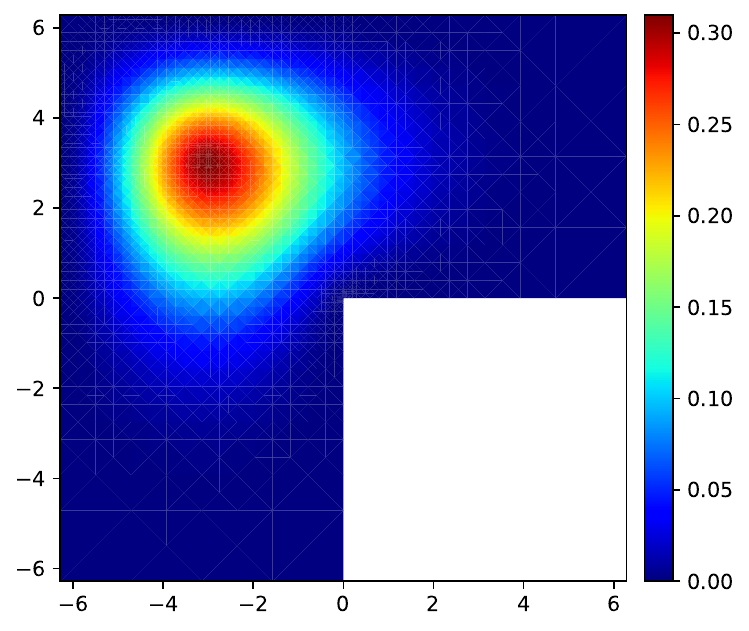}}
\subfigure[Solution obtained with $\xi$]
{\includegraphics[width=0.34\textwidth]{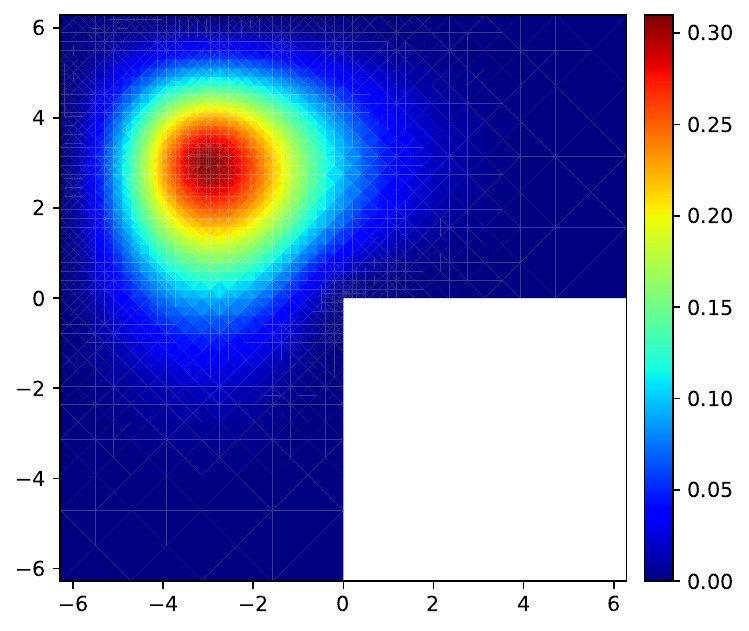}}
\vspace{-2mm}
\caption{Example \ref{exam2}: initial mesh and adaptive numerical solution.}\label{fig:exam2 initial mesh}
\end{figure}
\begin{figure}
\centering
\subfigure[$P_2$]
{\includegraphics[width=0.3\textwidth]{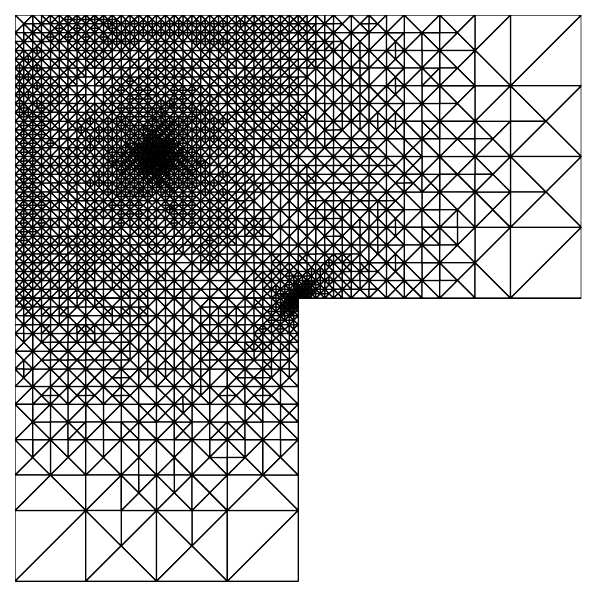}}\hspace{3mm}
\subfigure[$P_3$]
{\includegraphics[width=0.3\textwidth]{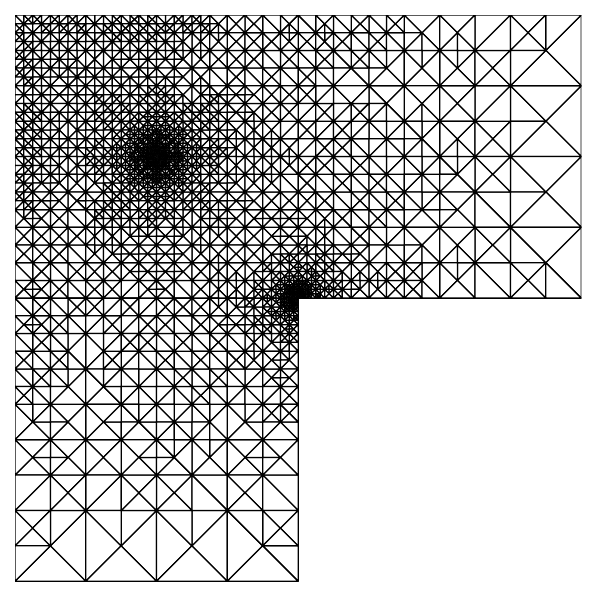}}\hspace{3mm}
\subfigure[$P_4$]
{\includegraphics[width=0.3\textwidth]{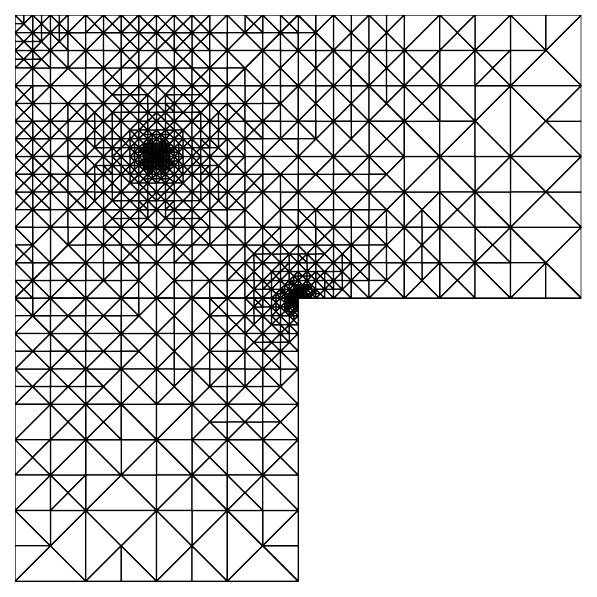}}
\vspace{-2mm}
\caption{Example \ref{exam2}: adaptive meshes generated by $\eta$.}\label{fig:exam2 test1 meshes}
\end{figure}

\begin{figure}
\centering
\subfigure[$P_2$]
{\includegraphics[width=0.3\textwidth]{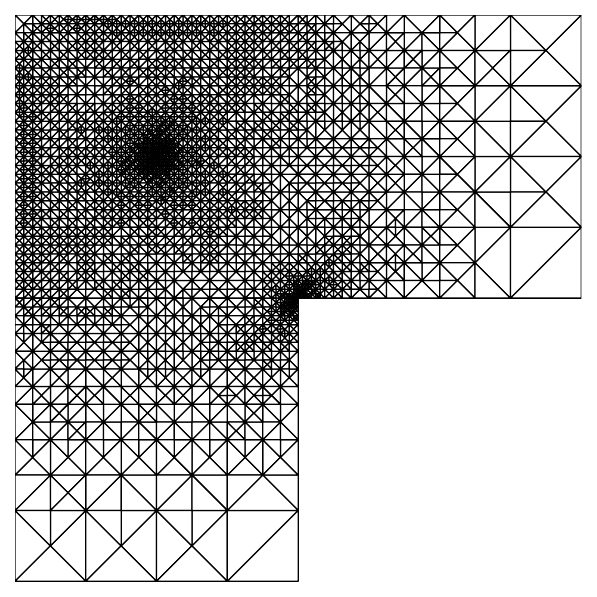}}\hspace{3mm}
\subfigure[$P_3$]
{\includegraphics[width=0.3\textwidth]{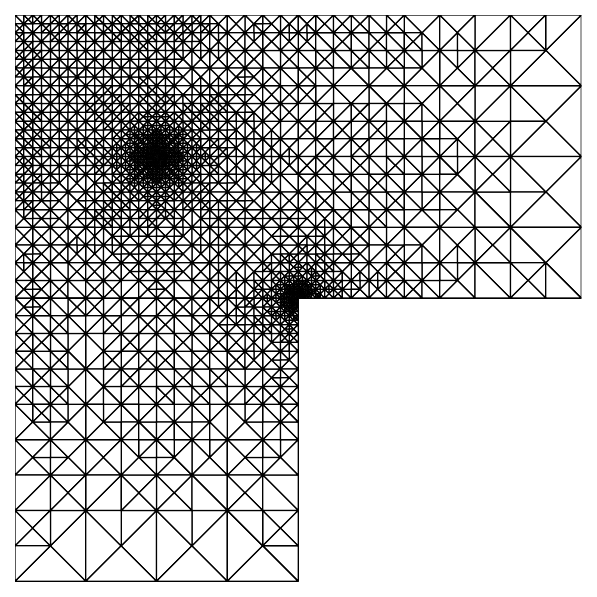}}\hspace{3mm}
\subfigure[$P_4$]
{\includegraphics[width=0.3\textwidth]{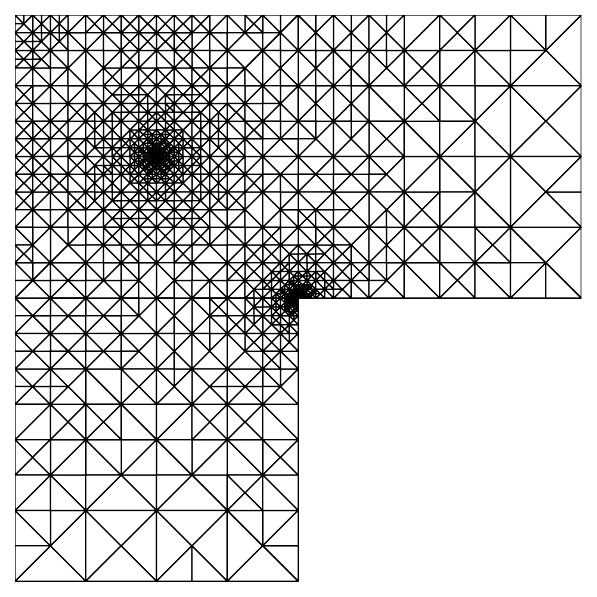}}
\vspace{-2mm}
\caption{Example \ref{exam2}: adaptive meshes generated by $\xi$.}\label{fig:exam2 test2 meshes}
\end{figure}

\begin{figure}
\centering
\subfigure[$P_2$]
{\includegraphics[width=0.32\textwidth]{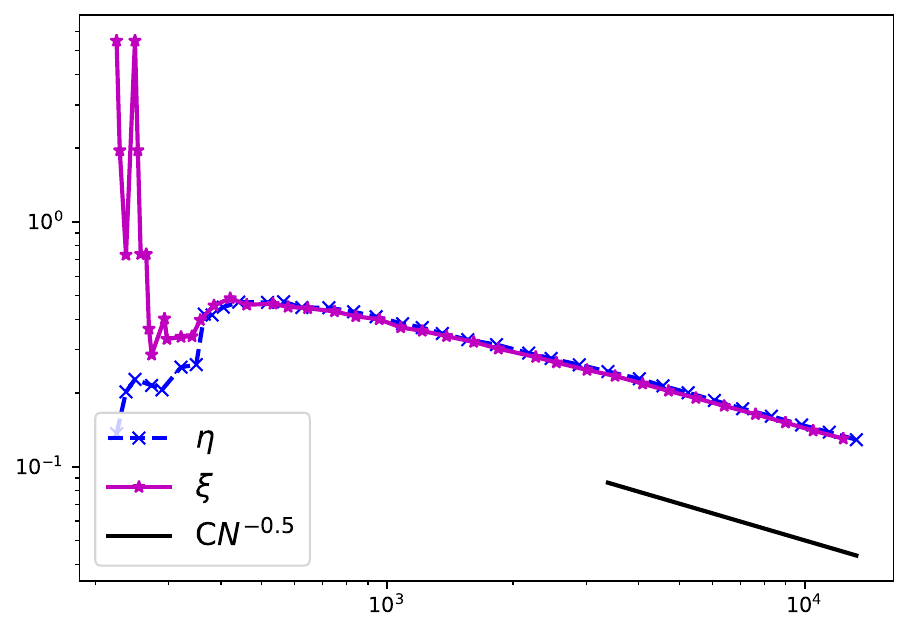}}
\subfigure[$P_3$]
{\includegraphics[width=0.32\textwidth]{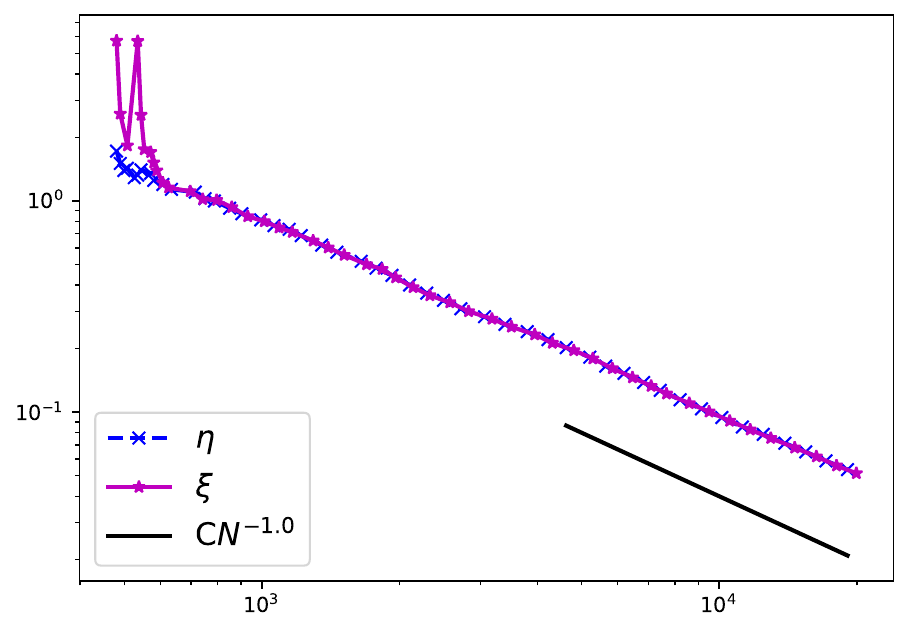}}
\subfigure[$P_4$]
{\includegraphics[width=0.32\textwidth]{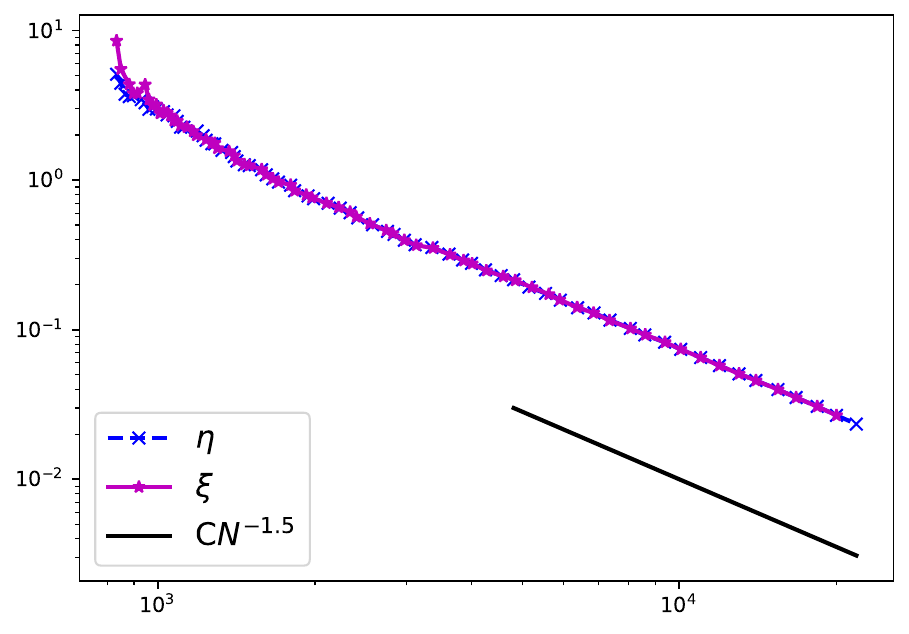}}
\vspace{-2mm}
\caption{Example \ref{exam2}: error estimators.}\label{fig:exam2 errors}
\end{figure}
\end{example}

\begin{example}[Non-homogeneous boundary]\label{exam1}
In this example, we consider a more general biharmonic equation 
\begin{align}\label{equation}
\Delta^2 u - \mu_1 \Delta u + \mu_2 u  = \delta_{\mathbf{x}_0} + f,
\end{align}
where the location of ${\mathbf{x}_0}$ in $\mathcal{T}_h$ are of three different types as shown in Case 1-3. The initial meshes of Cases 1-3 are reported in Figure \ref{fig:exam1 initialmeshes}.
\begin{itemize}
 \item [Case 1: ]{${\mathbf{x}_0} = (0,0)$ is a node of the triangulations.}
\item [Case 2: ]{${\mathbf{x}_0} = (-\sqrt{7},-\pi)$ belongs to an inter edge $e \in \mathcal{E}_I$.}
\item [Case 3: ]{${\mathbf{x}_0} = (\sqrt{5},\sqrt{8})$ is contained by one element $K \in \mathcal{T}_h$.}
\end{itemize}
\begin{figure}
\centering
\subfigure[Case1,\,$\mathbf{x}_0 = (0,0)$]
{\includegraphics[width=0.3\textwidth]{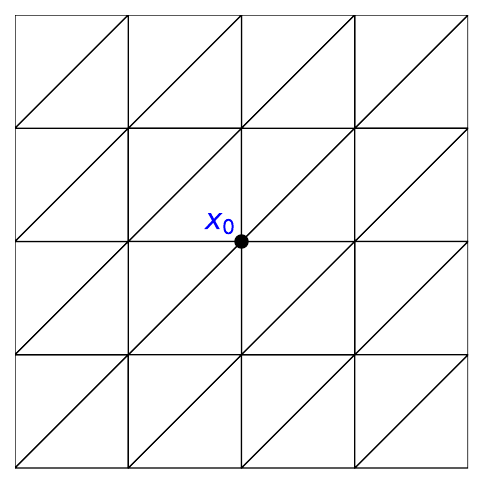}}\hspace{3mm}
\subfigure[Case2,\,$\mathbf{x}_0 = (-\sqrt{7},-\pi)$]
{\includegraphics[width=0.3\textwidth]{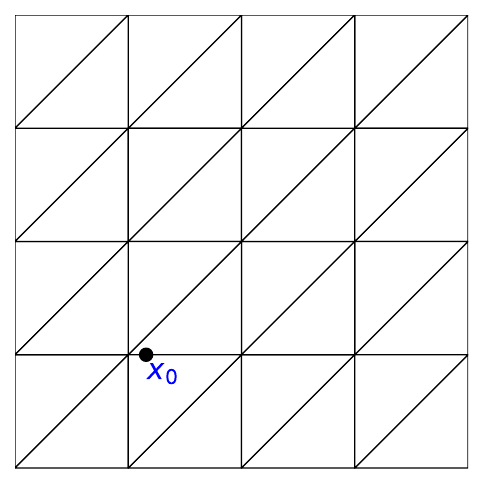}}\hspace{3mm}
\subfigure[Case3,\,$\mathbf{x}_0 = (\sqrt{5},\sqrt{8})$]
{\includegraphics[width=0.3\textwidth]{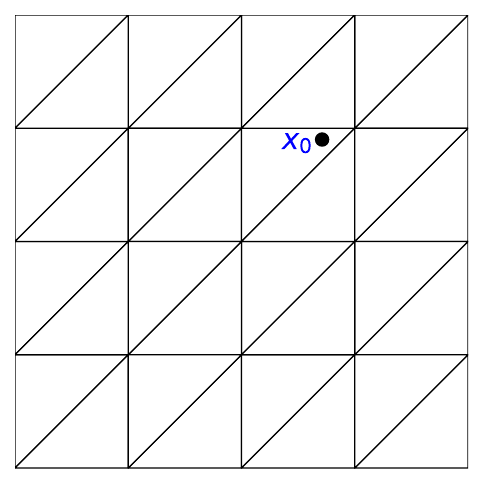}}
\vspace{-2mm}
\caption{Example \ref{exam1}: initial meshes.}\label{fig:exam1 initialmeshes}
\end{figure}
We take the function 
\begin{align}
f(x,y) =\mu_2 \frac{|\mathbf{x} - \mathbf{x}_0|^2}{8\pi}\ln |\mathbf{x} - \mathbf{x}_0| - \mu_1\left(\frac{\ln |\mathbf{x} - \mathbf{x}_0|}{2\pi} + \frac{1}{2\pi}\right),
\end{align}
then the exact solution of equation \eqref{equation} is given by
\begin{align}
u(x,y) = \frac{|\mathbf{x} - \mathbf{x}_0|^2}{8\pi}\ln |\mathbf{x} - \mathbf{x}_0|.\label{solution}
\end{align}
We consider \eqref{equation} with two different boundary conditions and the parameters in $u(x,y)$ are taken as $\mu_1 = \mu_2 = 1$. 

\noindent\textbf{Test 1.} We first consider non-homogeneous campled boundary conditions.
The convergence rates of the $C^0$ interior penalty method solutions based on $P_2,\,P_3$ and $P_4$ polynomials for Case 1-3 in quasi-uniform meshes are shown in Table \ref{table:exam1_text1}. The convergence rates are approximately $\mathcal{R} \approx 1$. While the convergence rates for $P_2$ are quasi-optimal, the rates for $P_3,\,P_4$ only achieve suboptimal convergence.
This is due to $u \in H^{3-\epsilon}(\Omega)$ for $\omega = \frac{\pi}{2}$. Given the low regularity of the solution $u$, these results are the best that can be achieved with quasi-uniform meshes.

To address this, we apply the adaptive $C^0$ interior penalty method based on the residual-based a posteriori error estimator $\eta_K$ in \eqref{local indicator2}. The corresponding numerical solutions of the adaptive algorithm using the error estimator $\eta$ are presented in \Cref{fig:exam1 test1 solution}. \Cref{fig:exam1 test1 meshes p3} and \Cref{fig:exam1 test1 meshes p4} show the adaptive meshes of $P_3,\,P_4$ approximations, respectively.  The error estimator effectively guides mesh refinements, particularly around the point $\mathbf{x}_0$. The convergence rates of the error estimator $\eta$ based on $P_3,\,P_4$ polynomials are illustrated in Figures \ref{fig:exam1 test1 errors}(a)-(b), respectively. These results suggest that the convergence rates of $\eta$ are quasi-optimal for all three cases. Furthermore, with mesh refinement, the convergence slopes of the error estimators for Cases 1-3 nearly coincide. This demonstrates the superior performance of the adaptive algorithm based on the a posteriori error indicator presented in this work, especially when compared with uniform refinement.

\noindent\textbf{Test 2.} We set the non-homogeneous Navier boundary conditions. Similar to Test 1, we apply both the $C^0$ interior penalty method and the adaptive $C^0$ interior penalty method to this problem. Numerical results on uniform meshes are listed in \Cref{table:exam1_test2}. We perform $P_3$ and $P_4$ polynomial approximations using the adaptive $C^0$ interior penalty method, with the numerical results displayed in Figures \ref{fig:exam1 test2 meshes p3}-\ref{fig:exam1 test2 errors}. The results obtained are similar to those in Test 1. As mentioned earlier, (i) the convergence rates on uniform meshes are $\mathcal{R} \approx 1$; (ii) the position of the Dirac point within the cell has a negligible effect on the convergence rates, whether the meshes are adaptively refined or uniformly refined; (iii) refinements are concentrated around the point $\mathbf{x}_0$; (iv) The convergence rates of $\eta$ are quasi-optimal.

\begin{table}
\centering
\vspace{-3mm}
\caption{Example \ref{exam1} Test 1: Convergence rate of numerical solution in uniform meshes.}
\begin{tabular}{|c|c|c|c|c||c|c|c|c||c|c|c|c|c|c|c|c|}
\hline
&  \multicolumn{3}{c}{$\qquad \quad P_2$} & &\multicolumn{3}{c}{$\qquad  \quad P_3$} & &\multicolumn{3}{c}{$\qquad  \quad P_4$}&\\
\hline
j &    4 &   5 &   6 &7 &    3 &    4 &   5 &   6  &     2 &    3 &    4 &   5\\
\hline
Case 1  & 0.99 & 0.99 & 1.00   & 1.00 & 1.60  & 1.41 & 1.18 & 1.06 & 1.97 & 1.59 & 1.14 & 1.02\\
\hline
Case 2  & 0.98 & 0.99 & 0.99 & 1.00 & 1.67 & 1.72 & 1.20  & 0.96 & 1.72 & 1.59 & 1.35 & 0.87\\
\hline
Case 3  & 0.98 & 0.99 & 0.99 & 1.00 & 1.68 & 1.72 & 1.33 & 1.23 & 1.93 & 1.60  & 1.32 & 0.94\\
\hline
\end{tabular}\label{table:exam1_text1}
\end{table}

\begin{table}
\centering
\vspace{-3mm}
\caption{Example \ref{exam1} Test 2: Convergence rate of numerical solution in uniform meshes.}
\begin{tabular}{|c|c|c|c|c||c|c|c|c||c|c|c|c|c|c|c|c|}
\hline
&  \multicolumn{3}{c}{$\qquad \quad P_2$} & &\multicolumn{3}{c}{$\qquad  \quad P_3$} & &\multicolumn{3}{c}{$\qquad  \quad P_4$}&\\
\hline
j &    4 &   5 &   6 &7 &    3 &    4 &   5 &   6  &     2 &    3 &    4 &   5\\
\hline
Case 1   & 0.99 & 0.99 & 1.00   & 1.00 & 1.60  & 1.41 & 1.18 & 1.06 & 1.97 & 1.59 & 1.14 & 1.02\\
\hline
Case 2  & 0.98 & 1.00   & 0.99 & 1.00 & 1.67 & 1.72 & 1.19 & 0.95 & 1.68 & 1.57 & 1.35 & 0.87\\
\hline
Case 3   & 0.99 & 0.99 & 0.99 & 1.00  & 1.71 & 1.75 & 1.33 & 1.23 & 1.94 & 1.58 & 1.32 & 0.94\\
\hline
\end{tabular}\label{table:exam1_test2}
\end{table}
\begin{figure}
\centering
\subfigure[Case 1,$x_0 = (0,0)$]
{\includegraphics[width=0.32\textwidth]{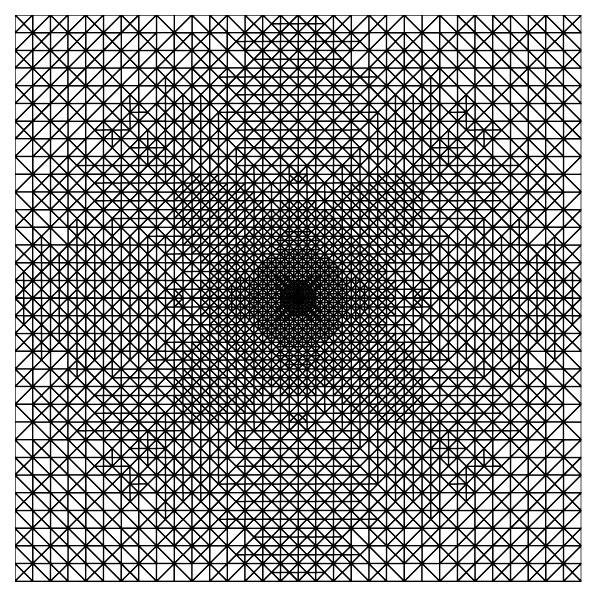}}
\subfigure[Case 2,$x_0 = (-\sqrt{7},-\pi)$]
{\includegraphics[width=0.32\textwidth]{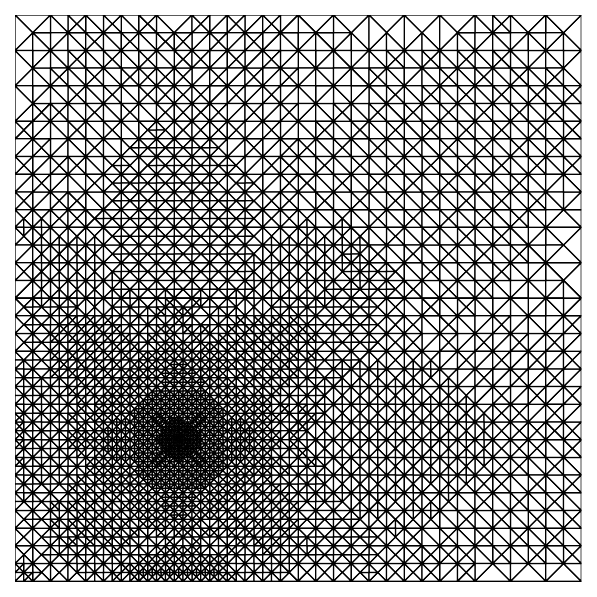}}
\subfigure[Case 3,$x_0 = (\sqrt{5},\sqrt{8})$]
{\includegraphics[width=0.32\textwidth]{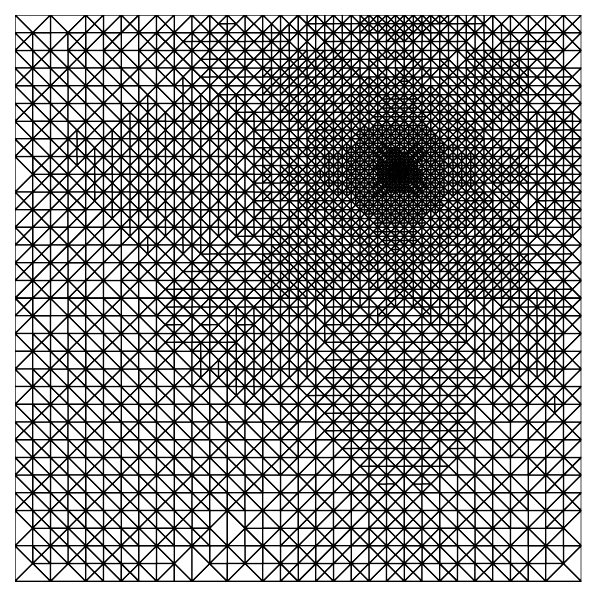}}
\vspace{-2mm}
\caption{Example \ref{exam1} Test 1: adaptive meshes for $P_3$.}\label{fig:exam1 test1 meshes p3}
\end{figure}

\begin{figure}
\centering
\subfigure[Case 1,$x_0 = (0,0)$]
{\includegraphics[width=0.32\textwidth]{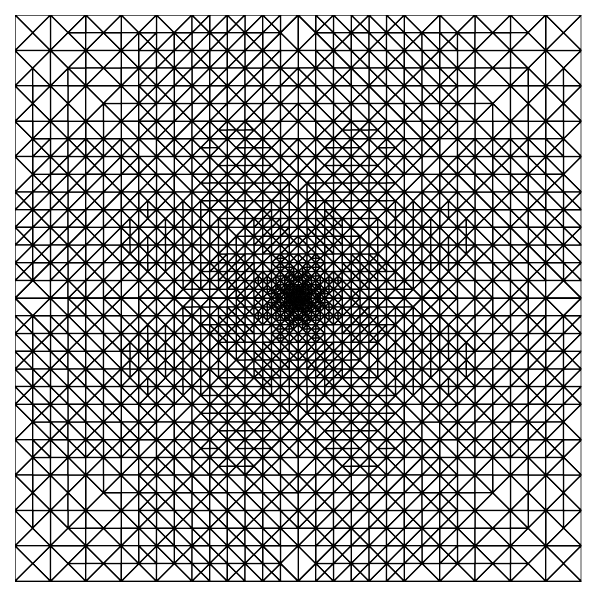}}
\subfigure[Case 2,$x_0 = (-\sqrt{7},-\pi)$]
{\includegraphics[width=0.32\textwidth]{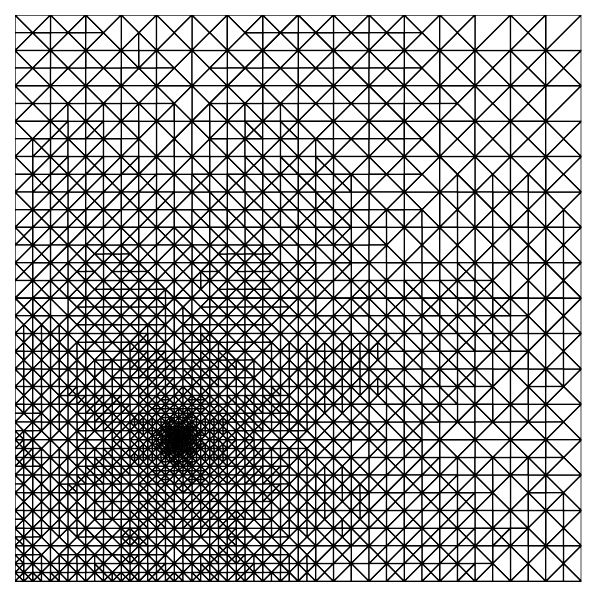}}
\subfigure[Case 3,$x_0 = (\sqrt{5},\sqrt{8})$]
{\includegraphics[width=0.32\textwidth]{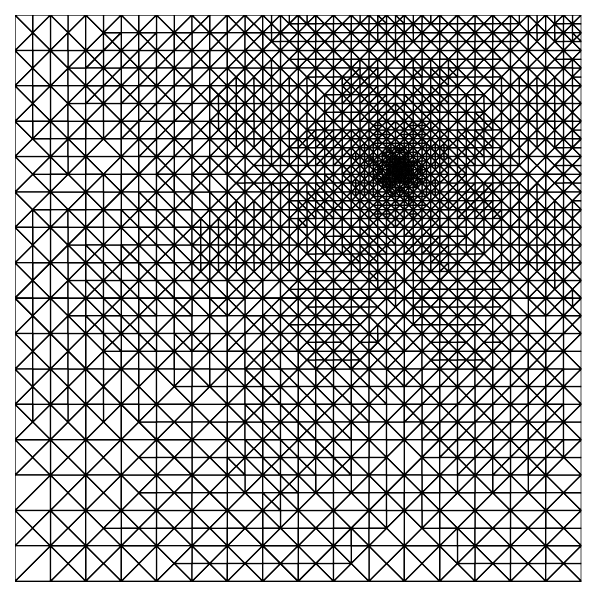}}
\vspace{-2mm}
\caption{Example \ref{exam1} Test 1: adaptive meshes for $P_4$.}\label{fig:exam1 test1 meshes p4}
\end{figure}

\begin{figure}
\centering
\subfigure[Case 1,$x_0 = (0,0)$]
{\includegraphics[width=0.32\textwidth]{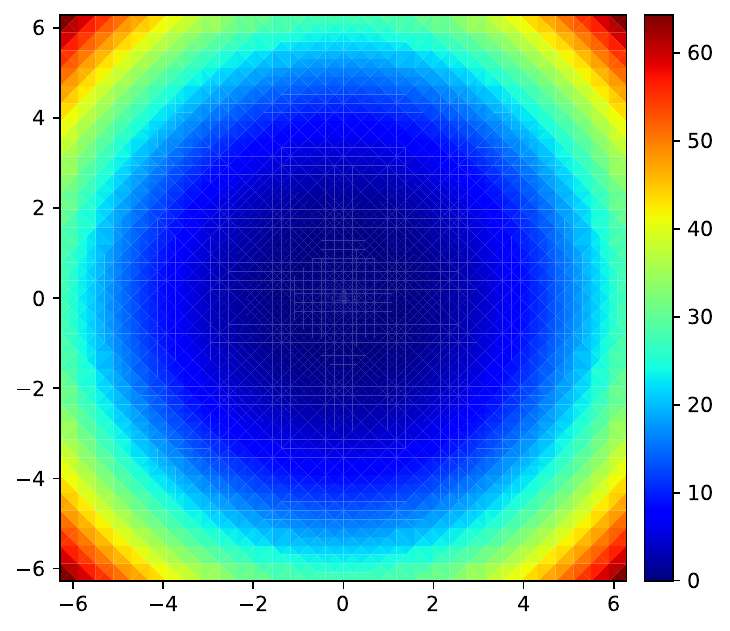}}
\subfigure[Case 2,$x_0 = (-\sqrt{7},-\pi)$]
{\includegraphics[width=0.32\textwidth]{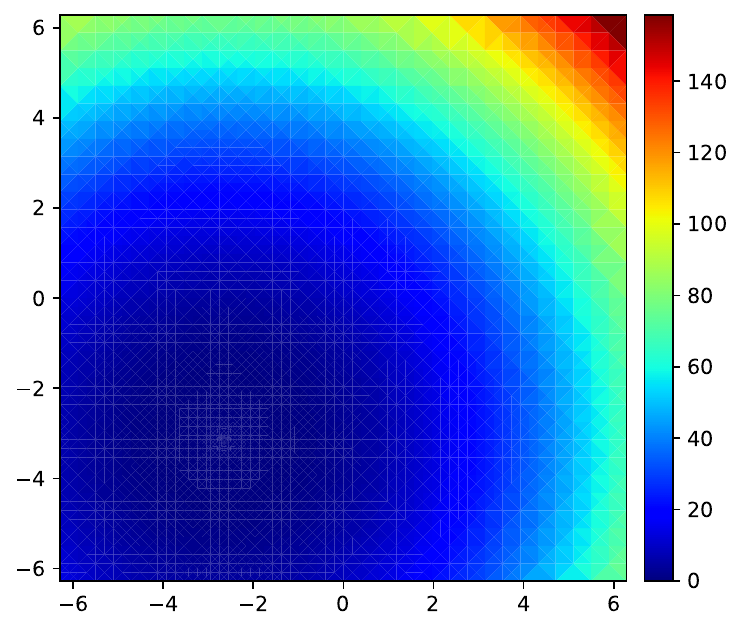}}
\subfigure[Case 3,$x_0 = (\sqrt{5},\sqrt{8})$]
{\includegraphics[width=0.32\textwidth]{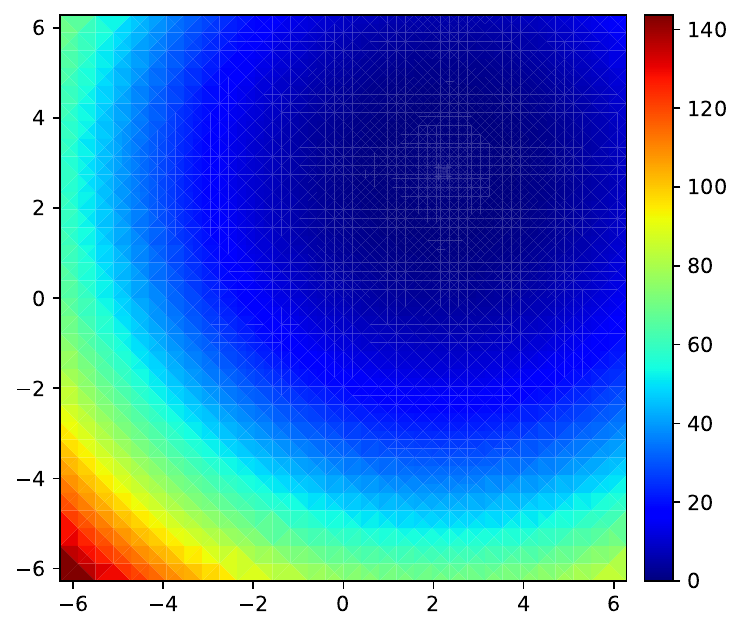}}
\vspace{-2mm}
\caption{Example \ref{exam1} Test 1:adaptive solution for $P_3$.}\label{fig:exam1 test1 solution}
\end{figure}

\begin{figure}
\centering
\subfigure[$P_3$]
{\includegraphics[width=0.42\textwidth]{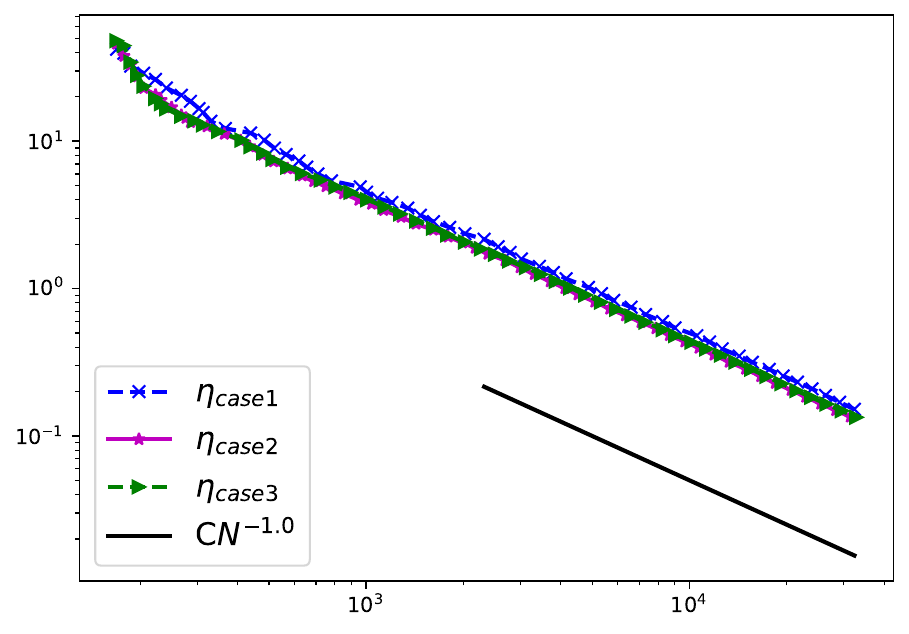}}
\subfigure[$P_4$]
{\includegraphics[width=0.42\textwidth]{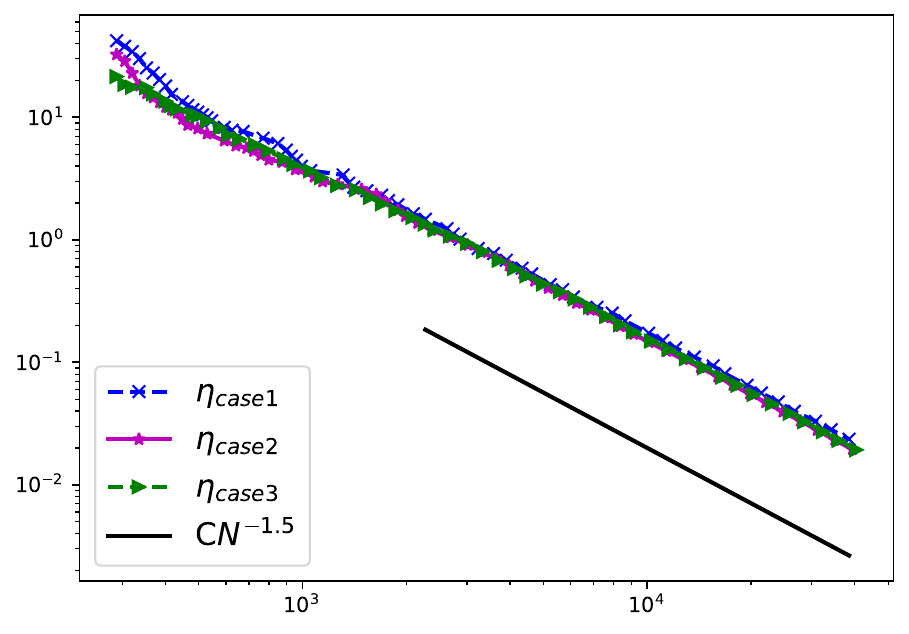}}
\vspace{-2mm}
\caption{Example \ref{exam1} Test 1: error estimators.}\label{fig:exam1 test1 errors}
\end{figure}

\begin{figure}
\centering
\subfigure[Case 1,$x_0 = (0,0)$]
{\includegraphics[width=0.32\textwidth]{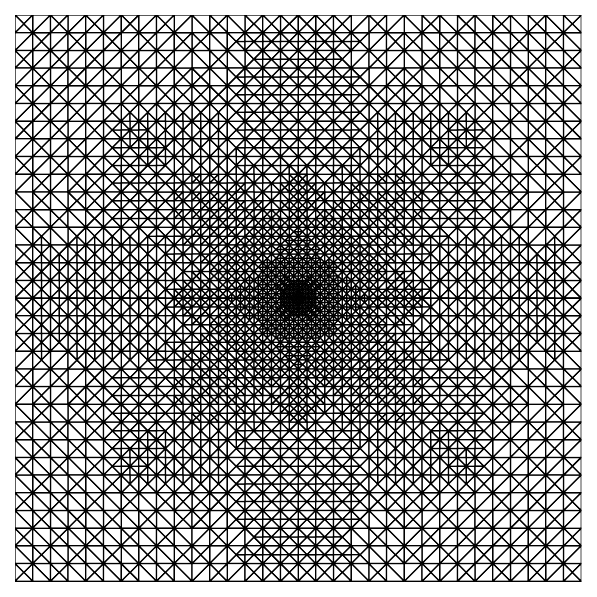}}
\subfigure[Case 2,$x_0 = (-\sqrt{7},-\pi)$]
{\includegraphics[width=0.32\textwidth]{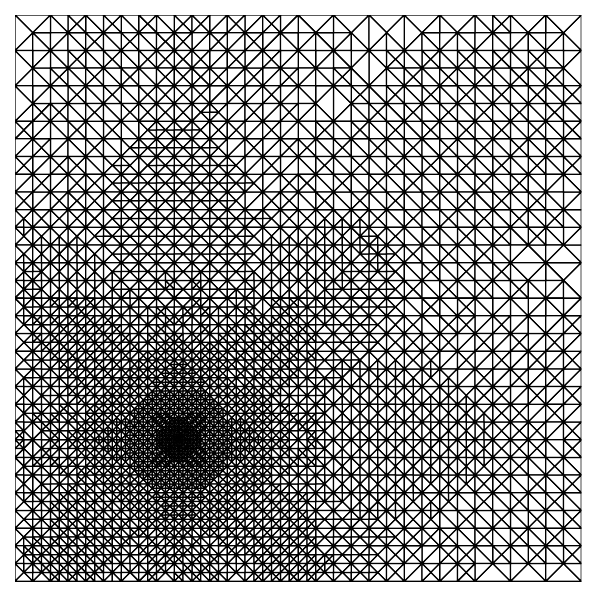}}
\subfigure[Case 3,$x_0 = (\sqrt{5},\sqrt{8})$]
{\includegraphics[width=0.32\textwidth]{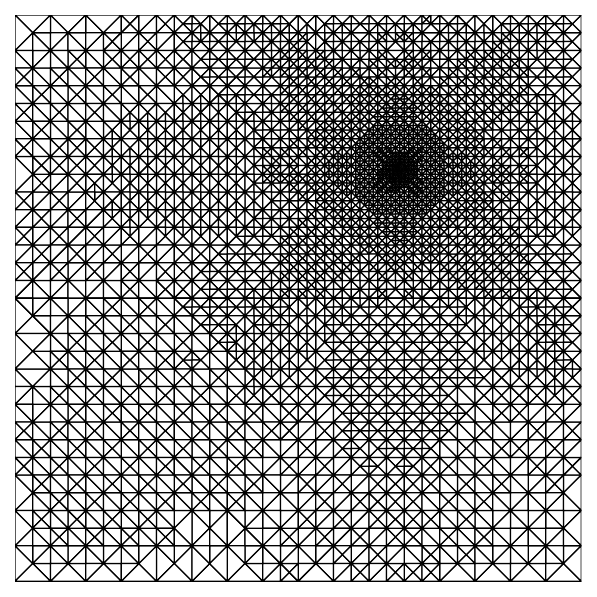}}
\vspace{-2mm}
\caption{Example \ref{exam1} Test 2: adaptive meshes for $P_3$.}\label{fig:exam1 test2 meshes p3}
\end{figure}

\begin{figure}
\centering
\subfigure[Case 1,$x_0 = (0,0)$]
{\includegraphics[width=0.32\textwidth]{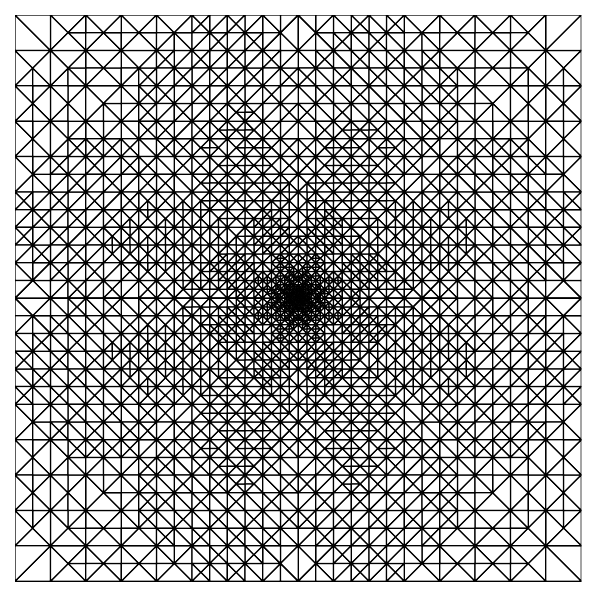}}
\subfigure[Case 2,$x_0 = (-\sqrt{7},-\pi)$]
{\includegraphics[width=0.32\textwidth]{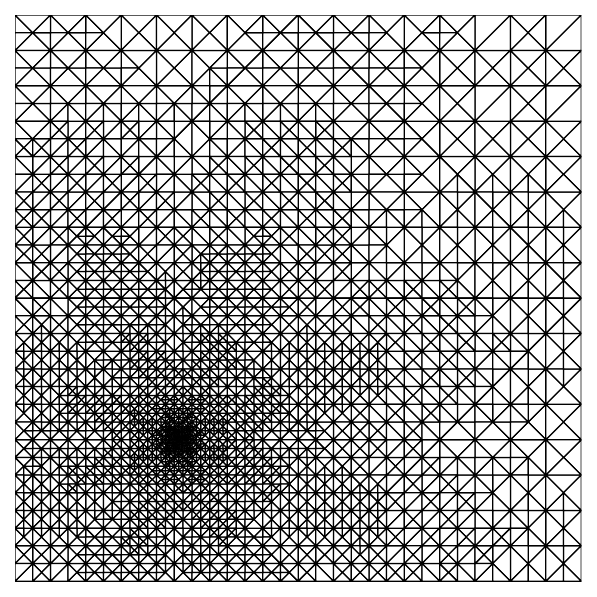}}
\subfigure[Case 3,$x_0 = (\sqrt{5},\sqrt{8})$]
{\includegraphics[width=0.32\textwidth]{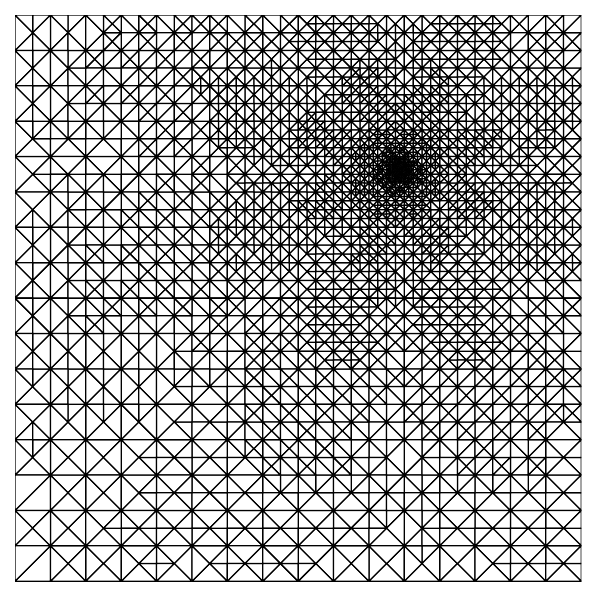}}
\vspace{-2mm}
\caption{Example \ref{exam1} Test 2: adaptive meshes for $P_4$.}\label{fig:exam1 test2 meshes p4}
\end{figure}

\begin{figure}
\centering
\subfigure[Case 1,$x_0 = (0,0)$]
{\includegraphics[width=0.32\textwidth]{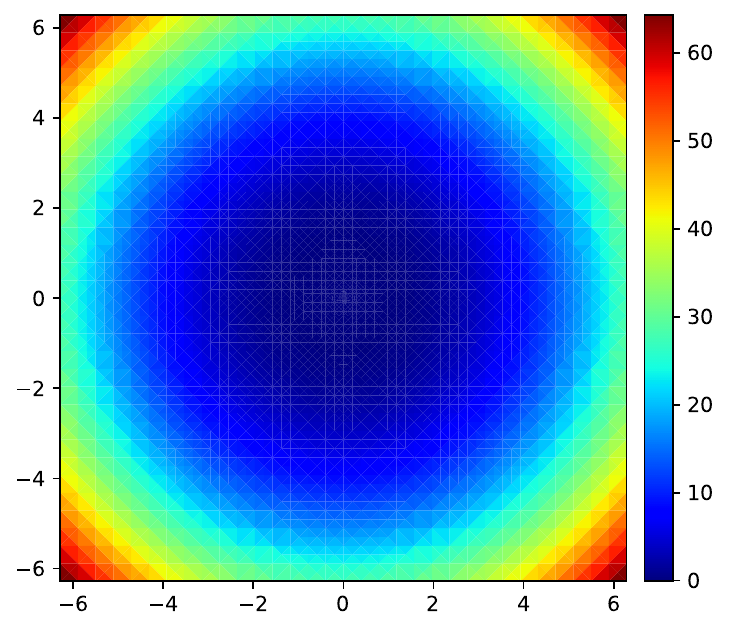}}
\subfigure[Case 2,$x_0 = (-\sqrt{7},-\pi)$]
{\includegraphics[width=0.32\textwidth]{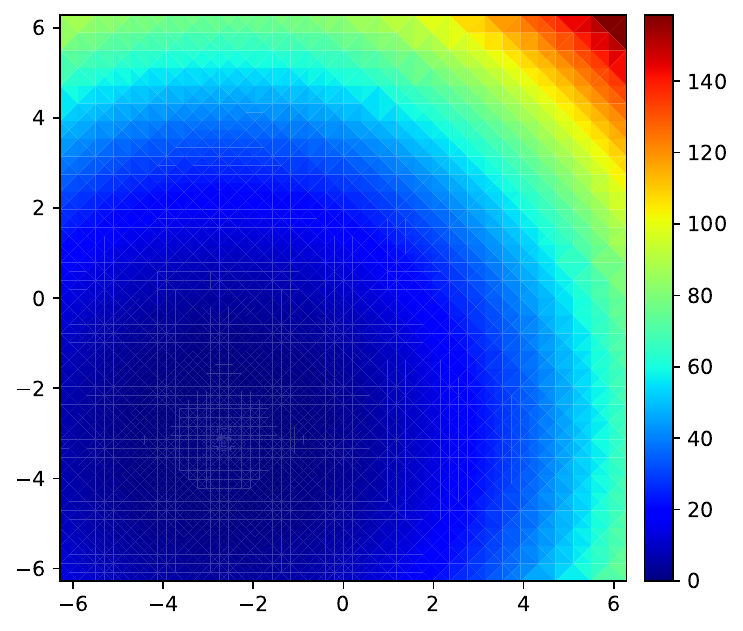}}
\subfigure[Case 3,$x_0 = (\sqrt{5},\sqrt{8})$]
{\includegraphics[width=0.32\textwidth]{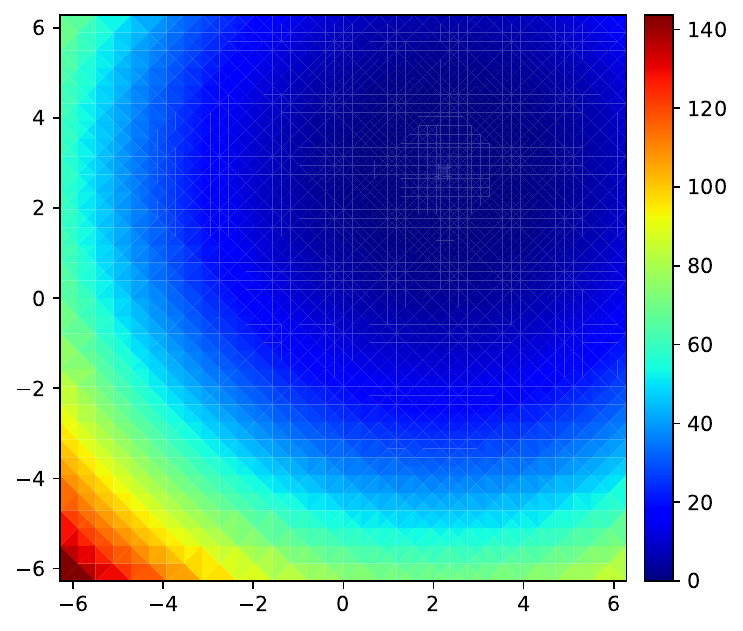}}
\vspace{-2mm}
\caption{Example \ref{exam1} Test 2: adaptive solution for $P_3$.}\label{fig:exam1 test2 solution}
\end{figure}

\begin{figure}
\centering
\subfigure[$P_3$]
{\includegraphics[width=0.42\textwidth]{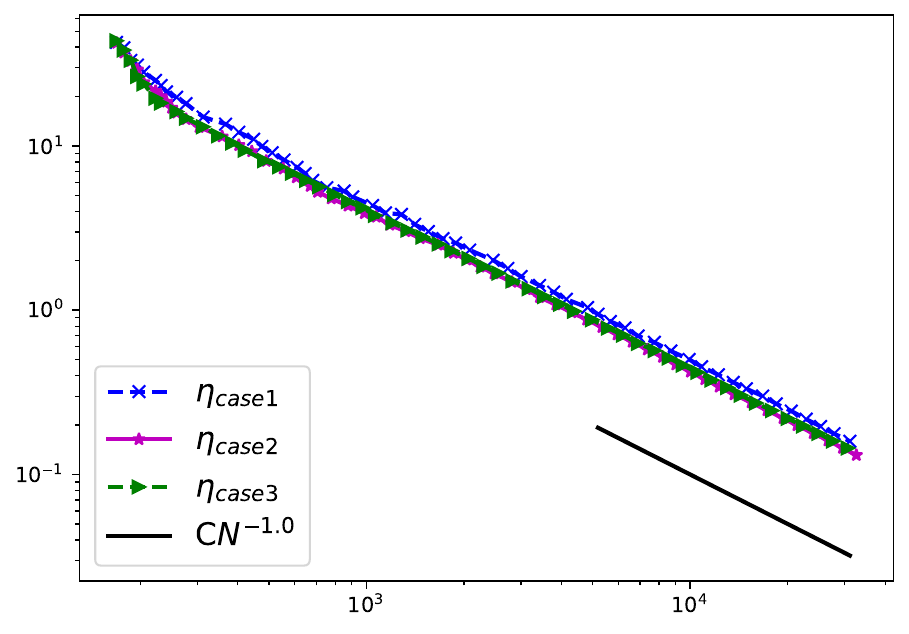}}
\subfigure[$P_4$]
{\includegraphics[width=0.42\textwidth]{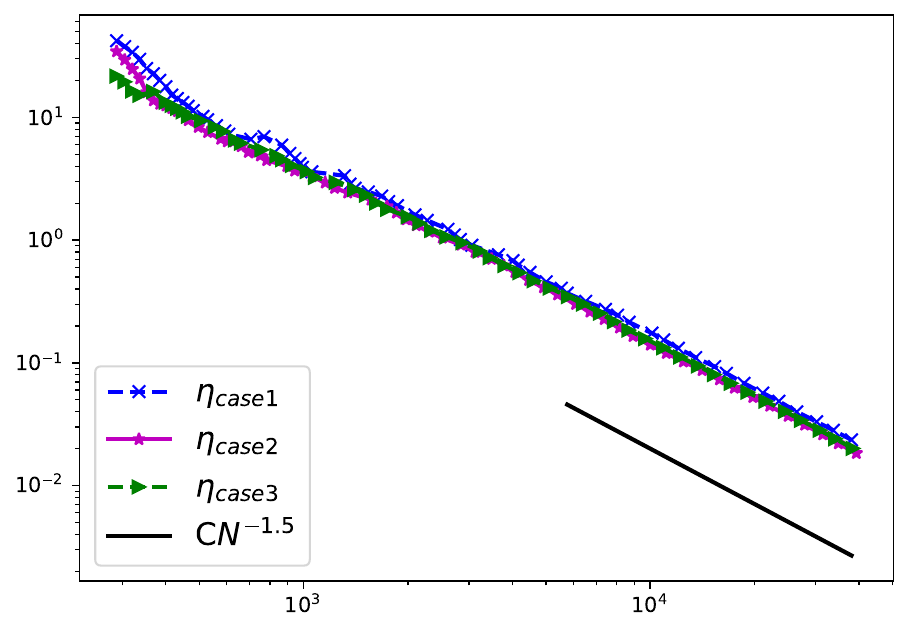}}
\vspace{-2mm}
\caption{Example \ref{exam1} Test 2: error estimators.}\label{fig:exam1 test2 errors}
\end{figure}
\end{example}

\begin{example}[Homogeneous Neumann boundary]\label{exam3}
For the final example, we consider problem \eqref{eq:bh6} with homogeneous Neumann boundary on convex domain $\Omega = (-2\pi,2\pi) \times (-2\pi,2\pi)$. The right-hand side function is given as $\delta_{\mathbf{x}_1} - \delta_{\mathbf{x}_0}$, where $\mathbf{x}_0 = (-\pi,0),\,\mathbf{x}_1 = (\pi,0)$. An initial uniform triangular mesh $\mathcal{T}_h^0$ is shown in Figure \ref{fig:exam3 meshes}(a). 

The convergence rates of the $C^0$ interior penalty method solutions based on $P_m,\,m = 2,3,4$ polynomials on quasi-uniform meshes are shown in Table \ref{table:exam3}. We observe that $\mathcal{R} \approx 1$. Due to the low global regularity of the solution, the convergence rates on quasi-uniform meshes can not reach the optimal order for $P_3,\ P_4$ polynomial approximations. For enhanced accuracy, the adaptive $C^0$ interior penalty method is better suited for this type of problem. The contour of adaptive $C^0$ interior penalty method approximation based on $P_3$ is shown in Figure \ref{fig:exam3 meshes}(a). Figure \ref{fig:exam3 meshes}(b)-(c) show the adaptive meshes of $P_3,\,P_4$ polynomial approximations, respectively. We can see clearly that the error estimator guides the mesh refinement densely around the points $\mathbf{x}_0$ and $\mathbf{x}_1$. The convergence rates of error estimator $\eta$ based on $P_3,\,P_4$ polynomials are shown Figure \ref{fig:exam3 errors}(b)-(c). These results align with expectations, indicating that the convergence rates of the error estimator are quasi-optimal.
\begin{table}
\centering
\vspace{-3mm}
\caption{Example \ref{exam3}: Convergence rate of numerical solution inuniform meshes.}
\begin{tabular}{|c|c|c|c|c||c|c|c|c||c|c|c|c|c|c|c|c|}
\hline
&  \multicolumn{3}{c}{$\qquad \quad P_2$} & &\multicolumn{3}{c}{$\qquad  \quad P_3$} & &\multicolumn{3}{c}{$\qquad  \quad P_4$}&\\
\hline
j &    4 &   5 &   6 &7 &    3 &    4 &   5 &   6  &     2 &    3 &    4 &   5\\
\hline
$\mathcal{R}$&  0.87 &  0.90  &  0.91 &  0.92 &  1.00   &  1.00   &  1.00   &  1.00  &  1.00 &  1.00 &  1.00 &  1.00\\
\hline
\end{tabular}\label{table:exam3}
\end{table}
\begin{figure}
\centering
\subfigure[Initial mesh]
{\includegraphics[width=0.3\textwidth]{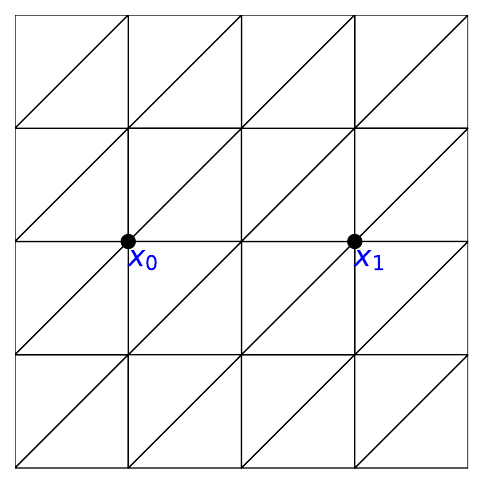}}\hspace{3mm}
\subfigure[Adaptive mesh of $P_3$]
{\includegraphics[width=0.3\textwidth]{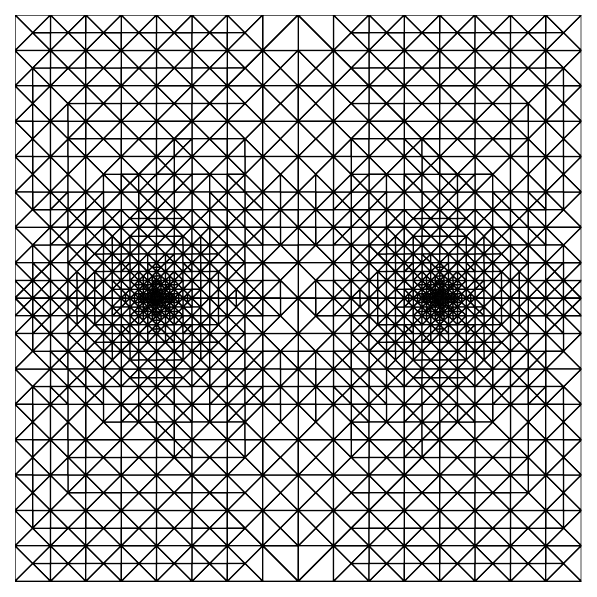}}\hspace{3mm}
\subfigure[Adaptive mesh of $P_4$]
{\includegraphics[width=0.3\textwidth]{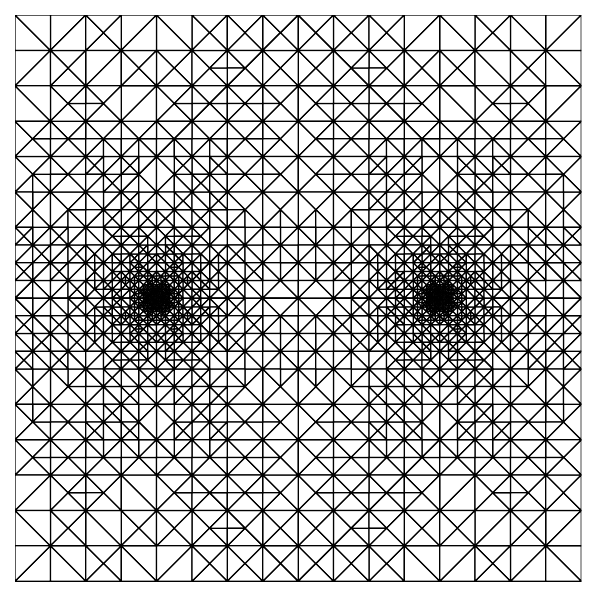}}
\vspace{-2mm}
\caption{Example \ref{exam3}: initial mesh and adaptive meshes.}\label{fig:exam3 meshes}
\end{figure}
\begin{figure}
\centering
\subfigure[Numerical solution]
{\includegraphics[width=0.32\textwidth]{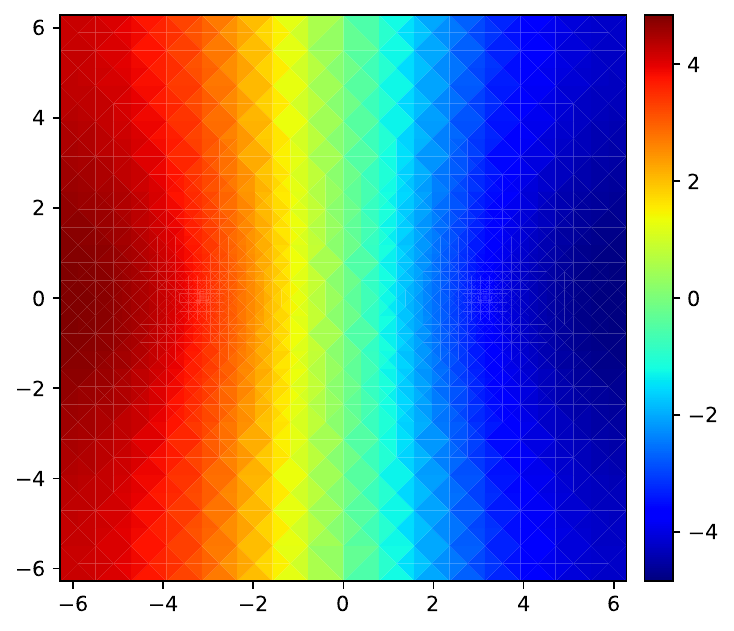}}
\subfigure[Error estimator of $P_3$]
{\includegraphics[width=0.32\textwidth]{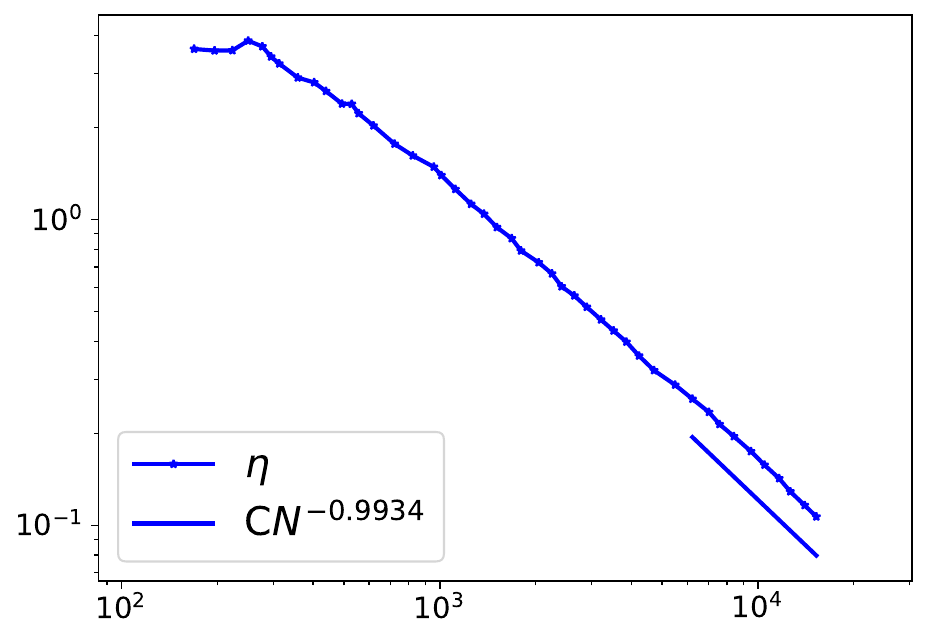}}
\subfigure[Error estimator of $P_4$]
{\includegraphics[width=0.32\textwidth]{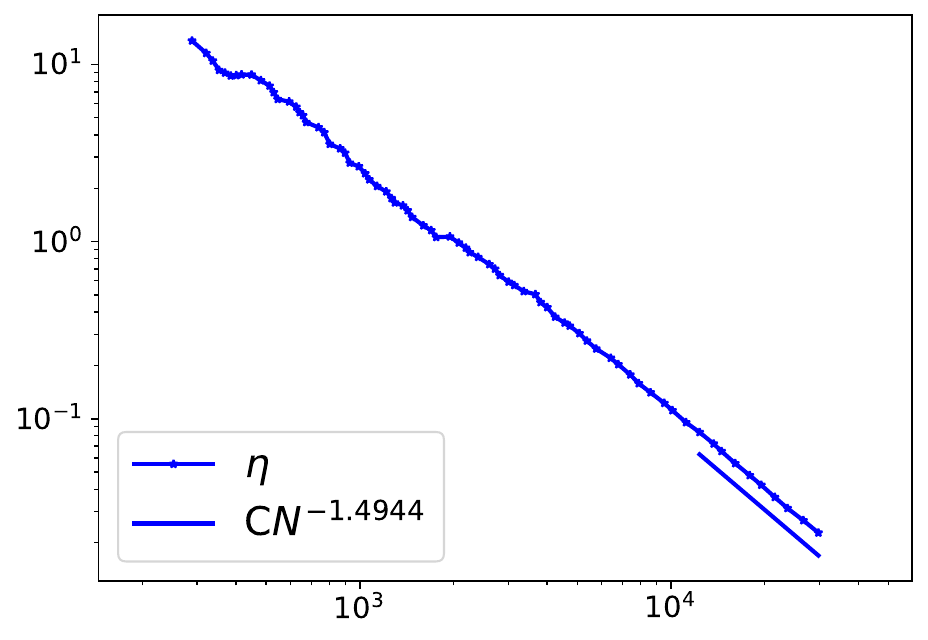}}
\vspace{-2mm}
\caption{Example \ref{exam2}: numerical solution and error estimators.}\label{fig:exam3 errors}
\end{figure}
\end{example}

\section{Conclusion}\label{sec6}
Two residual-based a posteriori estimators are proposed for biharmonic problem \eqref{eq:bh1}. The first estimator is directly derived from the model equation, while the second estimator is based on the projection of the Dirac delta function onto the discrete finite element space. 
The later one introduces an additional projection error, which is included in the error estimator, yielding a similar effect to the first estimator in guiding mesh refinement.
For both a posteriori estimators, we rigorously prove that these estimators are efficient and reliable. An adaptive $C^0$ interior penalty algorithm is provided based on the proposed a posteriori estimators. Extensions of the first estimator to more general fourth order elliptic equations are provided, and quasi-optimal convergence rates are numerically observed.

\section*{Acknowledgments}
Y. Huang was supported in part by NSFC Project (12431014), Project of Scientific Research Fund of Hunan Provincial Science and Technology Department (2020ZYT003).
N. Yi was supported by NSFC Project (12071400), Project of Scientiﬁc Research Fund of the Hunan Provincial Science and Technology Department (2024ZL5017), 
and Program for Science and Technology Innovative Research Team in Higher Educational Institutions of Hunan Province of China.

\section*{Data availability} 
Enquiries about data availability should be directed to the authors.

\section*{Declarations} 
The authors declare that they have no conflict of interest.


\begin{thebibliography}{10}
\bibitem{ABR2006} 
R. Araya, E. Behrens, and  R. Rodr\'{i}guez.
\newblock  A posteriori error estimates for elliptic problems with Dirac delta source terms.
\newblock {\em  Numerische Mathematik}, 105(2):193-216, 2006.

\bibitem{ABR2007} 
R. Araya, E. Behrens, and  R. Rodr\'{i}guez.
\newblock An adaptive stabilized finite element scheme for a water quality model.
\newblock {\em  Computer Methods in Applied Mechanics and Engineering}, 196(29-30):2800-2812, 2007.

\bibitem{A1968}
\newblock J.H. Argyris, I. Fried and D.W. Scharpf. 
\newblock The TUBA family of plate elements for the matrix displacement method.
\newblock {\em The Aeronautical Journal,} 72(692):701-709, 1968.

\bibitem{AGM2014}
J.P. Agnelli, E.M. Garau and P. Morin
\newblock A posteriori error estimates for elliptic problems with Dirac measure terms in weighted spaces. 
\newblock {\em  ESAIM: Mathematical Modelling and Numerical Analysis}, 48(6):1557-1581,2014.

\bibitem{bacuta2002}
C. Bacuta and J.H. Bramble and J.E. Pasciak.
\newblock  Shift theorems for the biharmonic Dirichlet problem.
\newblock {\em  Springer}: 1-26,2002.

\bibitem{bourlard1992}
M. Bourlard and M. Dauge and M.S. Lubuma and S. Nicaise.
\newblock  Coefficients of the singularities for elliptic boundary value problems on domains with conical points. {III}: Finite element methods on polygonal domains.
\newblock {\em  SIAM Journal on Numerical Analysis}, 29(1):136-155, 1992.

\bibitem{BGS2010} 
S.C. Brenner, T. Gudi, and  L.Y. Sung.
\newblock  An a posteriori error estimator for a quadratic $C^0$ -interior penalty method for the biharmonic problem.
\newblock {\em  IMA Journal of Numerical Analysis}, 30: 777–798, 2010.

\bibitem{BS2005}
S.C. Brenner and L.Y. Sung.
\newblock $C^0$ interior penalty methods for fourth order elliptic boundary value problems on polygonal domains.
\newblock {\em Journal of Scientific Computing}, 23(23): 83-118, 2005.

\bibitem{BGGS2012}
S.C. Brenner, S. Gu, T. Gudi and L.Y. Sung.
\newblock  A quadratic $C^0$ interior penalty method for linear fourth order boundary value problems with boundary conditions of the Cahn--Hilliard type.
\newblock {\em SIAM Journal on Numerical Analysis}, 50(4): 2088-2110, 2012.

\bibitem{BN2011}
S.C. Brenner and M. Neilan.
\newblock A $C^0$ interior penalty method for a fourth order elliptic singular perturbation problem.
\newblock {\em SIAM Journal on Numerical Analysis}, 49(2): 869-892, 2011.


\bibitem{BS2008}
S.~Brenner and L.~Scott.
\newblock The mathematical theory of finite element methods.
\newblock {\em Volume 15 of Texts in Applied Mathematics}, 3rd edn. Springer, New York, 2008.

\bibitem{Camp87}
C. V. Camp.
\newblock A solution of the nonhomogeneous biharmonic equation by the boundary element method.
\newblock {\em Ph.D. thesis, Oklahoma State University}, 1987.

\bibitem{Ciarlet74}
Philippe G. Ciarlet.
\newblock The Finite Element Method for Elliptic Problems.
\newblock {\em Universit\'{e} Pierre et Marie Curie, Paris, France,} 1974.

\bibitem{CH1958}
J. W. Cahn and J. E. Hilliard.
\newblock  Free energy of a nonuniform system-I: Interfacial free energy.
\newblock {\em The Journal of Chemical Physics}, 28(2): 258-267, 1958.

\bibitem{CL2011}
J.T. Chen,  H. Z. Liao and W.M. Lee
\newblock  An analytical approach for the Green's functions of biharmonic problems with circular and annular domains.
\newblock {\em Journal of Mechanics}, 25(1): 59-74,2011.


\bibitem{DDPS1979} 
J.J. Douglas, T. Dupont, P. Percell and R. Scott.
\newblock  A family of $C^1$ finite elements with optimal approximation properties for various Galerkin methods for 2nd and 4th order problems.
\newblock {\em RAIRO Analyse $num\acute{e}rique$}, 13(3):227–255, 1979.

\bibitem{EGHLMT2002} 
G. Engel, K. Garikipati, T.J.R. Hughes, M.G. Larson, L. Mazzei and R.L. Taylor.
\newblock  Continuous/discontinuous ﬁnite element approximations of fourth-order elliptic problems in structural and continuum mechanics with applications to thin beams and plates,
and strain gradient elasticity.
\newblock {\em Computer methods in applied mechanics and engineering}, 191:3669-3750, 2002.

\bibitem{Grisvard92}
P.~Grisvard.
\newblock {\em Singularities in Boundary Value Problems}, volume~22 of {\em
  Research Notes in Applied Mathematics}.
\newblock Springer-Verlag, New York, 1992.

\bibitem{GHZ2014}
W. Gong, M. Hinze and Z.J. Zhou.
\newblock A priori error analysis for the finite element approximations of parabolic optimal control problems with
pointwise control.
\newblock {\em SIAM Journal on Control and Optimization,} 52(1):97-119, 2014.

\bibitem{GHV2011}
E.H. Georgoulis, P. Houston and J. Virtanen.
\newblock An aposteriori error indicator for discontinuous Galerkin
approximations of fourth-order elliptic problems.
\newblock {\em IMA Journal of Numerical Analysis,} 31:281-298, 2011.

\bibitem{GMV2016}
F.D. Gaspoz, P.Morin and A, Veeser.
\newblock A posteriori error estimates with point sources in fractional sobolev spaces .
\newblock {\em Numerical Methods for Partial Differential Equation}, 33(4):1018-1042,2016.

\bibitem{GPJ2009} 
E.H. Georgoulis, P. Houston and J. Virtanen.
\newblock  An a posteriori error indicator for discontinuous Galerkin approximations of fourth-order elliptic problems.
\newblock {\em IMA Journal of Numerical Analysis}, 31(1): 281–298, 2009.

\bibitem{HW2012}
P. Houston and T. Wihler.
\newblock  Discontinuous Galerkin methods for problems with Dirac delta source.
\newblock {\em  ESAIM: Mathematical Modelling and Numerical Analysis}, 46(6):1467-1483, 2012.

\bibitem{J1975}
J. D. Jackson.
\newblock  Classical electrodynamics.
\newblock {\em  John Wiley and Sons, Inc., New York-London-Sydney}, second edition, 1975.

\bibitem{kozlov2001}
V.A. Kozlov, V.G. Maz$\acute{y}$a and J. Rossmann.
\newblock  Spectral problems associated with corner singularities of solutions to elliptic equations.
\newblock {\em  American Mathematical Society}, 85,2001.

\bibitem{KN2014}
S. B. G. Karakoc and M. Neilan.
\newblock  A $C^0$ finite element method for the Biharmonic problem without extrinsic penalization.
\newblock {\em  Numerical Methods for Partial Differential Equations}, 30(4),1254-1278,2014.

\bibitem{L2020}
D. Leykekhman.
\newblock  Pointwise error estimates for $C^0$ interior penalty approximation of biharmonic problems.
\newblock {\em Mathematics of Computation}, 90(327):41-63, 2020.

\bibitem{LV2013}
D. Leykekhman and B. Vexler.
\newblock  Optimal a priori error estimates of parabolic optimal control problems with pointwise control.
\newblock {\em SIAM Journal on Numerical Analysis}, 51(5): 2797-2821, 2013.

\bibitem{LWY22}
H. Li, C. D. Wickramasinghe and P. Yin.
\newblock  A $C^0$ finite element method for the biharmonic problem with Dirichlet boundary conditions in a polygonal domain.
\newblock {\em arXiv preprint arXiv:2207.03838}, 2022.


\bibitem{M20}
W. McLean.
\newblock {\em Strongly Elliptic Systems and Boundary Integral Equations}.
\newblock Cambridge University Press, 2000.

\bibitem{MMRZ2022}
F. Millar, I.Muga, S Rojas and K.G. Van der Zee.
\newblock  Projection in negative norms
and the regularization of rough linear functionals.
\newblock {\em Numerische Mathematik}, 150: 1087-1121, 2022.

\bibitem{RH2012} 
R. An and X.H. Huang.
\newblock  Constrained $C^0$ finite element methods for biharmonic problem.
\newblock {\em  Abstract and Applied Analysis}, 2012(137), 2012.


\bibitem{Scott73}
R. Scott.
\newblock Finite element convergence for singular data.
\newblock {\em Numerische Mathematik,} 21:317--327, 1973.

\bibitem{S1992}
B. Semper.
\newblock Conforming finite element approximations for a fourth-order singular perturbation problem.
\newblock {\em SIAM Journal on Numerical Analysis}, 29(4): 1043-1058, 1992.

\bibitem{SKF1999}
J.Y. Shu, W.E. King and N.A. Fleck.
\newblock Finite elements for materials with strain gradient effects.
\newblock {\em International Journal for Numerical Methods in Engineering}, 44(3): 373-391, 1999.




\bibitem{TSW1959}
S. Timoshenko, S. Woinowsky-krieger.
\newblock Theory of plates and shells.
\newblock {\em McGraw-hill New York}, 1959.

\bibitem{TZ2013}
F. Tan and Y. L. Zhang.
\newblock The regular hybrid boundary node method in the bending analysis of thin plate structures subjected to a concentrated load.
\newblock {\em European Journal of Mechanics A/Solids}, 38:79-89, 2013.

\end{thebibliography}
\end{document}